\documentclass[11pt]{amsart}

\usepackage{amssymb,amscd,verbatim}
\usepackage[all]{xy}

\theoremstyle{plain}
\newtheorem{prop}{Proposition}
\newtheorem{thm}[prop]{Theorem}
\newtheorem{lem}[prop]{Lemma}
\newtheorem{cor}[prop]{Corollary}
\theoremstyle{definition}
\newtheorem{defprop}[prop]{Definition-Proposition}
\newtheorem{rem}[prop]{Remark}
\newtheorem{ex}[prop]{Example}
\newtheorem{definition}[prop]{Definition}

\numberwithin{prop}{section}
\numberwithin{equation}{section}

\newcommand{\Integer}{\mathbb{Z}}
\newcommand{\Int}{\Integer_{\ge 0}}
\newcommand{\la}{\lambda}

\newcommand{\C}{\mathbb{C}}
\newcommand{\CC}{\mathcal{C}}
\newcommand{\ST}{\mathrm{ST}}
\newcommand{\CST}{\mathrm{CST}}
\newcommand{\Hom}{\mathrm{Hom}}
\newcommand{\ev}{\mathrm{ev}}

\newcommand{\at}{\widetilde{\alpha}}

\newcommand{\charge}{c}
\newcommand{\co}{\mathrm{co}}
\newcommand{\cc}{cc}

\newcommand{\shape}{\mathrm{shape}}
\newcommand{\word}{\mathrm{word}}
\newcommand{\phib}{\overline{\phi}}
\newcommand{\phit}{\widetilde{\phi}}

\newcommand{\LRT}{\mathrm{CLR}}
\newcommand{\SLRT}{\mathrm{SLR}}
\newcommand{\RLRT}{\mathrm{RLR}}
\newcommand{\CLRT}{\mathrm{LRT}}
\newcommand{\CO}{\mathrm{C}}
\newcommand{\RC}{\mathrm{RC}}
\newcommand{\LRTh}{\LRT\htt}
\newcommand{\RCh}{\RC\htt}

\newcommand{\Kn}{\cong_K}
\newcommand{\ck}{{}^{\vee}}
\newcommand{\htt}{{}^{\wedge}}

\newcommand{\db}{\overline{\delta}}
\newcommand{\dt}{\widetilde{\delta}}

\newcommand{\lt}{\widetilde{\ell}}
\newcommand{\ltt}{\widetilde{\lt}}
\newcommand{\lb}{\overline{\ell}}
\newcommand{\ltb}{\overline{\lt}}
\newcommand{\lbt}{\widetilde{\lb}}
\newcommand{\lbi}{\overline{s}}
\newcommand{\nut}{\widetilde{\nu}}
\newcommand{\nub}{\overline{\nu}}
\newcommand{\nutb}{\overline{\nut}}
\newcommand{\nubt}{\widetilde{\nub}}
\newcommand{\Jt}{\widetilde{\mathrm{J}}}
\newcommand{\Jb}{\overline{\mathrm{J}}}
\newcommand{\Jtb}{\overline{\Jt}}
\newcommand{\Jbt}{\widetilde{\Jb}}
\newcommand{\rk}{\mathrm{rk}}
\newcommand{\rkb}{\overline{\rk}}
\newcommand{\rkt}{\widetilde{\rk}}
\newcommand{\rku}{\underline{\rk}}

\newcommand{\rb}{\overline{r}}
\newcommand{\rt}{\widetilde{r}}
\newcommand{\rbt}{\widetilde{\rb}}
\newcommand{\rtb}{\overline{\rt}}

\newcommand{\cb}{\overline{c}}

\newcommand{\Rht}{R\htt}
\newcommand{\Rck}{R\ck}
\newcommand{\Rhtck}{{\Rht\ck}}

\newcommand{\ii}{\imath}
\newcommand{\il}{\ii^<}
\newcommand{\ig}{\ii^>}
\newcommand{\ip}{\ii^+}
\newcommand{\iht}{{\imath\htt}}
\newcommand{\ick}{{\imath\ck}}

\newcommand{\jj}{\jmath}
\newcommand{\jl}{\jj^<}
\newcommand{\jg}{\jj^>}
\newcommand{\jp}{\jj^+}
\newcommand{\jht}{{\jmath\htt}}
\newcommand{\jck}{{\jmath\ck}}

\newcommand{\ib}{-}
\newcommand{\iit}{D}
\newcommand{\jb}{\db}
\newcommand{\jbi}{\jb^{-1}}
\newcommand{\jt}{\dt}

\newcommand{\Rb}{\overline{R}}
\newcommand{\Rt}{\widetilde{R}}
\newcommand{\Rbt}{\widetilde{\Rb}}
\newcommand{\Rtb}{\overline{\Rt}}

\newcommand{\com}{\theta}
\newcommand{\comev}{\theta^{\ev}}

\newcommand{\Tb}{\overline{T}}	

\newcommand{\CW}{\mathrm{CW}}
\newcommand{\RW}{\mathrm{RW}}
\newcommand{\nuht}{{\nu\htt}}
\newcommand{\Jht}{{\mathrm{J}\htt}}
\newcommand{\nuck}{{\nu\ck}}
\newcommand{\Jck}{{\mathrm{J}\ck}}

\newcommand{\lcb}{{\overline{\lc}}}
\newcommand{\lc}{{\ell\ck}}

\newcommand{\Rp}{R^+}
\newcommand{\Rpt}{R^{+t}}

\newcommand{\Rl}{{R^<}}

\newcommand{\Rbl}{{\Rb^<}}

\newcommand{\Rg}{{R^>}}

\newcommand{\ZC}{\mathrm{ZC}}

\newcommand{\ZR}{\mathrm{ZR}}

\newcommand{\tr}{\mathrm{tr}}
\newcommand{\LRtr}{\mathrm{tr}_{\mathrm{LR}}}
\newcommand{\RCtr}{\mathrm{tr}_{\mathrm{RC}}}
\newcommand{\std}{\mathrm{std}}

\newcommand{\bm}{{b^-}}
\newcommand{\beq}{{b^=}}
\newcommand{\bp}{{b^+}}

\newcommand{\dom}{\trianglerighteq}

\newcommand{\Rhtb}{\overline{\Rht}}
\newcommand{\del}{\partial}

\newcommand{\bull}{\bullet}

\newcommand{\mt}{m^t}
\newcommand{\mht}{m\htt}
\newcommand{\nutrb}{\overline{\nu^t}}
\newcommand{\Jtrb}{\overline{J^t}}

\newcommand{\asc}{\mathrm{asc}}
\newcommand{\rows}{\mathrm{rows}}
\newcommand{\Ev}{\mathrm{Ev}}
\newcommand{\WR}{\mathrm{WR}}
\newcommand{\WC}{\mathrm{WC}}

\newcommand{\qbin}[2]{\genfrac{[}{]}{0pt}{}{#1}{#2}}
\newcommand{\qbins}[2]{{\textstyle\genfrac{[}{]}{0pt}{}{#1}{#2}}}

\begin{document}

\title{A bijection between Littlewood-Richardson \\ 
tableaux and rigged configurations}

\author{Anatol N. Kirillov}
\address{Division of Mathematics, Graduate School of Science\\
Hokkaido University\\
Sapporo, 060-0810\\
Japan}
\email{kirillov@math.sci.hokudai.ac.jp}
\author{Anne Schilling${}^*$}
\address{Instituut voor Theoretische Fysica\\
Universiteit van Amsterdam\\
Valckenierstraat 65\\
1018 XE Amsterdam\\
The Netherlands}
\email{schillin@wins.uva.nl}
\thanks{${}^*$ Supported by the 
``Stichting Fundamenteel Onderzoek der Materie''.}
\author{Mark Shimozono${}^\dag$}
\address{Department of Mathematics\\
Virginia Tech\\
Blacksburg, VA 24061-0123\\
U.S.A.}
\email{mshimo@math.vt.edu}
\thanks{${}^{\dag}$ Partially supported by NSF grant DMS-9800941.}

\date{December 1998}

\begin{abstract}
We define a bijection from Littlewood--Richardson tableaux to
rigged configurations and show that it preserves the appropriate 
statistics. This proves in particular a quasi-particle expression 
for the generalized Kostka polynomials $K_{\la R}(q)$ labeled by 
a partition $\la$ and a sequence of rectangles $R$. The generalized 
Kostka polynomials are $q$-analogues of multiplicities of the 
irreducible $GL(n,\C)$-module $V^\la$ of highest weight $\la$
in the tensor product $V^{R_1} \otimes \dots\otimes V^{R_L}$.
\end{abstract}
 
\maketitle

\markboth{A. N. KIRILLOV, A. SCHILLING, AND M. SHIMOZONO}{BIJECTION 
BETWEEN LR TABLEAUX AND RIGGED CONFIGURATIONS}

\section{Introduction}

In their study of the $XXX$ model using Bethe Ansatz techniques, 
Kirillov and Reshetikhin~\cite{KR} obtained an expression for the 
Kostka polynomials as the generating function of rigged configurations. 
Rigged configurations index the solutions of the Bethe Ansatz equations;
they are sequences of partitions obeying certain conditions together 
with quantum numbers or riggings labeling the parts of the partitions. 
In ref.~\cite{KR} a statistic-preserving bijection between column-strict 
Young tableaux and rigged configurations was given
which connects the generating function of rigged configurations to
the charge representation of the Kostka polynomials as obtained by 
Lascoux and Sch\"utzenberger~\cite{LS Kostka}.

In fact, the representation of the Kostka polynomials in terms of 
rigged configurations is exactly in quasi-particle form.
In recent years, much research has been devoted to the study of
quasi-particle representations of characters of conformal field
theories and configuration sums of exactly solvable lattice 
models (see for example \cite{Be,BLS,BMS,DF,FLW,HKKOTY,HKOTY,SW,Waa,Wab} 
and references therein).
These quasi-particle representations are physically interesting
\cite{KKMM1,KKMM2} because they reflect the particle structure of 
the underlying model.

Recently certain generalizations of the Kostka polynomials were
introduced and studied~\cite{KS,SW,S,S2,S3,ShW}. 
These generalized Kostka polynomials $K_{\la R}(q)$ are labeled by a 
partition $\la=(\la_1,\la_2,\ldots)$ and a sequence of rectangles
$R=(R_1,\ldots,R_L)$, that is, each $R_i=(\eta_i^{\mu_i})$ is a partition 
of rectangular shape. They are $q$-analogues of multiplicities of the 
irreducible $GL(n,\C)$-module $V^\la$ of highest weight $\la$
in the tensor product $V^{R_1} \otimes \dots\otimes V^{R_L}$.
This multiplicity is equal to the cardinality
of the set of Littlewood--Richardson tableaux \cite{FG}.
When all $R_i$ are single rows (in which case $R_i=(\eta_i)$), the
generalized Kostka polynomial $K_{\la R}(q)$ reduces to the Kostka 
polynomial $K_{\la\eta}(q)$. 
Conjecturally, the generalized Kostka polynomials coincide with
special cases of the spin generating functions over
ribbon tableaux of Lascoux, Leclerc and Thibon~\cite{LLT}.

In refs.~\cite{SW,S} the generalized Kostka polynomials were expressed 
as the generating function of Littlewood--Richardson tableaux with 
a generalized charge statistic, extending the result for the Kostka 
polynomials of Lascoux and Sch\"utzenberger~\cite{LS Kostka}.
A representation of the generalized Kostka polynomials in terms 
of rigged configurations was conjectured in~\cite{KS,SW}.

In this paper it is shown that the algorithm described 
in~\cite{K} gives a statistic-preserving  bijection
between Littlewood--Richardson tableaux and rigged configurations.
This bijection extends the bijection between column-strict Young tableaux
and rigged configurations~\cite{KR}, and in particular provides a proof
of the quasi-particle representation of the generalized Kostka 
polynomials as conjectured in~\cite{KS,SW}.

This paper is organized as follows. In Section~\ref{sec QP} the set of
Littlewood--Richardson (LR) tableaux and the set of rigged configurations
corresponding to a sequence of rectangles are defined. The charge 
expression of the generalized Kostka polynomials is recalled and
the quasi-particle representation is stated.
In Section~\ref{sec operations} several operations on rectangles
and their analogues on LR tableaux and rigged configurations are discussed.
These operations are crucial for the definition of the bijection
between LR tableaux and rigged configurations as given in 
Definition-Proposition~\ref{phi def}. The definition of the bijection 
is based on the operations of splitting off the last column of the last 
rectangle in $R$ and, if the last rectangle is a single column, removing 
one box from it. The algorithm of ref.~\cite{K} 
for the bijection, which is computationally simpler but less suitable
for the proofs, is stated in Section~\ref{sec algo}.
The Evacuation Theorem~\ref{ev} is proven in Section~\ref{sec ev}.
It states that under the bijection
the evacuation of LR tableaux corresponds to the complementation
of quantum numbers on rigged configurations.
In Section~\ref{sec row rec} another recurrence for $\phib_R$
is given based on operations involving rows instead of columns as
used in Definition-Proposition~\ref{phi def}.
This formulation of $\phib_R$ is in a sense transpose to the 
description of $\phib_{R^t}$ using the column operations where
$R^t=(R_1^t,\ldots,R_L^t)$ and $R_i^t$ is the transpose of $R_i$.
It is used in Section~\ref{sec trans} to prove the Transpose 
Theorem~\ref{transpose thm}. Like the LR coefficients,
the generalized Kostka polynomials have a transpose symmetry which
has been explained combinatorially by a transpose bijection on LR 
tableaux~\cite{SW,S2} and by a transpose bijection on
rigged configurations~\cite{KS}. The Transpose Theorem shows that
the bijection from LR tableaux to rigged configurations intertwines 
these two transpose bijections.
In refs.~\cite{SW,S2} families of statistic-preserving embeddings 
between sets of LR tableaux were defined. The Embedding 
Theorem~\ref{embed bij} of Section~\ref{sec embed} shows that these
embeddings on LR tableaux correspond to an inclusion on rigged 
configurations under the bijection between LR tableaux and rigged 
configurations. The proof of the Embedding Theorem relies on the 
Evacuation Theorem.
In Section~\ref{sec statistics} it is finally shown that the bijection
is statistic preserving. The proof uses the Transpose 
Theorem~\ref{transpose thm} and the Embedding Theorem~\ref{embed bij}
to reduce to the case that all rectangles in $R$ are single boxes.
In the single box case the property that the bijection is statistic
preserving can be checked explicitly.
Some technical proofs are relegated to Appendices A and B.

\section{Charge and quasi-particle representation of the generalized 
Kostka polynomials}
\label{sec QP}

In this section we define the set of Littlewood--Richardson tableaux
and the set of rigged configurations, and recall the charge 
representation of the generalized Kostka polynomials~\eqref{gK}.
The quasi-particle representation of the generalized Kostka polynomials, 
to be proven in this paper, is stated in Theorem~\ref{thm QP}.

\subsection{Littlewood--Richardson tableaux}
\label{sec LR}

Given a partition $\la$ and a sequence of partitions $R=(R_1,\dots,R_L)$,
define the tensor product multiplicity
\begin{equation*}
  c^R_\la = \dim \Hom_{GL(n)}(V^\la, V^{R_1} \otimes \dots\otimes V^{R_L})
\end{equation*}
where $V^\la$ is the irreducible $GL(n,\C)$-module of highest weight $\la$.
There are many well-known equivalent formulations of the
celebrated Littlewood-Richardson rule, which give
the multiplicity $c^R_\la$ as the cardinality of a certain
set of tableaux (see \cite{FG}). We refer to
sets of tableaux that have this cardinality as LR tableaux.
Here we recall various well-known notions of LR tableaux,
which count the multiplicity $c^R_\la$ when $R$ is a sequence of rectangles.
One of these is particularly well-suited for use with 
the bijection to rigged configurations.

First a few tableau definitions are necessary. The English convention
is used here for partition diagrams and tableaux.
For a partition $\la$ denote by $\ST(\la)$ the set of 
standard tableaux of shape $\la$.
Given a possibly skew column-strict tableau $T$ in an alphabet $A$
and $B$ a subinterval of $A$, denote by $T|_B$ the restriction of
$T$ to $B$, which is by definition the skew
column-strict tableau given by erasing from $T$ the letters that are
not in $B$.  Define the row-reading word of $T$ to be the concatenation
$\word(T) = \dots v^2 v^1$, where $v^r$ is the word given
by the $r$-th row of $T$ read from left to right.  Knuth \cite{Kn}
defined an equivalence relation on words denoted $v \Kn w$.
Given a word $w$ there is a unique column-strict tableau
$P(w)$ such that $\word(P(w)) \Kn w$; this is the Schensted $P$-tableau
of the word $w$.  Write $P(T):=P(\word(T))$ for the skew column strict
tableau $T$.

Let $R=(R_1,\dots,R_L)$ be a sequence of partitions.
For $1\le j \le L$ let $B_j$ be an interval of integers such that
if $i<j$, $x\in B_i$ and $y\in B_j$, then $x<y$.
Set $B = \bigcup_{j=1}^L B_j$.
Let $Z_j$ be any column-strict tableau of shape $R_j$ in the
alphabet $B_j$ for $1\le j \le L$.
Define the set $\SLRT(\la;Z)$ to be the set of column-strict tableaux $T$
of shape $\la$ in the alphabet $B$ such that $P(T|_{B_j})=Z_j$ for all $j$.
It is well-known that $|\SLRT(\la;Z)|=c^R_\la$.

\begin{ex} \label{KS LR}
Suppose $R_j$ has $\mu_j$ parts for all $j$.
Define $B_j=[\mu_1+\dots+\mu_{j-1}+1,\mu_1+\dots+\mu_{j-1}+\mu_j]$.
Let $Z_j$ be the column-strict tableau of shape $R_j$
whose $r$-th row is filled with copies of the $r$-th largest
letter of $B_j$, namely, $\mu_1+\dots+\mu_{j-1}+r$.  Then the set
$\SLRT(\la;Z)$ is equal to the set of LR tableaux $\CLRT(\la;R)$,
which was defined in \cite{KS,SW} for the case that each $R_j$ is a
rectangle.
\end{ex}

\begin{ex} \label{compat LR} Let $|R_j|=N_j$
and $N=N_1+\dots+N_L$.  Define the successive subintervals of
$[1,N]$ given by $B_j = [N_1+\dots+N_{j-1}+1,N_1+\dots+N_{j-1}+N_j]$
for $1\le j \le L$.  Let $Z_j$ be any standard tableau of shape
$R_j$ in the alphabet $B_j$.  Then the set $\SLRT(\la;Z)$ is
given by the set of standard tableaux of shape $\la$
that is compatible with a certain labeling of the
cells of the partitions $R_j$ (see \cite{RW,W}).
\end{ex}

\begin{rem} \label{funny compatible}
In the situation of standard tableaux as given in Example \ref{compat LR},
if partitions $R_j$ are rectangles, there is a much
simpler characterization of $\SLRT(\la;Z)$.  Namely, the standard
tableau $S$ of shape $\la$ is in $\SLRT(\la;Z)$ if and only if,
for every index $j$ and every pair of letters $x$ and $y$ in $B_j$,
if $y$ is immediately south (resp. west) of $x$ in $Z_j$, then
in $S$, $y$ is in a row (resp. column) strictly south (resp. west)
of that of $x$.
\end{rem}

\begin{ex} In the characterization of $\SLRT(\la;Z)$ given in
Remark \ref{funny compatible}, 
it is important that each $R_j$ be a rectangle.
Take $\la=(2,2)$ and $R=((1),(2,1))$.  Then for any choice
of $Z_1$ and $Z_2$, $|\SLRT(\la;Z)|=1$, but
there are no tableaux satisfying the criterion in
Remark \ref{funny compatible}.
\end{ex}

Let us fix two canonical choices for the $Z_j$.
The columnwise standard tableau
$\CW(\nu)$ of the partition shape $\nu$ is given by
placing the numbers $\nu^t_1+\dots+\nu^t_{c-1}+1$ through
$\nu^t_1+\dots+\nu^t_c$ from top to bottom in the $c$-th column
for all $c$.  Let $T+x$ denote the tableau obtained
by adding $x$ to every entry of the tableau $T$.
Given a sequence of rectangles $R$, define the sequence
of tableaux $\ZC_j$ ($1\le j\le L$) by
$\ZC_j = \CW(R_j) + |R_{j-1}| + \dots+|R_2|+|R_1|$.
Similarly, one can define the rowwise standard tableau $\RW(\nu)$ of shape
$\nu$ given by placing the numbers $\nu_1+\dots+\nu_{r-1}+1$
through $\nu_1+\dots+\nu_r$ from left to right in the $r$-th row,
for all $r$.  Define $\ZR_j=\RW(R_j)+|R_{j-1}|+\dots+|R_2|+|R_1|$
for $1\le j\le L$.
\begin{definition} Set
\begin{align*}
  \LRT(\la;R) &:= \SLRT(\la;(\ZC_1,\dots,\ZC_L)),\\
  \RLRT(\la;R) &:=\SLRT(\la;(\ZR_1,\dots,\ZR_L)).
\end{align*}
\end{definition}
The set $\LRT(\la;R)$ will be used for the bijection
with rigged configurations.

Observe that for any choice of $Z_j$ and $Z'_j$ (with $R_j$ a
rectangle for all $j$), a bijection $\SLRT(\la;Z)\rightarrow\SLRT(\la;Z')$
is given by relabeling.  Namely, let $S\in\SLRT(\la;Z)$.
Then for each $j$ and each cell $s$ in $R_j$,
replace the letter $Z_j(s)$ in $S$ by $Z'_j(s)$.
In particular there is a bijection
\begin{equation} \label{gamma}
\gamma_R:\LRT(\la;R)\rightarrow \RLRT(\la;R).
\end{equation}
Note that the ordinary transposition of standard tableaux
$\tr:\ST(\la)\to\ST(\la^t)$ restricts to a 
bijection
\begin{equation} \label{tr}
\tr:\LRT(\la;R) \to\RLRT(\la^t;R^t).
\end{equation}
Here $\la^t$ stands for the transpose of the partition $\la$
and $R^t=(R_1^t,\ldots,R_L^t)$ is the sequence of rectangles obtained
by transposing each rectangle of $R$.

\begin{rem} \label{LR defs}
A bijection $\beta_R:\CLRT(\la;R)\rightarrow \LRT(\la;R)$ is given
by a trivial relabeling.
Recall the alphabets $B_1$ through $B_L$ from Example \ref{KS LR}.
Let $T\in \CLRT(\la;R)$.  Then the tableau $\beta_R(T)$
is obtained from $T$ by replacing the occurrences of the
$r$-th largest letter of the subalphabet $B_j$, from left to right,
by the numbers in the $r$-th row of $\ZC_j$.  Alternatively,
$\beta_R = \gamma_R^{-1} \circ \std$ where the bijection
$\std:\CLRT(\la;R)\to\ST(\la)$ is Schensted's standardization map.  Also
$\beta_R = \tr \circ \std \circ \LRtr'$,
where $\LRtr'$ is the LR-transpose map
$\CLRT(\la;R)\rightarrow\CLRT(\la^t;R^t)$
defined in \cite{KS}.
\end{rem}

\begin{definition}\label{def LRtr}
The bijection $\LRtr:\LRT(\la;R)\rightarrow \LRT(\la^t;R^t)$ is given by
$\LRtr:=\tr\circ \gamma_R = \gamma_{R^t} \circ \tr$ where $\gamma_R$ and 
$\tr$ are defined in \eqref{gamma} and \eqref{tr} respectively.
\end{definition}

\begin{ex} \label{ex LR}
Let $R=((3,3),(2,2,2,2),(1,1,1))$ as in \cite[Example 10]{KS}.
We give the set $\LRT(\la;R)$.  $B_1=[1,6]$, $B_2=[7,14]$, $B_3=[15,17]$,
and
\begin{equation*}
  Z_1=\begin{matrix} 1 & 3 & 5 \\ 2 & 4 & 6 \end{matrix} \qquad
  Z_2=\begin{matrix} 7 & 11 \\ 8 & 12 \\ 9 & 13 \\ 10 & 14 \end{matrix}
\qquad Z_3=\begin{matrix} 15 \\ 16 \\ 17 \end{matrix}\,.
\end{equation*}
For $\la = (5,4,3,2,2,1)$, the set $\LRT(\la;R)$ is given below,
listed in order as the
images of the set $\CLRT(\la;R)$ in \cite[Example 11]{KS}
under the bijection $\beta_R$.
\begin{equation*}
\begin{split}
  &\begin{matrix}
    1&3&5&7&11\\
    2&4&6&12& \\
    8&13&15& & \\
    9&14& & & \\
    10&16& & & \\
    17& & & & 
  \end{matrix} \qquad
  \begin{matrix}
    1&3&5&7&11\\
    2&4&6&15& \\
    8&12&16& & \\
    9&13& & & \\
    10&14& & & \\
    17& & & &
  \end{matrix}  \\ \\
  &\begin{matrix}
    1&3&5&11&15\\
    2&4&6&12& \\
    7&13&16& & \\
    8&14& & & \\
    9&17& & & \\
    10& & & & 
  \end{matrix} \qquad
  \begin{matrix}
    1&3&5&11&15\\
    2&4&6&16& \\
    7&12&17& & \\
    8&13& & & \\
    9&14& & & \\
    10& & & & 
  \end{matrix}
\end{split}
\end{equation*}
\end{ex}

\begin{ex} \label{std ex}
When $R=((1)^N)$ and $\la$ is a partition of $N$, $\LRT(\la;R)=\ST(\la)$.
When $R=((\eta_1),\dots,(\eta_L))$, $\CLRT(\la;R)=\CST(\la;\eta)$,
the set of column-strict tableaux of shape $\la$ and content $\eta$,
and $\beta_R$ is Schensted's standardization map.
\end{ex}

It was shown in refs.~\cite{SW,S} that the set
$\CLRT(R)=\cup_{\la}\CLRT(\la;R)$ has the structure of a graded poset
with covering relation given by the $R$-cocyclage and grading function
given by the generalized charge, denoted $c_R$. 
The bijection $\beta_R$ also induces a ranked poset structure on
$\LRT(R)=\cup_{\la}\LRT(\la;R)$; by abuse of notation
we denote its grading function also by $c_R$.
The generalized Kostka polynomial is the generating function of 
LR tableaux with charge statistic~\cite{SW,S}
\begin{equation}\label{gK}
K_{\la R}(q)=\sum_{T\in\LRT(\la;R)} q^{c_R(T)}.
\end{equation}
This extends the charge representation of the Kostka polynomials
$K_{\la\eta}(q)$ of Lascoux and Sch\"utzenberger~\cite{LS Kostka,LS}.
The generalized Kostka polynomials $K_{\la R}(q)$ specialize to
$K_{\la\eta}(q)$ when $R=((\eta_1),\ldots,(\eta_L))$
is a sequence of single rows.

\subsection{Rigged configurations}
\label{sec RC}

Let $R=(R_1,\dots,R_L)$ be a sequence of rectangular partitions such that
$R_j$ has $\mu_j$ rows and $\eta_j$ columns for $1\le j \le L$;
this convention differs from~\cite{KS}.
Denote by $|\la|$ the size of the partition $\la$ and set
$|R|=|R_1|+\cdots+|R_L|$. For $|\la|=|R|$
a $(\la^t;R^t)$-configuration is a sequence of partitions
$\nu = (\nu^{(1)},\nu^{(2)},\dots)$ with the size constraints
\begin{equation} \label{config def}
  |\nu^{(k)}| = \sum_{j>k} \la^t_j -
  	\sum_{a=1}^L \mu_a \max(\eta_a-k,0).
\end{equation}
Define $m_n(\rho)$ as the number of parts of the partition $\rho$
equal to $n$ and $Q_n(\rho)=\rho^t_1+\rho^t_2+\dots+\rho^t_n$,
the size of the first $n$ columns of $\rho$.
Let $\xi^{(k)}(R)$ be the partition whose parts are the heights
of the rectangles in $R$ of width $k$.
The vacancy numbers for the $(\la^t;R^t)$-configuration $\nu$
are the numbers (indexed by $k \ge 1$ and $n\ge 0$) defined by
\begin{equation} \label{vacancy def}
  P^{(k)}_n(\nu) = Q_n(\nu^{(k-1)}) - 2 Q_n(\nu^{(k)}) +
   Q_n(\nu^{(k+1)}) + Q_n(\xi^{(k)}(R))
\end{equation}
where $\nu^{(0)}$ is the empty partition by convention.
In particular $P^{(k)}_0(\nu) = 0$ for all $k\ge 1$.
The $(\la^t;R^t)$-configuration $\nu$ is
admissible if $P^{(k)}_n(\nu) \ge 0$ for all $k,n \ge 1$, and
the set of admissible $(\la^t;R^t)$-configurations is denoted
by $\CO(\la^t;R^t)$. Set
\begin{equation*}
  \cc(\nu)=\sum_{k,n\ge 1} \alpha_n^{(k)}(\alpha_n^{(k)}-\alpha_n^{(k+1)})
\end{equation*}
where $\alpha_n^{(k)}$ is the size of the $n$-th column in
$\nu^{(k)}$. 
Finally, define the $q$-binomial as
\begin{equation*}
\qbin{m+n}{m}=\frac{(q)_{m+n}}{(q)_m(q)_n}
\end{equation*}
for $m,n\in\Int$ and zero otherwise where $(q)_m=(1-q)(1-q^2)\cdots
(1-q^m)$.

With this notation we can state the following quasi-particle
expression of the generalized Kostka polynomials conjectured
in~\cite{KS,SW}, stemming from the analogous expression of
Kirillov and Reshetikhin~\cite{KR} for the Kostka polynomial.

\begin{thm}[Quasi-particle representation] \label{thm QP}
For $R$ a sequence of rectangles and $\la$ a partition
such that $|\la|=|R|$
\begin{equation}\label{qp}
  K_{\la R}(q)=\sum_{\nu\in\CO(\la^t;R^t)} q^{\cc(\nu)}
   \prod_{k,n\ge 1} \qbin{P_n^{(k)}(\nu)+m_n(\nu^{(k)})}{m_n(\nu^{(k)})}.
\end{equation}
\end{thm}

Expression \eqref{qp} can be reformulated as the generating function
over rigged configurations. To this end we need
to define certain labelings of the rows of the
partitions in a configuration.
For this purpose one should view a partition as
a multiset of positive integers.
A rigged partition is by definition a finite multiset of
pairs $(n,x)$ where $n$ is a positive integer and
$x$ is a nonnegative integer.  The pairs $(n,x)$ are referred to
as strings; $n$ is referred to as the
length or size of the string and $x$ as the label or
quantum number of the string.  A rigged partition is
said to be a rigging of the partition $\rho$ if
the multiset consisting of the sizes of the strings,
is the partition $\rho$.  So a rigging of $\rho$
is a labeling of the parts of $\rho$ by nonnegative integers,
where one identifies labelings that differ only by
permuting labels among equal-sized parts of $\rho$.

A rigging $J$ of the
$(\la^t;R^t)$-configuration $\nu$ is a sequence of
riggings of the partitions $\nu^{(k)}$ such that
every label $x$ of a part of $\nu^{(k)}$ of size $n$,
satisfies the inequalities
\begin{equation} \label{rigging def}
  0 \le x \le P^{(k)}_n(\nu).
\end{equation}
The pair $(\nu,J)$ is called a rigged configuration.
The set of riggings of admissible $(\la^t;R^t)$-configurations
is denoted by $\RC(\la^t;R^t)$.
Let $(\nu,J)^{(k)}$ be the $k$-th rigged partition
of $(\nu,J)$.  A string $(n,x)\in (\nu,J)^{(k)}$
is said to be singular if $x=P^{(k)}_n(\nu)$, that is,
its label takes on the maximum value.

\begin{rem} \label{ordering} Observe that the definition
of the set $\RC(\la^t;R^t)$ is completely insensitive to
the order of the rectangles in the sequence $R$.
However the notation involving the sequence $R$
is useful when discussing the bijection 
$\phib_R:\LRT(\la;R)\to\RC(\la^t;R^t)$,
since the ordering on $R$ is essential in the definition
of $\LRT(\la;R)$.
\end{rem}

The set of rigged configurations is endowed with a natural
statistic $\cc$ \cite[(3.2)]{KS} defined by
\begin{equation} \label{RC charge}
  \cc(\nu,J)=\cc(\nu)+\sum_{k,n\ge 1} |J_n^{(k)}|
\end{equation}
for $(\nu,J)\in\RC(\la^t;R^t)$.
Here $J_n^{(k)}$ denotes the partition inside the rectangle
of height $m_n(\nu^{(k)})$ and width $P_n^{(k)}(\nu)$ given
by the labels of the parts of $\nu^{(k)}$ of size $n$.
Since the $q$-binomial $\qbins{P+m}{m}$ is the generating
function of partitions with at most $m$ parts each not
exceeding $P$, Theorem \ref{thm QP} is equivalent to
the following theorem.

\begin{thm}[Rigged configuration representation]\label{thm RC}
For $R$ a sequence of rectangles and $\la$ a partition
such that $|\la|=|R|$
\begin{equation}\label{rc}
  K_{\la R}(q)=\sum_{(\nu,J)\in\RC(\la^t;R^t)} q^{\cc(\nu,J)}.
\end{equation}
\end{thm}
The proof of this theorem follows from the bijection
$\phib_R:\LRT(\la;R)\to\RC(\la^t;R^t)$ of 
Definition-Proposition~\ref{phi def}
and Theorem~\ref{thm charge} below.

\section{Maps on rectangles, Littlewood--Richardson tableaux and
rigged configurations}
\label{sec operations}

In this section we define several operations on sequences
of rectangles and their counterparts on the sets of LR tableaux
and rigged configurations. These operations underly the recursive 
definition of the bijection $\phib_R:\LRT(\la;R)\to\RC(\la^t;R^t)$ 
as given in Definition-Proposition~\ref{phi def}.
A summary of the definitions and results of this section is given
in Table~\ref{table summarize}.

\subsection{Operations on sequences of rectangles}
\label{rect seq}
	
Let $R=(R_1,R_2,\dots,R_L)$ be a sequence of rectangles
such that $R_j=(\eta_j^{\mu_j})$ has $\mu_j$ rows and $\eta_j$
columns.
Let $\Rht$ be the sequence of rectangles obtained from $R$ by splitting
off the last column of $R_L$; formally, $\Rht_j = R_j$ for $1\le j\le L-1$,
$\Rht_L = ((\eta_L-1)^{\mu_L})$ and $\Rht_{L+1}=(1^{\mu_L})$.
Note that if the last rectangle of $R$ is a single column, then
(ignoring the empty rectangle) $\Rht=R$.
If the last rectangle of $R$ is a single column,
let $\Rb$ be given by removing one cell from the column $R_L$;
$\Rb_j=R_j$ for $1\le j \le L-1$ and $\Rb_L = (1^{\mu_L-1})$.
Let $\Rck$ be given by splitting off the first column of $R_1$;
if $R_1$ is a single column, then $\Rck=R$.
If $R_1$ is a single column, let $\Rt$ be given by removing
one cell from the column $R_1$.
Finally, let $R^\ev = (R_L,\dots,R_2,R_1)$ denote the reverse of $R$.

\begin{rem} \label{bij ind}
Given any sequence of rectangles, there is a unique
sequence of transformations of the form $R\rightarrow \Rht$ or
$R\rightarrow \Rb$ resulting in the empty sequence,
where $R\rightarrow \Rht$ is only used when the last rectangle
of $R$ has more than one column.
\end{rem}

\subsection{Maps between sets of LR tableaux}\label{sec maps LR}
For each operation on sequences of rectangles,
there is a natural (injective) map
on the corresponding sets of LR tableaux of a fixed shape.

Observe that there are inclusions
\begin{equation*}
\begin{split}
  \iht: \LRT(\la;R) &\rightarrow \LRT(\la;\Rht) \\
  \ick: \LRT(\la;R) &\rightarrow \LRT(\la;\Rck)
\end{split}
\end{equation*}
which correspond to the transformations $R\to\Rht$ and $R\to\Rck$
on rectangles.

Recall that $\ST(\la)$ denotes the set of standard tableaux of 
shape $\la$ and define
$\ST(\la^-) = \displaystyle\bigcup_{\rho\lessdot\la} \ST(\rho)$,
where $\rho$ and $\la$ are partitions and $\rho \lessdot\la$ means 
that $\rho\subset\la$ and $\la/\rho$ is a single cell.
There is a bijection
\begin{equation} \label{std minus}
\begin{split}
  -:\ST(\la) &\rightarrow \ST(\la^-) \\
  S &\mapsto S^-
\end{split}
\end{equation}
where $S^-$ is the standard tableau obtained by removing the
maximum letter from $S$.  Obviously $S$ is uniquely determined
by its shape and the tableau $S^-$.  If the last rectangle of $R$ is
a single column, write
\begin{equation*}
\LRT(\la^-;\Rb) = \displaystyle\bigcup_{\rho \lessdot \la} \LRT(\rho;\Rb).
\end{equation*}
The following result is an immediate consequence of the definitions.

\begin{prop} \label{LR minus inj} Suppose the last rectangle of $R$ is
a single column.
\begin{enumerate}
\item The map \eqref{std minus} restricts to an injection
$\ib:\LRT(\la;R) \rightarrow \LRT(\la^-;\Rb)$.
\item If $\mu_L=1$ then $\ib$ is bijective.
\item Suppose $\mu_L > 1$ and $T\in \LRT(\la^-;\Rb)$ such that
the cell $\la/\shape(T)$ is in the $r$-th row.
Then $T$ is in the image of $\ib$ if and only if the cell
$\shape(T)/\shape(T^-)$ is in the $r'$-th row with $r' < r$.
\end{enumerate}
\end{prop}
The injection $\ib:\LRT(\la;R) \rightarrow \LRT(\la^-;\Rb)$
corresponds to the operation $R\to\Rb$ on rectangles.
Next we describe a dual operation to $S\mapsto S^-$ giving
rise to the analogue of $R\to\Rt$ on LR tableaux.

Fix partitions $\sigma \subset \la$ such that the skew shape
$\la/\sigma$ has two cells.  Consider the set of saturated chains
in Young's lattice of partitions under inclusion, that have
maximum element $\la$ and minimum element $\sigma$:
\begin{equation}
  \CC[\sigma,\la] := \{ \sigma \lessdot \rho \lessdot \la \}.
\end{equation}
If $\la/\sigma$ is connected (that is, its two cells are adjacent) then
$\CC[\sigma,\la]$ is a singleton, whose intermediate partition
$\rho$ is obtained by adjoining the
inner of these two cells to $\sigma$ or removing the outer of the
two cells from $\la$.  If $\la/\sigma$ is
disconnected (that is, its two cells are not adjacent)
then $\CC[\sigma,\la]$ has exactly two elements, whose intermediate
partitions are obtained by
adjoining either of the two cells to $\sigma$ or removing either
of the two cells from $\la$.  Let $\tau=\tau_{\sigma,\la}$ be the
involution on $\CC[\sigma,\la]$ that has order two if
$|\CC[\sigma,\la]|=2$ (and of course must
be the identity when $|\CC[\sigma,\la]|=1$).

\begin{defprop} \label{D def} For $|\la|>0$,
there is a unique bijection $D:\ST(\la)\rightarrow \ST(\la^-)$ denoted
$S\mapsto S^D$ such that:
\begin{enumerate}
\item[(D1)] If $|\la|=1$ then $1^D=\emptyset$.
\item[(D2)] If $|\la|=N>1$, $S\in \ST(\la)$, then
$S^D$ is uniquely defined by the properties that
$S^{D-} = S^{-D}$, and the shape of $S^D$ is the intermediate partition
in the chain
\begin{equation*}
  \tau(\sigma \lessdot \shape(S^-) \lessdot\la)
\end{equation*}
where $\sigma=\shape(S^{-D})$.
\end{enumerate}
\end{defprop}
\noindent $S$ is uniquely determined by its shape and the tableau $S^D$.

Using the characterization of $D$ in Definition-Proposition
\ref{D def}, it can be shown~\cite{Sc} that $D$ is computed by the following 
well-known tableau algorithm.

\begin{lem} \label{D jeu} For $S\in\ST(\la)$ with
$|\la|=N$, $S^D$ is obtained by removing the number $1$
from $S$, subtracting one from each entry, and sliding the
resulting skew tableau to partition shape, that is,
$S^D = P(S|_{[2,N]}-1)$.  More generally for any $0\le i\le N$,
\begin{equation*}
  S^{D^i} = P(S|_{[i+1,N]}-i).
\end{equation*}
\end{lem}

Say that the index $i$ is a descent of the standard tableau $S$ if
$i+1$ appears in a later row in $S$ than $i$ does.

\begin{lem} \label{D descent} Suppose $S\in \ST(\la)$ and $|\la|=N \ge 2$.
Let $r$ and $r'$ be the rows of the cells
$\shape(S)/\shape(S^D)$ and $\shape(S^D)/\shape(S^{D^2})$.
Then $1$ is a descent of $S$ if and only if $r' < r$.
\end{lem}
\begin{proof} The partition $\shape(S^{D^2})$ is calculated by
Lemma \ref{D jeu} as the shape obtained by sliding the skew tableau
$S|_{[3,N]}$ to partition shape, first into the cell of $S$ containing
the letter $2$ (vacating the cell $s$ say) and then into the cell of $S$
containing the letter $1$ (vacating $s'$ say).  Moreover
$s=\shape(S)/\shape(S^D)$ and $s'=\shape(S^D)/\shape(S^{D^2})$.

However $\shape(S^{D^2})$ can be calculated another way.
Take the two-letter tableau $U=S|_{[1,2]}$ and slide it
to the southeast into the cells of $S|_{[3,N]}$ occupied by $3$, then $4$,
etc., producing the skew tableau $V$ say.  Then
$\shape(S^{D^2})$ is given by $\shape(S)-\shape(V)$.
It is clear that the cells of $V$ containing $1$ and $2$ are
$s'$ and $s$ respectively.  But sliding preserves Knuth equivalence,
so $2$ is in a later row than $1$ in $U$ (or $S$) if and only if
it is in $V$.
\end{proof}

Suppose the first rectangle of $R$ is a single column.  Define
\begin{equation*}
  \LRT(\la^-;\Rt) = \bigcup_{\rho \lessdot \la} \LRT(\rho;\Rt).
\end{equation*}
The analogue of $R\to\Rt$ is given by the following proposition.
\begin{prop} \label{LR D inj} Suppose the first rectangle of $R$ is a
single column.
\begin{enumerate}
\item The bijection $D:\ST(\la)\rightarrow\ST(\la^-)$ restricts to an
injection $D:\LRT(\la;R) \rightarrow \LRT(\la^-;\Rt)$.
\item If $\mu_1=1$ then $D$ is bijective.
\item Suppose that $\mu_1 > 1$ and $T\in \LRT(\la^-;\Rt)$ such that
the cell $\la/\shape(T)$ is in the $r$-th row.
Then $T$ is in the image of $D$ if and only if the cell
$\shape(T)/\shape(T^D)$ is in the $r'$-th row with $r' < r$.
\end{enumerate}
\end{prop}
\begin{proof} 1 follows by Lemma \ref{D jeu},
the definition of $\LRT$, and the fact that restriction to subintervals
preserves Knuth equivalence.  For 2 and 3,
let $T\in \LRT(\la^-;\Rt)$.  If there is a tableau $S\in \LRT(\la;R)$
such that $T=S^D$, then by Lemma \ref{D jeu} $S$ is obtained from $T$ by
sliding $T+1$ to the southeast into the
cell $\la/\shape(T)$, and placing a $1$ in the northwest corner.
If $\mu_1=1$ then immediately $S\in \LRT(\la;R)$, proving 2.
If $\mu_1>1$ then by construction $S\in \LRT(\la;R)$ if and only if
$1$ is a descent of $S$. By Lemma~\ref{D descent} point~3 follows.
\end{proof}

Let us discuss carefully the commutation of $-$ and $D$ in (D2)
in Definition-Proposition~\ref{D def}. Define
\begin{equation*}
  \ST(\la^{--}) = \bigcup_{\sigma\lessdot\rho\lessdot\la} \ST(\sigma),
\end{equation*}
where the right hand side is a disjoint union.  By definition,
for each chain $\CC = (\sigma \lessdot \rho\lessdot\la)$,
there is a distinguished copy of $\ST(\sigma)$ in
$\ST(\la^{--})$ denoted $\ST(\sigma)_\CC$, and for $T\in\ST(\sigma)$,
denote its copy in $\ST(\sigma)_\CC$ by $T_\CC$.  The point is that
when $\la/\sigma$ is disconnected there are two copies of
$\ST(\sigma)$ in $\ST(\la^{--})$ and they are distinguished by
the chain $\CC$.  By abuse of notation let $\tau$ denote the involution on
$\ST(\la^{--})$ that sends $T_\CC\rightarrow T_{\tau(\CC)}$.
In this notation, which concerns itself with
the intermediate partition in the chain $\CC$, (D2) is expressed as
\begin{equation} \label{std D and minus}
  - \circ D = \tau \circ D \circ -
\end{equation}
viewed as maps $\ST(\la)\rightarrow \ST(\la^{--})$.

Suppose both the first and last rectangles of $R$ are single columns.
Of course $\Rbt=\Rtb$.  Define
\begin{equation*}
  \LRT(\la^{--};\Rbt) = \bigcup_{\sigma\lessdot\rho\lessdot\la}
	  \LRT(\sigma;\Rbt).
\end{equation*}
The identity of maps \eqref{std D and minus} restricts to the
identity of maps $\LRT(\la;R)\rightarrow\LRT(\la^{--};\Rbt)$:
\begin{equation} \label{D and minus}
  - \circ D = \tau \circ D \circ -
\end{equation}

Define Sch\"utzenberger's evacuation map
\begin{equation} \label{ev def}
\begin{split}
  \ev: \ST(\la) &\rightarrow \ST(\la) \\
  S&\mapsto S^\ev
\end{split}
\end{equation}
where $S^\ev\in\ST(\la)$ is defined by
$S^\ev = S$ if $\la = \emptyset$ and $S^{\ev-} = S^{D\ev}$.

The following result is well-known and easy to prove.

\begin{prop} \label{ev and} $\ev$ is an involution satisfying
\begin{equation*}
\ev \circ - = D \circ \ev \qquad \text{and} \qquad
\ev \circ D = - \circ \ev.
\end{equation*}
\end{prop}

\begin{lem} \label{ev res} Let $|\la|=N$ and $S\in \ST(\la)$.  Then
for any $1\le i \le j \le N$,
\begin{equation*}
  P(S^\ev|_{[i,j]}-(i-1)) = P(S|_{[N+1-j,N+1-i]}-(N-j))^\ev.
\end{equation*}
\end{lem}
\begin{proof} By Lemma \ref{D jeu}, \eqref{std D and minus},
and the definitions,
\begin{equation*}
\begin{split}
  P(S|_{[i,j]}-(i-1)) &= P((S|_{[1,j]})|_{[i,j]}-(i-1)) \\
  &= S^{-^{N-j}D^{i-1}} = S^{D^{i-1}-^{N-j}}.
\end{split}
\end{equation*}
Using this and Proposition \ref{ev and}, we have
\begin{equation*}
\begin{split}
  P(S^\ev|_{[i,j]}-(i-1))
  &=S^{\ev D^{i-1}-^{N-j}}
   = S^{-^{i-1}D^{N-j}\ev} \\
  &= P(S|_{[N+1-j,N+1-i]}-(N-j))^\ev.
\end{split}
\end{equation*}
\end{proof}

Observe that the tableau in the singleton set $\LRT(R_1,(R_1))$
evacuates to itself.  Using this fact,
induction, Lemma \ref{ev res}, the definition of $\LRT$,
and the fact that Knuth equivalence is preserved under restriction
to subintervals, it follows that the evacuation map \eqref{ev def}
restricts to a bijection
\begin{equation*}
  \ev: \LRT(\la;R) \rightarrow \LRT(\la;R^\ev).
\end{equation*}

Of course ${\Rht}^\ev=(R^\ev)\ck$.
By definition the following diagram commutes
\begin{equation} \label{hat and ev}
\begin{CD}
  \LRT(\la;R) @>{\ev}>> \LRT(\la;R^\ev) \\
  @V{\iht}VV	@VV{\ick}V \\
  \LRT(\la;\Rht) @>>{\ev}> \LRT(\la;{\Rht}^\ev)
\end{CD}
\end{equation}
where $\ick$ is defined with respect to the sequence
of rectangles $R^\ev$.
Since $\ev$ is an involution, one may exchange the roles of
$\iht$ and $\ick$, of $\Rht$ and $\Rck$, and of $R$ and $R^\ev$.

\subsection{Maps between sets of rigged configurations}\label{sec maps RC}
For the various transformations of
sequences of rectangles, one has the following
maps between the corresponding sets of rigged configurations.

Define the map
\begin{equation*}
  \jht: \RC(\la^t;R^t) \rightarrow \RC(\la^t;(\Rht)^t)
\end{equation*}
by declaring that $\jht(\nu,J)$ is obtained from $(\nu,J)\in\RC(\la^t;R^t)$
by adding a singular string of length $\mu_L$ to each
of the first $\eta_L - 1$ rigged partitions.
Note that $\jht$ is the identity map if $R_L$ is a single column.

\begin{lem} \label{j hat} $\jht$ is a well-defined injection
that preserves the vacancy numbers of the underlying configurations.
\end{lem}
\begin{proof}
Let $(\nuht,J\htt) = \jht(\nu,J)$.
It is enough to show that $\nuht$ is a $(\la^t;(\Rht)^t)$-configuration
that has the same vacancy numbers as $\nu$.
First we verify that the partitions in $\nuht$ have the correct size.
To this end set $\chi(\mathrm{true})=1$ and $\chi(\mathrm{false})=0$. Then
\begin{equation*}
\begin{split}
  |\nuht^{(k)}| &= \chi(k < \eta_L)\, \mu_L + |\nu^{(k)} | \\
  &= \chi(k < \eta_L)\, \mu_L +
	\sum_{j>k} \la^t_j - \sum_{a=1}^L \mu_a \max(\eta_a-k,0) \\
  &= \sum_{j>k} \la^t_j - \sum_{a=1}^{L-1} \mu_a \max(\eta_a-k,0) \\
  &\quad - \mu_L \{\max(\eta_L-1-k,0) + \max(1-k,0)\}.
\end{split}
\end{equation*}
Next we check that the vacancy numbers remain the same.
Note that
\begin{equation*}
  Q_n(\nuht^{(k)}) - Q_n(\nu^{(k)}) = \chi(1 \le k < \eta_L) \min(\mu_L,n),
\end{equation*}
valid for $k\ge 0$.  Then for $k,n\ge 1$,
\begin{equation*}
\begin{split}
  &P^{(k)}_n(\nuht) - P^{(k)}_n(\nu) \\
  =\,\, &\min(\mu_L,n) \{\chi(1\le k-1 < \eta_L) - 2 \chi(1\le k < \eta_L) +
	\chi(1\le k+1 < \eta_L)\} \\
  &\,\,+\min(\mu_L,n)(\delta_{\eta_L-1,k} + \delta_{1,k}-\delta_{\eta_L,k}) \\
  =\,\,& 0
\end{split}
\end{equation*}
where $\delta_{a,b}=\chi(a=b)$.
\end{proof}

Define the map
\begin{equation*}
  \jck: \RC(\la^t;R^t) \rightarrow \RC(\la^t;(\Rck)^t)
\end{equation*}
by declaring that $\jck(\nu,J)$ is obtained from $(\nu,J)\in\RC(\la^t;R^t)$
by adding a string with label zero and length $\mu_1$ to each 
of the first $\eta_1 - 1$ rigged partitions.  If $R_1$ is a single column
then $\jck$ is the identity.

By Remark \ref{ordering}, the proof of Lemma \ref{j hat}
also shows that:

\begin{lem} \label{j check} $\jck$ is a well-defined injection
that preserves the vacancy numbers of the underlying configurations.
\end{lem}

Define the involution
\begin{equation}\label{com def}
  \com_R:\RC(\la^t;R^t) \to \RC(\la^t;R^t)
\end{equation}
by $\com_R(\nu,J)=(\nu,P-J)$, which is shorthand for
saying that $\com_R$ preserves the
underlying configuration and replaces the label $x$ of
a string $(n,x)\in (\nu,J)^{(k)}$
by its colabel, which is by definition the number
$P^{(k)}_n(\nu) - x$, its complement with respect to the
vacancy number of the string.
In light of Remark \ref{ordering}, one may also define
\begin{equation*}
  \comev_R: \RC(\la^t;R^t) \to \RC(\la^t;R^{\ev t})
\end{equation*}
which complements the quantum numbers as $\com_R$ and in
addition reverses the sequence of rectangles.

By definition the following diagram commutes
\begin{equation} \label{hat and theta}
\begin{CD}
  \RC(\la^t;R^t) @>{\comev_R}>> \RC(\la^t;R^{\ev t}) \\
  @V{\jht}VV	@VV{\jck}V \\
  \RC(\la^t;(\Rht)^t) @>>{\comev_{\Rht}}> \RC(\la^t;(\Rht)^{\ev t})
\end{CD}
\end{equation}
where $\jck$ acts with respect to the reversed sequence
of rectangles $R^\ev$.  Since $\comev_R$ is an involution
one may exchange the roles of $\jht$ and $\jck$, of
$\Rht$ and $\Rck$, and of $R$ and $R^\ev$.

Suppose the last rectangle of $R$ is a single column.
Define the set
\begin{equation*}
  \RC(\la^{-t};\Rb^t) = \bigcup_{\rho\lessdot \la} \RC(\rho^t;\Rb^t).
\end{equation*}
The key algorithm on rigged configurations is given by the map
\begin{equation*}
  \jb: \RC(\la^t;R^t) \rightarrow \RC(\la^{-t};\Rb^t),
\end{equation*}
defined as follows.  Let $(\nu,J) \in \RC(\la^t;R^t)$.
Define $\lb^{(0)} = \mu_L$.  By induction select the singular
string in $(\nu,J)^{(k)}$ whose length 
$\lb^{(k)}$ is minimal such that $\lb^{(k-1)} \le \lb^{(k)}$.
Let $\rkb(\nu,J)$ denote the smallest $k$ for which no such
string exists, and set $\lb^{(k)}=\infty$ for $k\ge \rkb(\nu,J)$.
Then $\jb(\nu,J)=(\nub,\Jb)$ is obtained from $(\nu,J)$ by
shortening each of the selected singular strings by one,
changing their labels so that they remain singular,
and leaving the other strings unchanged.

Let us compute the change in vacancy numbers under $\jb$.
Recalling that $\nu^{(0)}=\emptyset$, observe that
\begin{equation*}
  Q_n(\nu^{(k)}) - Q_n(\nub^{(k)}) =
  \chi(k \ge 1) \chi(n \ge \lb^{(k)})
\end{equation*}
for $k\ge 0$.  Then for $k \ge 2$ and $n \ge 0$ we have
\begin{equation*}
\begin{split}
  P^{(k)}_n(\nu) - P^{(k)}_n(\nub)
  =\,\, &\chi(n\ge\lb^{(k-1)})-2\chi(n\ge\lb^{(k)})+\chi(n\ge\lb^{(k+1)}) \\
  =\,\, &\chi(\lb^{(k-1)} \le n < \lb^{(k)}) -
  \chi(\lb^{(k)} \le n < \lb^{(k+1)}).
\end{split}
\end{equation*}
Recall that $\lb^{(0)} = \mu_L$.  For $k=1$ and $n\ge0$,
\begin{equation*}
\begin{split}
  P^{(1)}_n(\nu) - P^{(1)}_n(\nub)
  =\,\,& -2\chi(n\ge\lb^{(1)})+\chi(n\ge\lb^{(2)}) + \chi(n\ge \mu_L) \\
  =\,\,& \chi(\lb^{(0)} \le n < \lb^{(1)}) -
  \chi(\lb^{(1)} \le n < \lb^{(2)}).
\end{split}
\end{equation*}
Therefore for all $k\ge1$ and $n\ge 0$ we have
\begin{equation} \label{db vacancy}
  P^{(k)}_n(\nu) - P^{(k)}_n(\nub) =
  \chi(\lb^{(k-1)} \le n < \lb^{(k)}) -
  \chi(\lb^{(k)} \le n < \lb^{(k+1)}).
\end{equation}

As before $m_n(\rho)$ denotes the number of parts of size $n$ in the 
partition $\rho$.
One may easily verify (see also \cite[Appendix]{KS}) that for all $k,n\ge 1$
\begin{equation} \label{vacancy ddiff}
\begin{split}  &- P^{(k)}_{n-1}(\nu) + 2 P^{(k)}_n(\nu) - P^{(k)}_{n+1}(\nu) \\
  = \quad&m_n(\nu^{(k-1)}) - 2 m_n(\nu^{(k)}) + m_n(\nu^{(k+1)})
     +m_n(\xi^{(k)}(R)) \\
  \ge \quad&m_n(\nu^{(k-1)}) - 2 m_n(\nu^{(k)}) + m_n(\nu^{(k+1)}).
\end{split}
\end{equation}
In particular the vacancy numbers have the partial convexity property
(see \cite[(11.1)]{KS})
\begin{equation} \label{convex}
  P^{(k)}_n(\nu) \ge 1/2 (P^{(k)}_{n-1}(\nu) + P^{(k)}_{n+1}(\nu))
\quad\text{ if $m_n(\nu^{(k)})=0$.}  
\end{equation}
Repeated use of \eqref{convex} leads to the following result.

\begin{lem} \label{lower bound}
Let $\nu$ be an admissible configuration and
$m_n(\nu^{(k)}) = 0$ for $a < n < b$.
\begin{enumerate}
\item $P^{(k)}_n(\nu) \ge \min(P^{(k)}_a(\nu),P^{(k)}_b(\nu))$ for
$a\le n\le b$.
\item If $P^{(k)}_c(\nu) = 0$ for some $a < c < b$ then
$P^{(k)}_n(\nu)=0$ for all $a\le n \le b$.
\item If $P^{(k)}_c(\nu)=P^{(k)}_{c+1}(\nu)=1$ for some $c$
such that $a\le c <b$ then $P^{(k)}_n(\nu)=1$ for all $a < n < b$.
\end{enumerate}
\end{lem}

The following proposition is important for the definition of the bijection
$\phib_R$ between rigged configurations and LR tableaux to be defined
in Section \ref{bij sec}.
\begin{prop} \label{j bar}
Let $(\nu,J) \in \RC(\la^t;R^t)$ where the last rectangle of $R$ is a single
column.
\begin{enumerate}
\item The map $\jb$ is a well-defined injection such that
$\jb(\nu,J) \in \RC(\rho^t;\Rb^t)$ where $\rho$ is obtained from
$\la$ by removing the corner cell in the column of index
$\rkb(\nu,J)$.
\item If $\mu_L=1$ then $\db$ is bijective.
\item If $\mu_L>1$ then $(\nu',J')\in \RC(\rho^t;\Rb^t)$
is in the image of $\jb$ if and only if $\rkb(\nu',J') \ge \rkb(\nu,J)$.
\end{enumerate}
\end{prop}

\begin{proof}
To prove that $\jb$ is well-defined it needs to be shown that 
$(\nub,\Jb)=\jb(\nu,J)$ is an admissible rigged 
$(\rho^t;\Rb^t)$-configuration. Let us first show that
$\rho$ obtained from $\la$ by removing the corner cell in column
of index $\rkb(\nu,J)$ is indeed a partition.
Assume the contrary. This means that $\la^t_r=\la_{r+1}^t$ where
$r=\rkb(\nu,J)$. By~\cite[(11.2)]{KS}
$P_n^{(k)}(\nu)=\la_k^t-\la_{k+1}^t$ for large $n$,
so that $P_n^{(r)}(\nu)=0$. Let $\ell$ be the size of the largest part 
of $\nu^{(r)}$. Then by Lemma~\ref{lower bound} it follows that
$P_n^{(r)}(\nu)=0$ for all $n\ge \ell$. Since $m_n(\nu^{(r)})=0$
for $n>\ell$ and $P_n^{(k)}(\nu)\ge 0$ for all $k\ge 1$, $n\ge 0$, 
inequality~\eqref{vacancy ddiff} implies in particular that 
$m_n(\nu^{(r-1)})=0$ for $n>\ell$. This means that $1\le \lb^{(r-1)}
\le \ell$. Since $P_\ell^{(r)}(\nu)=0$ and $m_\ell(\nu^{(r)})>0$
there is a singular string of length $\ell$ in $\nu^{(r)}$ weakly
bigger than $\lb^{(r-1)}$. However, this contradicts the assumption that
$r=\rkb(\nu,J)$ which would imply that there is no such singular 
string. Hence $\rho$ is a partition.
By the definition of $\rkb(\nu,J)$ it is clear that 
$|\nub^{(k)}|=|\nu^{(k)}|-\chi(k<\rkb(\nu,J))$.
Since $\eta_L=1$ it follows from \eqref{config def} that $\nub$ is 
a $(\rho^t;\Rb^t)$-configuration.

Next we need to show that $(\nub,\Jb)$ is admissible.
Denote by $J_n^{(k)}(\nu,J)$ the maximal rigging occurring in the 
strings of length $n$ in $(\nu,J)^{(k)}$ (which is set to zero if $n$ 
does not appear as a part in $\nu^{(k)}$). Then 
to prove the admissibility of $(\nub,\Jb)$ we need to show for all 
$n,k \ge 1$ that 
\begin{equation} \label{adm cond}
  0 \le J_n^{(k)}(\nub,\Jb) \le P_n^{(k)}(\nub).
\end{equation}
Fix $k \ge 1$. Since only one string of size $\lb^{(k)}$ changes in the
transformation $(\nu,J)^{(k)}\to (\nub,\Jb)^{(k)}$ one finds that 
$J_n^{(k)}(\nub,\Jb)=J_n^{(k)}(\nu,J)$ for $1\le n<\lb^{(k)}-1$, 
$J_n^{(k)}(\nub,\Jb)=P_n^{(k)}(\nub)$ for $n=\lb^{(k)}-1$ and 
$0\le J_n^{(k)}(\nub,\Jb)\le J_n^{(k)}(\nu,J)$ for $n\ge \lb^{(k)}$.
Hence by \eqref{db vacancy} the inequality \eqref{adm cond} can only be 
violated when $\lb^{(k-1)}\le n < \lb^{(k)}$. 
By the construction of $\lb^{(k)}$ there are no singular strings of length 
$n$ in $(\nu,J)^{(k)}$ for $\lb^{(k-1)}\le n<\lb^{(k)}$. This means that 
$J_n^{(k)}(\nu,J)\le P_n^{(k)}(\nu)-1$ if $n$ occurs as a part in $\nu^{(k)}$,
that is $m_n(\nu^{(k)})>0$. Hence due to \eqref{db vacancy} the condition 
\eqref{adm cond} is fulfilled for these $n$.
It remains to prove that $P_n^{(k)}(\nub)\ge 0$ for all $n$ such that 
$m_n(\nu^{(k)})=0$ and $\lb^{(k-1)}\le n<\lb^{(k)}$. Note that 
$m_n(\nub^{(k)})=0$ if $m_n(\nu^{(k)})=0$ for $\lb^{(k-1)}\le n <\lb^{(k)}-1$.
By \cite[Lemma 10]{KS} it suffices to prove \eqref{adm cond} for all 
$k$ and $n$ such that $m_n(\nub^{(k)})>0$. Therefore the only remaining case 
for which \eqref{adm cond} might be violated occurs when
\begin{equation*}
  m_{\ell-1}(\nu^{(k)})=0, \quad P_{\ell-1}^{(k)}(\nu)=0, \quad
   \lb^{(k-1)}<\ell \quad \text{and $\ell$ finite}
\end{equation*}
where $\ell=\lb^{(k)}$. We show that these conditions cannot be met
simultaneously. Let $p<\ell$ be maximal such that $m_p(\nu^{(k)})>0$;
if no such $p$ exists set $p=0$. By Lemma \ref{lower bound} 
$P_{\ell-1}^{(k)}(\nu)=0$ is only possible if $P_n^{(k)}(\nu)=0$
for all $p\le n\le \ell$. By \eqref{vacancy ddiff} we find that 
$m_n(\nu^{(k-1)})=0$ for $p<n<\ell$. Since $\lb^{(k-1)}<\ell$ this implies
that $\lb^{(k-1)}\le p$. If $p=0$ this contradicts the condition
$\lb^{(k-1)}\ge 1$. Hence assume that $p>0$. Since $P_p^{(k)}(\nu)=0$
and $m_p(\nu^{(k)})>0$ there is a singular string of length $p$ in 
$(\nu,J)^{(k)}$ and therefore $\lb^{(k)}=p$. However, this contradicts 
$p<\ell=\lb^{(k)}$. This concludes the proof of the admissibility of 
$(\nub,\Jb)$ and also the proof of the well-definedness of $\jb$.

For the proof of the injectivity of $\jb$, and points 2 and 3 we require
the algorithm $\jbi$ defined on $\RC(\la^{-t};\Rb^t)$ as follows.
Recall that $(\nub,\Jb)\in\RC(\la^{-t};\Rb^t)$ means that 
$(\nub,\Jb)\in\RC(\rho^t;\Rb^t)$ for some $\rho\lessdot\la$. Suppose that
the cell $\la/\rho$ has column index $c$ in $\la$. Set $\lbi^{(k)}=\infty$
for $k\ge c$. For $1\le k<c$ select by downward induction the singular 
string in $(\nub,\Jb)^{(k)}$ whose length $\lbi^{(k)}$ is maximal such that
$\lbi^{(k)}\le \lbi^{(k+1)}$; set $\lbi^{(k)}=0$ if no such string exists.
Then $(\nu,J)=\jbi(\nub,\Jb)$ is obtained from $(\nub,\Jb)$ by adding
one box to the selected strings (and adding a string of length one if
$\lbi^{(k)}=0$) with labels such that they remain singular, and leaving
all other strings unchanged.

It is obvious from the constructions and \eqref{db vacancy}
that for $(\nu,J)\in\RC(\la^t;R^t)$, $\jbi\circ\jb(\nu,J)=(\nu,J)$ 
since $\lbi^{(k)}=\lb^{(k)}-1$.
This proves that $\jb$ is an injection and concludes the proof of point 1.

It follows immediately from the definition of $\jb$ that 
$\rkb(\jb(\nu,J))\ge \rkb(\nu,J)$ for $(\nu,J)\in\RC(\la^t;R^t)$ where
$\mu_L>1$. Hence, if $(\nu',J')\in\RC(\rho^t;\Rb^t)$ is in the image of 
$\jb$, then $\rkb(\nu',J')\ge \rkb(\nu,J)$. To prove the reverse and point 2
we introduce the set $\RC'(\la^{-t};\Rb^t)\subset\RC(\la^{-t};\Rb^t)$
as follows. For $\mu_L=1$ set $\RC'(\la^{-t};\Rb^t)=\RC(\la^{-t};\Rb^t)$.
For $\mu_L>1$, $(\nub,\Jb)\in\RC'(\la^{-t};\Rb^t)$ if $\rkb(\nub,\Jb)\ge c$
where $c$ is the column index of the cell $\la/\rho$ and $\rho\lessdot\la$
is the partition corresponding to $(\nub,\Jb)$. 

It will be shown that
\begin{equation*}
  \jbi:\RC'(\la^{-t};\Rb^t)\to\RC(\la^t;R^t)
\end{equation*}
is well-defined. Then set $(\nu,J)=\jbi(\nub,\Jb)$ for 
$(\nub,\Jb)\in\RC'(\la^{-t};\Rb^t)$.
Notice that the condition $\rkb(\nub,\Jb)\ge c$ implies that 
$\lb^{(k)}\le \lbi^{(k)}$ where $\lb^{(k)}$ and $\lbi^{(k)}$ are the lengths
of the selected strings in $(\nub,\Jb)$ under $\jb$ and $\jbi$, respectively.
In particular, $\lb^{(0)}\le \lb^{(1)}\le \lbi^{(1)}$ so that
$\lbi^{(0)}:=\mu_L-1\le \lbi^{(1)}$.
Using this one may easily verify that for all $k\ge 1$ and $n\ge 0$ the 
change in vacancy numbers under $\jbi$ is given by
\begin{equation}\label{dbi vacancy}
  P_n^{(k)}(\nu)-P_n^{(k)}(\nub)=\chi(\lbi^{(k-1)}<n\le \lbi^{(k)})
   -\chi(\lbi^{(k)}<n\le \lbi^{(k+1)}).
\end{equation}

It follows from the constructions of $\jb$ and $\jbi$
and \eqref{dbi vacancy} that $\jb\circ\jbi(\nub,\Jb)=(\nub,\Jb)$
for $(\nub,\Jb)\in\RC'(\la^{-t};\Rb^t)$.
This implies that the image of $\jb$ is given by $\RC'(\la^{-t};\Rb^t)$
proving point 3. Since for $\mu_L=1$, 
$\RC'(\la^{-t};\Rb^t)=\RC(\la^{-t};\Rb^t)$ it follows that
in this case $\jb$ is a bijection proving point 2.

We are left to prove that $\jbi$ is well-defined, that is, for every
$(\nub,\Jb)\in\RC'(\la^{-t};\Rb^t)$ the rigged configuration
$(\nu,J)=\jbi(\nub,\Jb)$ is admissible. This can be shown in a very 
similar fashion to the proof of the well-definedness of $\jb$.
Hence we only highlight the main arguments. By construction there are
no singular strings of length $n$ in $(\nub,\Jb)^{(k)}$
for $\lbi^{(k)}<n\le \lbi^{(k+1)}$. Hence by \eqref{dbi vacancy}
and \cite[Lemma 10]{KS} we need to show this time that 
for $\ell=\lbi^{(k)}$ the conditions 
\begin{equation}\label{pc}
m_{\ell+1}(\nub^{(k)})=0, \quad P_{\ell+1}^{(k)}(\nub)=0, \quad
\ell<\lbi^{(k+1)} \quad \text{and $\ell$ finite}
\end{equation}
cannot all be met simultaneously for all $k\ge 1$. Fix $k\ge 1$.
Let $\ell<p$ be minimal such that $m_p(\nub^{(k)})>0$; if no such
$p$ exists set $p=\infty$. 
By Lemma \ref{lower bound} the condition $P_{\ell+1}^{(k)}(\nub)=0$ 
implies that $P_n^{(k)}(\nub)=0$ for all $\ell\le n\le p$.
The condition that $\ell$ is finite requires $k<c$ where recall that
$c$ is the column index of the cell $\la/\rho$.

First assume that $k=c-1$. If $p$ is finite, i.e., $\ell$ is not the 
largest part in $\nub^{(k)}$ then there exists a singular string of 
length $p$ since $P_p^{(k)}(\nub)=0$ as argued above. But this means 
$\ell<\lbi^{(k)}$ which contradicts our assumptions. Hence assume that 
$\ell$ is the largest part in $\nub^{(k)}$. By \cite[(11.2)]{KS} 
one finds that $P_n^{(k)}(\nub)=\rho_k^t-\rho_{k+1}^t$ for large $n$. 
Since $P_n^{(k)}(\nub)=0$ for $n\ge \ell$ this requires 
$\rho_k^t=\rho_{k+1}^t$. Since $\la$ is a partition this implies that
$c\neq k+1$ which contradicts the assumption.
Hence \eqref{pc} cannot occur for $k=c-1$.

Now assume that $k<c-1$ which implies that $\lbi^{(k+1)}$ is finite.
Since $P_n^{(k)}(\nub)=0$ for $\ell\le n\le p$ one finds from
\eqref{vacancy ddiff} that $m_n(\nub^{(k+1)})=0$ for $\ell<n<p$ which 
means that $\nub^{(k+1)}$ does not contain any parts of size $n$.
When $\ell$ is not the largest part in $\nub^{(k)}$, that is, $p$ is finite,
we conclude that $\lbi^{(k+1)}\ge p$. However, since $P_p^{(k)}(\nub)=0$,
there exists a singular string of size $p$ in $(\nub,\Jb)^{(k)}$ and
hence $\lbi^{(k)}\ge p>\ell$. This contradicts $\lbi^{(k)}=\ell$.
When $\ell$ is the largest part in $\nub^{(k)}$, that is, $p$ is infinite,
we infer that the largest part in $\nub^{(k+1)}$ can be at most of size 
$\ell$. This implies $\lbi^{(k+1)}\leq \ell$ which contradicts the assumption
$\ell<\lbi^{(k+1)}$.
This concludes the proof that the conditions \eqref{pc} cannot occur and
shows that $\jbi$ is well-defined.
\end{proof}

Suppose the first rectangle of $R$ is a single column.
Define the set
\begin{equation*}
  \RC(\la^{-t};\Rt^t) = \bigcup_{\rho\lessdot\la} \RC(\rho^t;\Rt^t).
\end{equation*}
Define the map $\jt:\RC(\la^t;R^t)\rightarrow \RC(\la^{-t};\Rt^t)$
such that the diagram commutes:
\begin{equation} \label{jt def}
\begin{CD}
  \RC(\la^t;R^t) @>{\comev_R}>> \RC(\la^t;R^{\ev t}) \\
  @V{\jt}VV	@VV{\jb}V \\
  \RC(\la^{-t};\Rt^t) @>>{\comev_{\Rt}}> \RC(\la^{-t};\Rt^{\ev t}).
\end{CD}
\end{equation}
More precisely, for $(\nu,J)\in\RC(\la^t;R^{\ev t})$,
let $\jb(\nu,J)\in \RC(\rho^t;\Rt^{\ev t})$ for $\rho\lessdot\la$.
Then for $(\nu,J)\in\RC(\la^t;R^t)$, $\jt(\nu,J)=(\comev_{\Rt^\ev}\circ 
\jb \circ \comev_R)(\nu,J)\in \RC(\rho^t;\Rt^t)$ and define 
$\rkt(\nu,J)=\rkb(\comev_R(\nu,J))$. 
Observe that by definition and Proposition \ref{j bar} $\jt$ is an injection, 
with image given by $(\nu',J')\in\RC(\rho^t;\Rt^t)$ such that 
$\rkt(\nu',J') \ge c$ where $c$ is the column of the cell $\la/\rho$.
The map $\jt$ is given explicitly by the following algorithm.
Define $\lt^{(0)} = \mu_1$.  Inductively select a string
in $(\nu,J)^{(k)}$ with label $0$ and with length $\lt^{(k)}$
minimal such that $\lt^{(k-1)} \le \lt^{(k)}$.
Then $\rkt(\nu,J)$ is the minimum index $k$ for which such a string
does not exist; set $\lt^{(k)} = \infty$ for $k \ge \rkt(\nu,J)$.
Then $\jt(\nu,J)=(\nut,\Jt)$ is given by shortening the selected strings
by one and keeping their labels zero, and changing the
labels on all other strings such that their colabels do not change.
The computation for \eqref{db vacancy} yields
\begin{equation} \label{dt vacancy}
  P^{(k)}_n(\nu) - P^{(k)}_n(\nut) =
  \chi(\lt^{(k-1)} \le n < \lt^{(k)}) -
  \chi(\lt^{(k)} \le n < \lt^{(k+1)})
\end{equation}
for $k\ge1$ and $n\ge0$.

\begin{lem} \label{bar tilde}
Suppose the first and last rectangles of $R$ are single columns.  Then
$\jb$ and $\jt$ commute.  Moreover, for all $(\nu,J)\in \RC(\la^t;R^t)$
one of the following conditions holds, and the second can only hold if
$r:=\rkb(\nu,J)=\rkt(\nu,J)$ and $\la^t_r=\la^t_{r-1}$.
\begin{enumerate}
\item $\rkb(\jt(\nu,J))=\rkb(\nu,J)$ and $\rkt(\jb(\nu,J))=\rkt(\nu,J)$.
\item $\rkb(\jt(\nu,J))=\rkb(\nu,J)-1$ and $\rkt(\jb(\nu,J))=\rkt(\nu,J)-1$.
\end{enumerate}
\end{lem}
The proof of Lemma \ref{bar tilde} is rather technical
and is placed in Appendix \ref{App A}.  The statement regarding
$\rkb$ and $\rkt$ requires some explanation. Define
\begin{equation*}
  \RC(\la^{--t};\Rbt^t) = \bigcup_{\sigma\lessdot\rho\lessdot\la}
	  \RC(\sigma^t;\Rbt^t).
\end{equation*}
Similarly as for $\LRT(\la^{--};\Rbt)$, define
the involution $\tau$ on $\RC(\la^{--t};\Rbt^t)$ by
$(\nu,J)_\CC\mapsto (\nu,J)_{\tau(\CC)}$.  Then Lemma \ref{bar tilde}
says precisely that
\begin{equation} \label{j bar and tilde}
  \jt \circ \jb = \tau \circ \jb \circ \jt
\end{equation}
as maps $\RC(\la^t;R^t)\rightarrow \RC(\la^{--t};\Rbt^t)$.

For the reader's convenience the analogous maps on LR tableaux and rigged
configurations as defined in this Section 
and their main relations are listed in Table~\ref{table summarize}.
\begin{table}
\renewcommand{\arraystretch}{1.3}
\begin{tabular}{|l|c|c|}
\hline
Rectangles & LR tableaux & Rigged configurations \\ \hline
$R\to\Rht$   & $\iht$ & $\jht$ \\
$R\to\Rck$   & $\ick$ & $\jck$ \\
$R\to\Rb$    & $\ib$  & $\jb$ \\
$R\to\Rt$    & $D$    & $\jt$ \\
$R\to R^\ev$ & $\ev$  & $\comev_R$ \\
             & $\ib\circ D=\tau\circ D\circ\ib$ 
               & $\jb\circ\jt=\tau\circ\jt\circ\jb$ \\
             & $\ev\circ D=\ib\circ\ev$ 
               & $\comev_{\Rt}\circ\jt=\jb\circ\comev_R$ \\
             & $\ev\circ\iht=\ick\circ\ev$
               & $\comev_{\Rht}\circ\jht=\jck\circ\comev_R$ \\
\hline
\end{tabular}
\renewcommand{\arraystretch}{1}
\caption{\label{table summarize}
Maps defined in Section~\ref{sec operations} and some of their relations}
\end{table}

\section{The Bijection}
\label{bij sec}

\subsection{Definition}
We require two bijections $\phib_R$ and $\phit_R$ between
Littlewood--Richardson tableaux and rigged configurations.
The quantum number bijection $\phib_R:\LRT(\la;R)\rightarrow\RC(\la^t;R^t)$ 
is defined inductively based on Remark~\ref{bij ind} in 
Definition-Proposition~\ref{phi def} below. Recall
that $\com_R$ complements the quantum numbers of a rigged configuration
(see~\eqref{com def}). The coquantum number bijection 
$\phit_R:\LRT(\la;R)\rightarrow\RC(\la^t;R^t)$ is defined as
\begin{equation}\label{phi coquantum}
  \phit_R:=\com_R\circ\phib_R.
\end{equation}
It is $\phit_R$ that preserves the statistics (see Theorem~\ref{thm charge}).

\begin{defprop} \label{phi def} There is a unique family of bijections
$\phib_R:\LRT(\la;R)\rightarrow \RC(\la^t;R^t)$ indexed by $R$, such that:
\begin{enumerate}
\item If the last rectangle of $R$ is a single column, then
the following diagram commutes:
\begin{equation} \label{bij bar}
\begin{CD}
  \LRT(\la;R) @>{\ib}>> \LRT(\la^-;\Rb) \\
  @V{\phib_R}VV	@VV{\phib_{\Rb}}V \\
  \RC(\la^t;R^t) @>>{\jb}> \RC(\la^{-t};\Rb^t).
\end{CD}
\end{equation}
\item The following diagram commutes:
\begin{equation} \label{bij hat}
\begin{CD}
  \LRT(\la;R) @>{\iht}>> \LRT(\la;\Rht) \\
  @V{\phib_R}VV @VV{\phib_{\Rht}}V \\
  \RC(\la^t;R^t) @>>{\jht}> \RC(\la^t;{\Rht}^t).
\end{CD}
\end{equation}
\end{enumerate}
\end{defprop}

\begin{proof}
The proof proceeds by the induction on $R$ given by
Remark \ref{bij ind}. Suppose first that $R$ is the empty sequence.
Then both $\LRT(\la;R)$ and $\RC(\la^t;R^t)$ are the empty set
unless $\la$ is the empty partition, in which case
$\LRT(\la;R)$ is the singleton consisting of the empty tableau,
$\RC(\la^t;R^t)$ is the singleton consisting of the empty rigged
configuration, and $\phib_\emptyset$ is the unique bijection
$\LRT(\emptyset;\emptyset)\rightarrow \RC(\emptyset;\emptyset)$.

Suppose that the last rectangle of $R$ is a single column.
Then 2 holds trivially since $\Rht=R$ and $\iht$ and $\jht$ are the
identity maps. Consider 1.
By induction, for every partition $\rho$ such that
$\rho\lessdot \la$ the result holds for the pair $(\rho;\Rb)$.
Any map $\phib_R$ satisfying \eqref{bij bar} is injective by
definition and unique by induction.  For the existence and surjectivity
of $\phib_R$ it suffices to show that the bijection
$\phib_{\Rb}:\LRT(\la^-;\Rb)\rightarrow\RC(\la^{-t};\Rb^t)$
maps the image of $\ib:\LRT(\la;R)\rightarrow\LRT(\la^-;\Rb)$ onto
the image of $\jb:\RC(\la^t;R^t)\rightarrow\RC(\la^{-t};\Rb^t)$.
If $\mu_L=1$ then $-$ and $\jb$ are bijections by 
Propositions~\ref{LR minus inj}
and \ref{j bar} respectively proving the assertion. Hence assume $\mu_L>1$.
Let $c$ and $c'$ be the column indices of the cells $\la/\rho$ and 
$\rho/\sigma$ for the partitions $\sigma\lessdot\rho\lessdot\la$,
respectively. By Proposition \ref{LR minus inj}, $T\in\LRT(\rho;\Rb)$ 
with $\shape(T^-)=\sigma$ is in the image of 
$-:\LRT(\la;R)\to \LRT(\la^-;\Rb)$ if and only if $c\le c'$.
Similarly by Proposition~\ref{j bar}, $(\nu,J)\in\RC(\rho^t;\Rb^t)$ with
$\rkb(\nu,J)=c'$ is in the image of $\jb:\RC(\la^t;R^t)\to\RC(\la^{-t};\Rb^t)$
if and only if $c\le c'$. This proves the assertion about the images
of $-$ and $\jb$.

Suppose the last rectangle of $R$ has more than one column, that is
$\eta_L>1$. Any map $\phib_R$ satisfying 2 is injective by definition
and unique by induction.  For existence and surjectivity
it is enough to show that the bijection
$\phib_{\Rht}$ maps the image $\LRTh(\la;\Rht)$
of $\iht:\LRT(\la;R)\rightarrow\LRT(\la;\Rht)$
onto the image $\RCh(\la^t;{\Rht}^t)$
of $\jht:\RC(\la^t;R^t)\rightarrow\RC(\la^t;{\Rht}^t)$.
A rigged configuration $(\nu,J)$ is in $\RCh(\la^t;{\Rht}^t)$ if and
only if $(\nu,J)\in\RC(\la^t;{\Rht}^t)$ and $(\nu,J)^{(k)}$ contains
a singular string of length $\mu_L$ for all $1\le k<\eta_L$.
An LR tableau $T$ is in $\LRTh(\la;\Rht)$ if and only if
$T\in\LRT(\la;\Rht)$ and the column index $c_i$ of the cell
$\shape(T_i)/\shape(T_{i-1})$ and the column index 
$\cb_i$ of the cell $\shape(\Tb_i)/\shape(\Tb_{i-1})$ 
where $T_i=T^{-^{\mu_L-i}}$ and $\Tb_i=(\iht(T_0))^{-^{\mu_L-i}}$ 
obey $\cb_i<c_i$ for all $1\le i\le \mu_L$.
This follows from the definition of $\LRT$ and 
Remark~\ref{funny compatible}.

Set $(\nu_i,J_i)=\jb^{\mu_L-i}(\nu,J)$ for $0\le i\le \mu_L$, 
$c_i=\rkb(\nu_i,J_i)$ and denote the length of the 
string in $(\nu_i,J_i)^{(k)}$ selected by $\jb$ by $\ell_i^{(k)}$.
Similarly set $(\nub,\Jb)=\jht(\nu_0,J_0)$ and
$(\nub_i,\Jb_i)=\jb^{\mu_L-i}(\nub,\Jb)$ for $0\le i\le \mu_L$, 
define $\cb_i=\rkb(\nub_i,\Jb_i)$ and denote the length of the string in 
$(\nub_i,\Jb_i)^{(k)}$ selected by $\jb$ by $\lb_i^{(k)}$.
Hence to prove that the bijection $\phib_{\Rht}$ maps the image of $\iht$
onto the image of $\jht$ one needs to show that
\begin{quote}
\textbf{Claim.}\\
For $(\nu,J)\in\RC(\la^t;{\Rht}^t)$, $\ell^{(k)}_{\mu_L}=\mu_L$ for
$0\le k<\eta_L$ if and only if $\cb_i<c_i$ for $1\le i\le \mu_L$.
\end{quote}

We begin by showing that $\cb_i<c_i$ if $\ell_{\mu_L}^{(k)}=\mu_L$ for 
$0\le k<\eta_L$. First of all notice that 
\begin{equation}\label{ell restr}
\ell_i^{(k)}=i \quad \text{for $0\leq k<\eta_L$ and} \quad 
\lb_i^{(k)}=i \quad \text{for $0\leq k<\eta_L-1$}.
\end{equation}
For $i=\mu_L$ this equation holds since both $(\nu,J)=(\nu_{\mu_L},J_{\mu_L})$
and $(\nub,\Jb)=(\nub_{\mu_L},\Jb_{\mu_L})$ are in the image of $\jht$.
Now assume \eqref{ell restr} to be true at $i$.
This means that there are singular strings of length $i$
in the first $\eta_L-1$ ($\eta_L-2$) partitions of $\nu_i$ 
($\nub_i$). Hence by construction these turn into singular strings of 
length $i-1$. Since by definition $\ell_{i-1}^{(0)}=\lb_{i-1}^{(0)}=i-1$ 
this implies \eqref{ell restr} at $i-1$. We claim that for $k\ge 1$ 
and $0\le i\le \mu_L$
\begin{equation}\label{Pki}
  P_n^{(k)}(\nub_i)\geq J_n^{(k)}(\nub_i,\Jb_i)
   +\sum_{m=1}^i\chi(\ell_m^{(k)}\leq n<\ell_m^{(k+1)})
\end{equation}
where $J_n^{(k)}(\nu,J)$ is the maximal label occurring in the 
strings of length $n$ in $(\nu,J)^{(k)}$ and $J_n^{(k)}(\nu,J)=0$ if there is
no string of length $n$ in $\nu^{(k)}$.
Equations \eqref{ell restr} and \eqref{Pki} imply
\begin{equation}\label{l}
  \ell_i^{(k+1)} \le \lb_i^{(k)} \qquad \text{for all $k\ge 0$}.
\end{equation}
This can be shown by induction on $k$.
The initial condition follows immediately from \eqref{ell restr}. By
construction $\lb_i^{(k-1)}\leq \lb_i^{(k)}$ which becomes by
induction hypothesis $\ell_i^{(k)}\leq \lb_i^{(k)}$.
Inequality \eqref{Pki} implies that there are no singular strings of
length $n$ in $(\nub_i,\Jb_i)^{(k)}$ for $\ell_i^{(k)}\le n<\ell_i^{(k+1)}$
which proves \eqref{l}.
The condition \eqref{l} immediately implies that $\cb_i<c_i$.

It remains to show \eqref{Pki}, whose proof proceeds by
descending induction on $i$, with base case $i=\mu_L$.
To establish the base case it is shown by induction on 
$j$ ($0\le j\le \mu_L$) that
\begin{equation}\label{P ind}
  P_n^{(k)}(\nu_j)\geq J_n^{(k)}(\nu_j,J_j)+
   \sum_{m=j+1}^{\mu_L} \chi(\ell_m^{(k)}\le n<\ell_m^{(k+1)}) 
   \quad \text{for all $k\geq 1$.}
\end{equation}
Since $(\nu,J)=(\nu_{\mu_L},J_{\mu_L})$ is an admissible rigged configuration, 
$P_n^{(k)}(\nu)\ge J_n^{(k)}(\nu,J)$ which implies \eqref{P ind} at $j=\mu_L$.
Now we assume \eqref{P ind} to be true at $j+1$ and show its validity at $j$.
Since by Proposition \ref{j bar} all $(\nu_j,J_j)$ are admissible, 
$P_n^{(k)}(\nu_j)\ge J_n^{(k)}(\nu_j,J_j)$. This settles \eqref{P ind} for 
$n<\ell_{j+1}^{(k)}$ because in this case the sum over $m$ vanishes 
since by construction 
\begin{equation}\label{ell ineq}
  \ell_1^{(k)}<\ell_2^{(k)}<\cdots<\ell_{\mu_L}^{(k)}.
\end{equation}
Let us now consider the case $n\ge \ell_{j+1}^{(k)}$. The only string that 
changes in the transformation $(\nu_{j+1},J_{j+1})^{(k)}\to (\nu_j,J_j)^{(k)}$
is one singular string of length $\ell_{j+1}^{(k)}$. The riggings of all other 
strings remain unchanged. In particular, 
$J_n^{(k)}(\nu_j,J_j)=J_n^{(k)}(\nu_{j+1},J_{j+1})$ for $n>\ell_{j+1}^{(k)}$
and $J_n^{(k)}(\nu_j,J_j)\le J_n^{(k)}(\nu_{j+1},J_{j+1})$ for 
$n=\ell_{j+1}^{(k)}$. Hence for $n\ge \ell_{j+1}^{(k)}$ we find
\begin{align*}
  P_n^{(k)}(\nu_j)&=P_n^{(k)}(\nu_{j+1})+
   \chi(\ell_{j+1}^{(k)}\le n< \ell_{j+1}^{(k+1)})\\
  &\ge J_n^{(k)}(\nu_{j+1},J_{j+1})+\sum_{m=j+1}^{\mu_L} 
   \chi(\ell_m^{(k)}\le n< \ell_m^{(k+1)})\\
  &\ge J_n^{(k)}(\nu_j,J_j)+\sum_{m=j+1}^{\mu_L} 
   \chi(\ell_m^{(k)}\le n< \ell_m^{(k+1)})
\end{align*}
where the first line follows from \eqref{db vacancy} and 
$n\ge \ell_{j+1}^{(k)}$, and in the second line the induction 
hypothesis is used. This concludes the inductive 
proof of \eqref{P ind}.

Since by definition $(\nub,\Jb)=\jht(\nu_0,J_0)$ and $\jht$ preserves the
vacancy numbers by Lemma \ref{j hat}, \eqref{P ind} at $j=0$ implies 
\eqref{Pki} at $i=\mu_L$.
Now assume that \eqref{Pki} holds at $i$.
Because of the admissibility of $(\nub_{i-1},\Jb_{i-1})$ and
\eqref{ell ineq} at $k+1$, \eqref{Pki} holds at $i-1$ for 
$n\ge \ell_{i-1}^{(k+1)}$. Hence assume $n<\ell_{i-1}^{(k+1)}$. Since 
$\ell_{i-1}^{(k+1)}<\ell_i^{(k+1)}\le \lb_i^{(k)}$ by \eqref{ell ineq}
and \eqref{l}, in particular $n\le \lb_i^{(k)}-2$. Hence we find for
$n<\ell_{i-1}^{(k+1)}$
\begin{align*}
  P_n^{(k)}(\nub_{i-1})&=P_n^{(k)}(\nub_i)-
   \chi(\lb_i^{(k-1)}\le n< \lb_i^{(k)})\\
  &\ge J_n^{(k)}(\nub_i,\Jb_i)+\sum_{m=1}^{i-1} 
   \chi(\ell_m^{(k)}\le n< \ell_m^{(k+1)})
\end{align*}
where the first line follows from \eqref{db vacancy} and $n\le \lb_i^{(k)}-2$,
and in the second line the induction hypothesis and \eqref{l} are used.
Since $J_n^{(k)}(\nub_{i-1},\Jb_{i-1})=J_n^{(k)}(\nub_i,\Jb_i)$ for 
$n\le \lb_i^{(k)}-2$, this yields \eqref{Pki} at $i-1$.
This concludes the proof of \eqref{Pki} and also that of the forward 
direction of the claim.

Next we prove the reverse direction of the claim. More precisely,
we show that
\begin{equation}\label{lk}
  \ell_i^{(k)}=i \qquad \text{for $0\leq k<\eta_L$ and $1\leq i\leq \mu_L$}
\end{equation}
if $\cb_i<c_i$.
Since $(\nub,\Jb)=(\nub_{\mu_L},\Jb_{\mu_L})$ is in the image of $\jht$
it follows by construction that
\begin{equation}\label{lkb}
  \lb_i^{(k)}=i \qquad \text{for $0\leq k<\eta_L-1$ and $1\leq i\leq \mu_L$.}
\end{equation}
We claim that for $k\ge 1$ and $1\le i\le \mu_L+1$
\begin{equation}\label{Pk}
  P_n^{(k)}(\nu_{i-1})\ge J_n^{(k)}(\nu_{i-1},J_{i-1})+
   \sum_{m=i}^{\mu_L} \chi(\lb_m^{(k-1)}\le n<\lb_m^{(k)}).
\end{equation}
The condition $\cb_i<c_i$ and \eqref{Pk} imply by induction on $k$ that
\begin{equation}\label{ll}
  \ell_i^{(k)}\leq \lb_i^{(k-1)} \qquad 
   \text{for all $1\le i\le \mu_L$ and $k\geq 1$.}
\end{equation}
Before proving this notice the following. Let $\jbi$ be the inverse
algorithm of $\jb$ as introduced in the proof of Proposition \ref{j bar}
and let $s_i^{(k)}$ be the length of the singular string in 
$(\nu_{i-1},J_{i-1})^{(k)}$ selected by $\jbi$. Note that 
$\jbi\circ\jb(\nu_i,J_i)=(\nu_i,J_i)$ and $s_i^{(k)}=\ell_i^{(k)}-1$. 
If $(\nu_{i-1},J_{i-1})^{(k)}$ does not contain singular strings of length
$n$ for $a\le n<b$ and $s_i^{(k)}<b$ then by construction 
$s_i^{(k)}<a$ and by $s_i^{(k)}=\ell_i^{(k)}-1$ also $\ell_i^{(k)}\le a$.
Now we prove \eqref{ll} by induction on $k$.
For $k>\cb_i$ equation \eqref{ll} is true because in this case 
$\lb_i^{(k-1)}=\infty$. Now consider $k=\cb_i$. Then \eqref{Pk} implies that
there are no singular strings of length $n\ge\lb_i^{(\cb_i-1)}$
in $(\nu_{i-1},J_{i-1})^{(\cb_i)}$ since $\lb_i^{(\cb_i)}=\infty$. 
Since $\cb_i<c_i$ the variable $\ell_i^{(\cb_i)}$ is finite. This implies that 
$\ell_i^{(\cb_i)}\le \lb_i^{(\cb_i-1)}$ which is \eqref{ll} for $k=\cb_i$.
Now assume \eqref{ll} to be true at $k+1$. By construction 
$\ell_i^{(k)}\le \ell_i^{(k+1)}$ which implies by induction hypothesis that 
$\ell_i^{(k)}\le \lb_i^{(k)}$. Because of \eqref{Pk} there is no singular 
string of length $n$ in $(\nu_{i-1},J_{i-1})^{(k)}$ with
$\lb_i^{(k-1)}\le n<\lb_i^{(k)}$. Hence $\ell_i^{(k)}$ has to obey
condition \eqref{ll}. Together with $\ell_i^{(0)}=i$, \eqref{ll} and 
\eqref{lkb} immediately imply \eqref{lk}.

We are left to prove \eqref{Pk}. Using $\jbi$ one proves in an
analogous fashion to \eqref{P ind} that
\begin{equation*}
  P_n^{(k)}(\nub_j)\ge J_n^{(k)}(\nub_j,\Jb_j)
   +\sum_{m=1}^j \chi(\lb_m^{(k-1)}\le n<\lb_m^{(k)})
\end{equation*}
for all $k\ge 1$ and $0\le j\le \mu_L$. Since 
$(\nub_{\mu_L},\Jb_{\mu_L})=(\nub,\Jb)=\jht(\nu_0,J_0)$ and $\jht$
preserves the vacancy numbers by Lemma \ref{j hat}, this inequality
at $j=\mu_L$ implies \eqref{Pk} at $i=1$.
Now assume \eqref{Pk} to be true at $i-1$. 
Using \eqref{db vacancy} and $\ell_{i-1}^{(k)}\le \lb_{i-1}^{(k-1)}$
equation \eqref{Pk} at $i$ follows by similar arguments to those used
to prove \eqref{Pki}.
\end{proof}

\subsection{Direct algorithm for the bijection}
\label{sec algo}

Here we state the original algorithm for the bijection $\phib_R$
between Littlewood--Richardson tableaux and rigged configurations 
as given in~\cite{K}. It combines points 1 and 2 of the 
Definition-Proposition~\ref{phi def} which has the advantage that
it is computationally simpler. For the proofs the formulation
of Definition-Proposition~\ref{phi def} is however more convenient.

Let $T\in\LRT(\la;R)$ be an LR tableaux for a partition $\la$ and
a sequence of rectangles $R=(R_1,\ldots,R_L)$. Set $N=|\la|$
and $B_j=[|R_1|+\cdots+|R_{j-1}|+1,|R_1|+\cdots+|R_j|]$.
To obtain $(\nu,J)=\phib_R(T)$ one recursively constructs a rigged 
configuration $(\nu,J)_{(x)}$ for each letter $1\le x\le N$ occurring in
$T$. Set $(\nu,J)_{(0)}=\emptyset$. Suppose that $x\in B_j$, and
denote the column index of $x$ in $T$ by $c$ and the column index 
of $x$ in $\ZC_j$ by $c'$. Define the numbers $s^{(k)}$ for
$c'\le k<c$ as follows. Let $s^{(c-1)}$ be the length of the longest
singular string in $(\nu,J)_{(x-1)}^{(c-1)}$. Now select inductively
a singular string in $(\nu,J)_{(x-1)}^{(k)}$ for $k=c-2,c-3,\ldots,c'$ whose 
length $s^{(k)}$ is maximal such that $s^{(k)}\le s^{(k+1)}$; if no
such string exists set $s^{(k)}=0$.
Then $(\nu,J)_{(x)}$ is obtained from $(\nu,J)_{(x-1)}$ by adding one
box to the selected strings with labels such that they remain
singular, leaving all other strings unchanged.
Then the image of $T$ under $\phib_R$ is given by $(\nu,J)=(\nu,J)_{(N)}$.

For the above algorithm it is necessary to be able to compute
the vacancy numbers of an intermediate configuration $\nu_{(x)}$.
Suppose $x$ occurs in $\ZC_j$ in column $c'$. In general 
$R_{(x)}=(R_1,\ldots,R_{j-1},\shape(\ZC_j|_{[1,x]}))$ is not a sequence
of rectangles. If $\shape(\ZC_j|_{[1,x]})$ is not a rectangle
one splits it into two rectangles, one of width $c'$ and 
one of width $c'-1$. The vacancy numbers are calculated with respect to 
this new sequence of rectangles.

\begin{ex}\label{ex bij}
The non-trivial steps of the above algorithm applied to the third 
tableau of Example~\ref{ex LR} are given in Table~\ref{table ex}.
A rigged partition is represented by its Ferrers diagram where
to the right of each part the corresponding rigging is indicated.
The vacancy numbers are given to the left of each part.
For example $R_{(12)}=((3,3),(2,2,1,1))$ so that the vacancy numbers
of $(\nu,J)_{(12)}$ are calculated with respect to the sequence of 
rectangles $((3,3),(2,2),(1,1))$.
\begin{table}
\begin{center}
\begin{equation*}
\begin{array}{|r|llll|} \hline &&&&\\[-3mm] 
 x & (\nu,J)_{(x)}^{(1)} & (\nu,J)_{(x)}^{(2)} & (\nu,J)_{(x)}^{(3)} 
   & (\nu,J)_{(x)}^{(4)}
 \\[1mm] \hline &&&&\\
 11 &
    & \begin{array}{r|c|l} \cline{2-2} 0& &0 \\ \cline{2-2} \end{array}
    & \begin{array}{r|c|l} \cline{2-2} 0& &0 \\ \cline{2-2} \end{array}
    & \\ &&&&\\
 12 &
    & \begin{array}{r|c|c|l} \cline{2-3} 0& & &0 \\ \cline{2-3} \end{array}
    & \begin{array}{r|c|c|l} \cline{2-3} 0& & &0 \\ \cline{2-3} \end{array} 
    & \\ &&&&\\
 14 &
    & \begin{array}{r|c|c|l} \cline{2-3} 0& & &0 \\ \cline{2-3} \end{array}
    & \begin{array}{r|c|c|l} \cline{2-3} 0& & &0 \\ \cline{2-3} \end{array} 
    & \\ &&&&\\
 15 & \begin{array}{r|c|l} \cline{2-2} 1& &1 \\ \cline{2-2} \end{array}
    & \begin{array}{r|c|c|l} \cline{2-3} 0&&&0\\ \cline{2-3}
      0&&\multicolumn{2}{l}{0}\\ \cline{2-2} \end{array}
    & \begin{array}{r|c|c|l} \cline{2-3} 0&&&0\\ \cline{2-3}
      0&&\multicolumn{2}{l}{0}\\ \cline{2-2} \end{array}
    & \begin{array}{r|c|l} \cline{2-2} 0& &0 \\ \cline{2-2} \end{array} \\ 
    &&&&\\
 16 & \begin{array}{r|c|c|l} \cline{2-3} 1& & &1 \\ \cline{2-3} \end{array}
    & \begin{array}{r|c|c|c|l} \cline{2-4} 0&&&&0\\ \cline{2-4}
      0&&\multicolumn{3}{l}{0}\\ \cline{2-2} \end{array}
    & \begin{array}{r|c|c|l} \cline{2-3} 0&&&0\\ \cline{2-3}
      0&&\multicolumn{2}{l}{0}\\ \cline{2-2} \end{array}
    & \begin{array}{r|c|l} \cline{2-2} 0& &0 \\ \cline{2-2} \end{array}\\ 
    &&&&\\
 17 & \begin{array}{r|c|c|c|l} \cline{2-4} 1&&&&1 \\ \cline{2-4} \end{array}
    & \begin{array}{r|c|c|c|l} \cline{2-4} 1&&&&0\\ \cline{2-4}
      0&&\multicolumn{3}{l}{0}\\ \cline{2-2} \end{array}
    & \begin{array}{r|c|c|l} \cline{2-3} 0&&&0\\ \cline{2-3}
      0&&\multicolumn{2}{l}{0}\\ \cline{2-2} \end{array}
    & \begin{array}{r|c|l} \cline{2-2} 0& &0 \\ \cline{2-2} \end{array}\\ 
    &&&& \\ \hline
\end{array}
\end{equation*}
\end{center}
\caption{\label{table ex}
 Example for the bijection algorithm (see Example~\ref{ex bij})}
\end{table}
\end{ex}

It is relatively straightforward to see that the above described
algorithm is indeed an algorithm for the bijection of 
Definition-Proposition~\ref{phi def}. By induction it suffices
to study the effect of the last $\mu_L$ letters of the LR tableaux $T$,
that is $N-\mu_L<x\le N$.
Recall that $\jht$, which corresponds to splitting off the last column 
of $R_L$, adds a singular string of length $\mu_L$ to each of the first 
$\eta_L-1$ rigged partitions and leaves the vacancy numbers invariant. 
By construction, these extra singular strings are removed by $\jb^{\mu_L}$. 
This has the same effect as restricting the removal/addition of boxes
to the partitions $\nu^{(k)}_{(x)}$ with $\eta_L=c'\le k$ as is the case for
the algorithm of this section. One may also show that the vacancy numbers 
of the intermediate rigged partitions obtained from the the algorithm 
and the recursive definition of Definition-Proposition~\ref{phi def} 
are the same for $\nu^{(k)}_{(x)}$ with $k\neq \eta_L-1$. For 
$\nu^{(\eta_L-1)}_{(x)}$ the vacancy numbers from the algorithm cannot 
be smaller than the corresponding ones coming from 
Definition-Proposition~\ref{phi def}, and are the same when $R_{(x)}$ is 
a sequence of rectangles. This difference for non-rectangular $R_{(x)}$ 
is harmless since no strings in $\nu^{(\eta_L-1)}_{(x)}$ are changed by 
the algorithm.

\section{Evacuation Theorem}
\label{sec ev}

In this Section we prove the Evacuation Theorem~\ref{ev} which
relates the evacuation of LR tableaux to the complementation
of quantum numbers on rigged configurations. The proof requires
intertwining relations of $\phib_R$ with $R\to\Rck$ and
$R\to\Rt$ which are derived in Section~\ref{sec intertwine}.

\subsection{Intertwining of $\phib_R$ with $R\to\Rck$ and $R\to\Rt$}
\label{sec intertwine}

\begin{lem} \label{j check hat} $\jht \circ \jck = \jck \circ \jht$.
\end{lem}
\begin{proof} 
If $R$ consists of more than one rectangle or $R$ is a single rectangle
with more than two columns, the commutativity of
$\jht$ and $\jck$ is obvious. 
If $R$ is a single column then both $\jht$ and $\jck$
are the identity and obviously commute.  
So it may be assumed that $R$ is a single rectangle with exactly two columns.
Then $\Rht=\Rck=((1^{\mu_1}),(1^{\mu_1}))$.  This means that the
outer function in both $\jht\circ \jck$ and $\jck\circ \jht$
acts as the identity.  So ${\Rht}\ck = \Rht$, ${\Rck}\htt=\Rck$, and
for all $(\nu,J)\in \RC(\la^t;R^t)$, 
$\jht(\jck(\nu,J)) = \jck(\nu,J)$ and
$\jck(\jht(\nu,J)) = \jht(\nu,J)$.  So it must be shown that
$\jck(\nu,J)=\jht(\nu,J)$.  Since $\jht$ (resp. $\jck$)
adds a string of length $\mu_1$ with singular (resp. zero) label
to the first rigged partition of $(\nu,J)$, it must be shown that
$P^{(1)}_{\mu_1}(\nu) = 0$ for all $(\nu,J) \in \RC(\la^t;R^t)$.
It may be assumed that $\la=R_1$ for otherwise
$\RC(\la^t;R^t)$ is empty.  Then $\RC(\la^t;R^t)$ is the singleton set
consisting of the empty rigged configuration
$(\emptyset,\emptyset)$.  One computes the vacancy number
$P^{(1)}_{\mu_1}(\emptyset) = 0$.
\end{proof}

\begin{lem} \label{bar check}
Suppose the last rectangle of $R$ is a single column.  Then
$\jb \circ \jck = \jck \circ \jb$.
\end{lem}
\begin{proof} Let $\lcb^{(k)}$ be the lengths of strings selected
by $\db$ acting on $\jck(\nu,J)=(\nuck,\Jck)$.  
To prove the lemma it suffices to show that $\lcb^{(k)}=\lb^{(k)}$ for 
all $k\ge 1$.
By Lemma \ref{j check}, $P^{(k)}_n(\nu)=P^{(k)}_n(\nuck)$ for all $k,n \ge 1$.
Since $(\nuck,\Jck)^{(k)}$ is obtained from $(\nu,J)^{(k)}$ by
adding the string $(\mu_1,0)$ for $1\le k \le \eta_1-1$,
it is clear that $\lcb^{(k)} \le \lb^{(k)}$ for all $k$.
Let $k$ be minimal such that $\lcb^{(k)} < \lb^{(k)}$.
Then the string $(\mu_1,0)$ was selected in $(\nuck,\Jck)^{(k)}$,
so that $\lcb^{(k)} = \mu_1$ and
$P^{(k)}_{\mu_1}(\nuck)=P^{(k)}_{\mu_1}(\nu)=0$.  Now
$\lb^{(k-1)}=\lcb^{(k-1)}\le \lcb^{(k)}=\mu_1$ and
$\mu_1=\lcb^{(k)} < \lb^{(k)}$. By \eqref{db vacancy} we have
$P^{(k)}_{\mu_1}(\nub) = P^{(k)}_{\mu_1}(\nu) -
\chi(\lb^{(k-1)} \le \mu_1 < \lb^{(k)}) +
\chi(\lb^{(k)} \le \mu_1 < \lb^{(k+1)}) = -1$, which is a contradiction.
Thus there is no such $k$.
\end{proof}

For future use let us formalize a general nonsense lemma.
\begin{lem} \label{nonsense} Suppose we have a diagram of the following
kind:
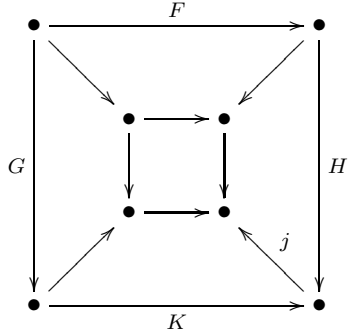
\begin{equation*}
\xymatrix{
 {\bull} \ar[rrr]^{F} \ar[ddd]_{G} \ar[dr] & & &
	{\bull} \ar[ddd]^{H} \ar[dl] \\
 & {\bull} \ar[r] \ar[d] & {\bull} \ar[d] & \\
 & {\bull} \ar[r]   & {\bull}  & \\
 {\bull} \ar[rrr]_{K} \ar[ur] & & & {\bull} \ar[ul]_{j}
}
\end{equation*}
Viewing this diagram as a cube with front face given by the
large square, suppose the square diagrams given by all the faces
of the cube except the front, commute.  Assume also that
the map $j$ is injective.  Then the front face must also commute.
\end{lem}
\begin{proof} The commuting faces yield the equality
$j \circ K \circ G = j \circ H \circ F$.
Since $j$ is injective this implies $K \circ G = H \circ F$ as desired.
\end{proof}

The following Lemma gives the intertwining of $\phib_R$ with
$R\to\Rck$.
\begin{lem} 
The following diagram commutes:
\begin{equation} \label{bij check}
\begin{CD}
  \LRT(\la;R) @>{\ick}>> \LRT(\la;\Rck) \\
   @V{\phib_R}VV	@VV{\phib_{\Rck}}V \\
  \RC(\la^t;R^t) @>>{\jck}> \RC(\la^t;(\Rck)^t).
\end{CD}
\end{equation}
\end{lem}
\begin{proof} The proof proceeds by the induction
that defines the bijection $\phib$.

Suppose the last rectangle of $R$ has more than one column;
this subsumes the base case $R=(R_1)$.
Clearly ${\Rck}\htt = {\Rht}\ck$.  Consider the diagram
\begin{equation*}
\xymatrix{
{\LRT(\la;R)} \ar[rrr]^{\ick} \ar[ddd]_{\phib_R} \ar[dr]^{\iht}
& & & {\LRT(\la;\Rck)} \ar[ddd]^{\phib_{\Rck}} \ar[dl]_{\iht} \\
 & {\LRT(\la;\Rht)} \ar[r]^{\ick} \ar[d]_{\phib_{\Rht}} &
	 {\LRT(\la;\Rhtck)} \ar[d]^{\phib_{\Rhtck}} & \\
 & {\RC(\la^t;{\Rht}^t)} \ar[r]_{\jck} & {\RC(\la^t;(\Rhtck)^t)} & \\
{\RC(\la^t;R^t)} \ar[ur]_{\jht} \ar[rrr]_{\jck} & & & {\RC(\la^t;(\Rck)^t)}
	\ar[ul]^{\jht}
}
\end{equation*}
We wish to show the front face commutes.  By Lemma \ref{nonsense}
and the injectivity of $\jht$,
it is enough to show that all the other faces commute.
The back face is assumed to commute by induction,
the left and right faces commute by the definition of the
bijections $\phib$ (see \eqref{bij hat}), the commutation
of the top face is obvious, and the commutativity of the bottom
face follows from Lemma~\ref{j check hat}.

For the remaining case suppose the last rectangle of $R$
is a single column.  Clearly $\Rb\ck=\overline{\Rck}$.  Consider the diagram
\begin{equation*}
\xymatrix{
{\LRT(\la;R)} \ar[rrr]^{\ick} \ar[ddd]_{\phib_R} \ar[dr]^{\ib}
& & & {\LRT(\la;\Rck)} \ar[ddd]^{\phib_{\Rck}} \ar[dl]_{\ib} \\
 & {\LRT(\la^-;\Rb)} \ar[r]^{\ick} \ar[d]_{\phib_{\Rb}} &
 {\LRT(\la^-;\Rb\ck)} \ar[d]^{\phib_{\Rb\ck}} & \\
 & {\RC(\la^{-t};\Rb^t)} \ar[r]_{\jck} & {\RC(\la^{-t};(\Rb\ck)^t)} & \\
{\RC(\la^t;R^t)} \ar[ur]_{\jb} \ar[rrr]_{\jck} & & & {\RC(\la^t;(\Rck)^t)}
	\ar[ul]^{\jb}
}
\end{equation*}
$\jb$ is injective, so again by Lemma \ref{nonsense} it suffices
to check that all faces but the front, commute.
The back face commutes by induction.  The left and right faces commute
by the definition of $\phib$.  The top face obviously commutes.
The bottom face commutes by Lemma \ref{bar check}.
\end{proof}

The intertwining relation of $\phib_R$ with $R\to\Rt$ is stated
in the next Lemma.
\begin{lem} Suppose the first rectangle of $R$ is a single column.
Then the following diagram commutes:
\begin{equation} \label{bij tilde}
\begin{CD}
  \LRT(\la;R) @>{\iit}>> \LRT(\la^-;\Rt) \\
   @V{\phib_R}VV	@VV{\phib_{\Rt}}V \\
  \RC(\la^t;R^t) @>>{\jt}> \RC(\la^{-t};\Rt^t).
\end{CD}
\end{equation}
\end{lem}
\begin{proof} Suppose that the last rectangle of $R$ is also a single
column; this case subsumes the base case that $R=(R_1)$ is a single
column.  Clearly $\Rbt=\Rtb$.
\begin{equation*}
\xymatrix{
{\LRT(\la;R)} \ar[rrr]^{\iit} \ar[ddd]_{\phib_R} \ar[dr]^{\ib}
& & & {\LRT(\la^-;\Rt)} \ar[ddd]^{\phib_{\Rt}} \ar[dl]_{\ib} \\
 & {\LRT(\la^-;\Rb)} \ar[r]^{\iit} \ar[d]_{\phib_{\Rb}} &
 {\LRT(\la^{--};\Rbt)} \ar[d]^{\phib_{\Rbt}} & \\
 & {\RC(\la^{-t};\Rb^t)} \ar[r]_{\jt} & {\RC(\la^{--t};\Rbt^t)} & \\
{\RC(\la^t;R^t)} \ar[ur]_{\jb} \ar[rrr]_{\jt} & & & {\RC(\la^{-t};\Rt^t)}
	\ar[ul]^{\jb}
}
\end{equation*}
In this diagram there is a map
\begin{equation}
  \phib_{\Rbt}: \LRT(\la^{--};\Rbt)\rightarrow \RC(\la^{--t};\Rbt^t).
\end{equation}
This is to be understood in the most obvious way, namely, that
given a chain of partitions $\CC=(\sigma\lessdot\rho\lessdot\la)$,
and $T\in\LRT(\sigma;\Rbt)$, then
$\phib_{\Rbt}(T_\CC) = \phib_{\Rbt}(T)_\CC$, that is,
the copy of $T\in\LRT(\sigma;\Rbt)$ labeled by $\CC$
is mapped to the copy of $\phib_{\Rbt}(T)\in \RC(\sigma^t;\Rbt^t)$
indexed by the same chain $\CC$.

We use the approach of Lemma \ref{nonsense}.
The back face commutes by induction,
the left and right faces by the definition of $\phib$,
the top face commutes up to the involution $\tau$
on $\LRT(\la^{--};\Rbt)$ by~\eqref{D and minus},
and the bottom face commutes up to the involution $\tau$
on $\RC(\la^{--t};\Rbt^t)$ by \eqref{j bar and tilde}.
Based on the above commutation up to $\tau$, it still follows
that $\jb \circ \jt \circ \phib_R = \jb \circ \phib_{\Rt} \circ \iit$
as maps into $\RC(\la^{--t};\Rbt^t)$.  Moreover
$\jb$ is injective as a map into the disjoint union
$\RC(\la^{--t};\Rbt^t)$ so it follows that the front face commutes.

In the remaining case, the last rectangle of $R$ has more than one column.
Then $\Rt\htt = \widetilde{\Rht}$.  We have the following diagram:
\begin{equation*}
\xymatrix{
{\LRT(\la;R)} \ar[rrr]^{\iit} \ar[ddd]_{\phib_R} \ar[dr]^{\iht}
& & & {\LRT(\la^-;\Rt)} \ar[ddd]^{\phib_{\Rt}} \ar[dl]_{\iht} \\
 & {\LRT(\la;\Rht)} \ar[r]^{\iit} \ar[d]_{\phib_{\Rht}} &
 {\LRT(\la^-;\Rt\htt)} \ar[d]^{\phib_{\Rt\htt}} & \\
 & {\RC(\la^t;{\Rht}^t)} \ar[r]_{\jt} & {\RC(\la^{-t};(\Rt\htt)^t)} & \\
{\RC(\la^t;R^t)} \ar[ur]_{\jht} \ar[rrr]_{\jt} & & & {\RC(\la^{-t};\Rt^t)}
	\ar[ul]^{\jht}
}
\end{equation*}
$\jht$ is injective, the back face commutes by induction and
the left and right faces by the definition of $\phib$.
The commutation of the top face is obvious.
The bottom face commutes by conjugating the result of
Lemma \ref{bar check} by $\comev$.  By Lemma \ref{nonsense} the
front face commutes.
\end{proof}

\subsection{Proof of the Evacuation Theorem}

\begin{thm}[Evacuation Theorem] \label{ev} The following diagram commutes:
\begin{equation*}
\begin{CD}
  \LRT(\la;R) @>{\ev}>> \LRT(\la;R^{\ev}) \\
  @V{\phib_R}VV @VV{\phib_{R^{\ev}}}V \\
  \RC(\la^t;R^t) @>>{\comev_R}> \RC(\la^t;R^{\ev t}).
\end{CD}
\end{equation*}
\end{thm}
\begin{proof} If $R$ is empty the result holds trivially.
Suppose that the last rectangle of $R$ has more than one column.
Obviously ${\Rht}^\ev = (R^{\ev})\ck$ and ${\Rck}^{\ev} = (R^{\ev})\htt$.
Consider the diagram:
\begin{equation*}
\xymatrix{
{\LRT(\la;R)} \ar[rrr]^{\ev} \ar[ddd]_{\phib_R} \ar[dr]^{\iht}
& & & {\LRT(\la;R^\ev)} \ar[ddd]^{\phib_{R^\ev}} \ar[dl]_{\ick} \\
 & {\LRT(\la;\Rht)} \ar[r]^{\ev} \ar[d]_{\phib_{\Rht}} &
	 {\LRT(\la;{\Rht}^\ev)} \ar[d]^{\phib_{{\Rht}^\ev}} & \\
 & {\RC(\la^t;{\Rht}^t)} \ar[r]_{\comev_{\Rht}} & {\RC(\la^t;{\Rht}^{\ev t})} 
 & \\
{\RC(\la^t;R^t)} \ar[ur]_{\jht} \ar[rrr]_{\comev_R} & & & 
{\RC(\la^t;R^{\ev t})}
	\ar[ul]^{\jck}
}
\end{equation*}
$\jck$ is injective, the back face commutes by induction,
the top and bottom faces are given by the commutative diagrams
\eqref{hat and ev} and \eqref{hat and theta},
the left face commutes by the definition of $\phib_R$,
and the right face commutes by \eqref{bij check}
with $R^\ev$ in place of $R$.
So by Lemma \ref{nonsense} the front face commutes.

Suppose that the last rectangle of $R$ is a single column.
Obviously $\Rb^{\ev} = \widetilde{R^{\ev}}$ and
$\Rt^{\ev} = \overline{R^{\ev}}$.  Consider the diagram:
\begin{equation*}
\xymatrix{
{\LRT(\la;R)} \ar[rrr]^{\ev} \ar[ddd]_{\phib_R} \ar[dr]^{\ib}
& & & {\LRT(\la;R^\ev)} \ar[ddd]^{\phib_{R^\ev}} \ar[dl]_{\iit} \\
 & {\LRT(\la^-;\Rb)} \ar[r]^{\ev} \ar[d]_{\phib_{\Rb}} &
	 {\LRT(\la^-;\Rb^\ev)} \ar[d]^{\phib_{\Rb^\ev}} & \\
 & {\RC(\la^{-t};\Rb^t)} \ar[r]_{\comev_{\Rb}} & 
{\RC(\la^{-t};\Rb^{\ev t})} & \\
{\RC(\la^t;R^t)} \ar[ur]_{\jb} \ar[rrr]_{\comev_R} & & & {\RC(\la^t;R^{\ev t})}
	\ar[ul]^{\jt}
}
\end{equation*}
$\jt$ is injective, the back face commutes by induction,
the top face commutes by Proposition~\ref{ev and},
the bottom face commutes by~\eqref{jt def} replacing $R$
by $R^\ev$ and using that $\comev_R$ is an involution,
the left face commutes by the definition of $\phib_R$,
and the right face commutes by \eqref{bij tilde}
with $R^\ev$ in place of $R$.
So by Lemma \ref{nonsense} the front face commutes.
\end{proof}

\section{Another recurrence for $\phib_R$}
\label{sec row rec}

The bijection $\phib_R$ is defined by a recurrence that
removes columns from the last rectangle.  
In this section it is shown that $\phib_R$ may be defined 
by an analogous recurrence which removes rows from the last rectangle.
This recurrence shall be used to prove some properties
of the transpose maps on LR tableaux and rigged configurations.

\subsection{Splitting off the first or last row}
Let $\Rl$ be obtained from $R$ by splitting off the
first row from the first rectangle so that the first rectangle in
$\Rl$ is $(\eta_1)$.
Similarly, let $\Rg$ be obtained from $R$ by splitting off the
last row from the last rectangle.
Recall the transpose map $\LRtr$ of Definition~\ref{def LRtr}.
Note that $(R^t)\ck = \Rl^t$. Define 
$\il:\LRT(\la;R)\to\LRT(\la;\Rl)$ as
\begin{equation*}
 \il:=\LRtr\circ\ick\circ\LRtr.
\end{equation*}
Note that $\Rg^t=(R^t)\htt$.
Similarly define $\ig:\LRT(\la;R)\to\LRT(\la;\Rg)$ as
\begin{equation}\label{ig and iht}
\ig:=\LRtr\circ\iht\circ\LRtr.
\end{equation}
Observing that $\Rg^\ev = R^{\ev<}$, it follows from the definitions
of $\il$ and $\ig$, \eqref{hat and ev} and the commutativity of
$\LRtr$ and $\ev$ that the following diagram commutes:
\begin{equation} \label{great and ev}
\begin{CD}
  \LRT(\la;R) @>{\ig}>> \LRT(\la;\Rg) \\
  @V{\ev}VV @VV{\ev}V \\
  \LRT(\la;R^\ev) @>>{\il}> \LRT(\la;\Rg^\ev).
\end{CD}
\end{equation}

Let $\jg:\RC(\la^t;R^t)\rightarrow\RC(\la^t;\Rg^t)$ be the
inclusion map. Then define $\jl$ by the following commutative diagram:
\begin{equation} \label{great and theta}
\begin{CD}
  \RC(\la^t;R^t) @>{\jg}>> \RC(\la^t;\Rg^t) \\
  @V{\comev_R}VV @VV{\comev_{\Rg}}V \\
  \RC(\la^t;R^{\ev t}) @>>{\jl}> \RC(\la^t;\Rg^{\ev t}).
\end{CD}
\end{equation}
For $(\nu,J)\in\RC(\la^t;R^t)$ set $(\nu^>,J^>)=\jg(\nu,J)$.
Note that
\begin{equation*}
P_n^{(k)}(\nu^>)=P_n^{(k)}(\nu)+\chi(k=\eta_L)\chi(1\le n<\mu_L).
\end{equation*}
Since $\comev_R$ reverses the sequence of rectangles and complements the
quantum numbers this implies that
$\jl(\nu,J)$ is obtained from $(\nu,J)$ by replacing 
$(n,x)\in(\nu,J)^{(\eta_1)}$ by $(n,x+1)$ for $1\le n<\mu_1$, and leaving
all other riggings invariant. In particular, $\jl$ preserves the colabels.

\begin{lem}\label{lem less bij}
The following diagram commutes
\begin{equation} \label{less bij}
\begin{CD}
  \LRT(\la;R) @>{\il}>> \LRT(\la;\Rl) \\
  @V{\phib_R}VV @VV{\phib_{\Rl}}V \\
  \RC(\la^t;R^t) @>>{\jl}> \RC(\la^t;\Rl^t).
\end{CD}
\end{equation}
\end{lem}
\begin{proof} Suppose first that $L \ge 2$ and that $R_L$ has more than 
one column. Note that $\Rl\htt={R\htt}^<$.
\begin{equation*} 
\xymatrix{
{\LRT(\la;R)} \ar[rrr]^{\il} \ar[ddd]_{\phib_R} \ar[dr]^{\iht}
& & & {\LRT(\la;\Rl)} \ar[ddd]^{\phib_{\Rl}} \ar[dl]_{\iht} \\
 & {\LRT(\la;\Rht)} \ar[r]^{\il} \ar[d]_{\phib_{\Rht}} &
	 {\LRT(\la;\Rl\htt)} \ar[d]^{\phib_{\Rl\htt}} & \\
 & {\RC(\la^t;{\Rht}^t)} \ar[r]_{\jl} & {\RC(\la^t;(\Rl\htt)^t)} & \\
{\RC(\la^t;R^t)} \ar[ur]_{\jht} \ar[rrr]_{\jl} & & & {\RC(\la^t;\Rl^t)}
	\ar[ul]^{\jht}
}
\end{equation*}
The top face commutes since $\il$ and $\iht$ are
both relabelings that replace different subalphabets.
The bottom face commutes since $\jht$ adds singular strings
and $\jl$ preserves colabels and hence preserves singularity of strings.
The left and right faces commute by the definition of $\phib$,
and the back face commutes by induction.  Since $\jht$ is
injective, by Lemma \ref{nonsense} the front face commutes as desired.

Suppose $L\ge 2$ and the last rectangle is a single column.
\begin{equation*} 
\xymatrix{
{\LRT(\la;R)} \ar[rrr]^{\il} \ar[ddd]_{\phib_R} \ar[dr]^{\ib}
& & & {\LRT(\la;\Rl)} \ar[ddd]^{\phib_{\Rl}} \ar[dl]_{\ib} \\
 & {\LRT(\la^-;\Rb)} \ar[r]^{\il} \ar[d]_{\phib_{\Rb}} &
	 {\LRT(\la;\Rbl)} \ar[d]^{\phib_{\Rbl}} & \\
 & {\RC(\la^{-t};\Rb^t)} \ar[r]_{\jl} & {\RC(\la^t;\Rbl^t)} & \\
{\RC(\la^t;R^t)} \ar[ur]_{\jb} \ar[rrr]_{\jl} & & & {\RC(\la^t;\Rl^t)}
	\ar[ul]^{\jb}
}
\end{equation*}
The commutation of the top face is obvious.
Since $\jl$ preserves colabels it is straightforward to verify that
the bottom face commutes.  The left and right faces
commute by the definition of $\phib$ and the
back face by induction.  Since $\jb$ is injective by
Lemma \ref{nonsense} the front face commutes.

The remaining case is $L\le 1$. It may be assumed
that $\LRT(\la;R)\not=\emptyset$ for otherwise there is nothing to
show.  By \cite[Prop. 33]{S} $\LRT(\la;R)$ is a singleton
if $R$ has at most two rectangles.  But
$\il:\LRT(\la;R)\rightarrow \LRT(\la;\Rl)$ is an embedding and
$\Rl$ has at most two rectangles so that $\LRT(\la;\Rl)$ is a singleton.
Since $\phib_R$ and $\phib_{\Rl}$ are bijections,
$\RC(\la^t;R^t)$ and $\RC(\la^t;\Rl^t)$ are also singletons.
So the embedding $\jl$ must send the unique
element of $\RC(\la^t;R^t)$ to the unique
element of $\RC(\la^t;\Rl^t)$ and the required commutation follows.
\end{proof}

\begin{lem} \label{ig and jg}
The following diagram commutes:
\begin{equation} \label{great bij}
\begin{CD}
  \LRT(\la;R) @>{\ig}>> \LRT(\la;\Rg) \\
  @V{\phib_R}VV @VV{\phib_{\Rg}}V \\
  \RC(\la^t;R^t) @>>{\jg}> \RC(\la^t;\Rg^t).
\end{CD}
\end{equation}
\end{lem}
\begin{proof} Consider the diagram 
\begin{equation*} 
\xymatrix{
{\LRT(\la;R)} \ar[rrr]^{\ig} \ar[ddd]_{\phib_R} \ar[dr]^{\ev}
& & & {\LRT(\la;\Rg)} \ar[ddd]^{\phib_{\Rg}} \ar[dl]_{\ev} \\
 & {\LRT(\la;R^\ev)} \ar[r]^{\il} \ar[d]_{\phib_{R^\ev}} &
	 {\LRT(\la;R^{\ev<})} \ar[d]^{\phib_{R^{\ev<}}} & \\
 & {\RC(\la^t;R^{\ev t})} \ar[r]_{\jl} & {\RC(\la^t;R^{\ev<t})} & \\
{\RC(\la^t;R^t)} \ar[ur]_{\comev_R} \ar[rrr]_{\jg} & & & {\RC(\la^t;R^{>t})}
	\ar[ul]^{\comev_{\Rg}}
}
\end{equation*}
The back face commutes by \eqref{less bij} for $R^\ev$,
the top and bottom faces commute by \eqref{great and ev} and
\eqref{great and theta}, the left and right faces commute by
the Evacuation Theorem~\ref{ev} for $R$ and $\Rg$, and $\comev_{\Rg}$ is 
injective. Therefore the front face commutes by Lemma \ref{nonsense}.
\end{proof}

\subsection{Removal of a cell from a single row}

Suppose $R_L$ is a single row.  Then $\Rht$ is given by
splitting off one cell from the end of the row $R_L$,
and $\Rhtb$ is obtained from $R$ by removing the cell at the
end of $R_L$.  As usual define
$\LRT(\la^-;\Rhtb):=\displaystyle\bigcup_{\rho \lessdot \la}
\LRT(\rho;\Rhtb)$ and
$\RC(\la^{-t};\Rhtb^t):=\displaystyle\bigcup_{\rho \lessdot \la}
\RC(\rho^t;\Rhtb^t)$.
The bijection $-:\ST(\la)\rightarrow \ST(\la^-)$ restricts to an
injection
\begin{equation*}
	-:\LRT(\la;R)\rightarrow \LRT(\la^-;\Rhtb).
\end{equation*}
Define the map
$\del:\RC(\la^t;R^t)\rightarrow\RC(\la^{-t};\Rhtb^t)$ by
the following algorithm.  Let $(\nu,J)\in \RC(\la^t;R^t)$.
Define the integers $\ell^{(0)} = \dots = \ell^{(\eta_L-1)} = 1$.
For $k \ge \eta_L$, inductively define $\ell^{(k)}$ to be minimal
such that $\ell^{(k)} \ge \ell^{(k-1)}$ and
there is a singular string of length $\ell^{(k)}$ in
$(\nu,J)^{(k)}$.  Let $\rku(\nu,J)$ be the minimal index such that
such a singular string does not exist, and set
$\ell^{(k)}=\infty$ for $k\ge \rku(\nu,J)$.  Then
$\del(\nu,J)$ is obtained from $(\nu,J)$ by shortening each of the
selected singular strings in
$(\nu,J)^{(k)}$ (for $\eta_L\le k < \rku(\nu,J)$)
by one and keeping them singular,
and leaving all other strings unchanged.

\begin{lem} Suppose the last rectangle of $R$ is a single row.
The following diagram commutes:
\begin{equation} \label{del def}
\begin{CD}
  \LRT(\la;R) @>{-}>> \LRT(\la^-;\Rhtb) \\
  @V{\phib_R}VV @VV{\phib_{\Rhtb}}V \\
  \RC(\la^t;R^t) @>>{\del}> \RC(\la^{-t};\Rhtb^t).
\end{CD}
\end{equation}
\end{lem}
\begin{proof} Consider the diagram
\begin{equation*} 
\xymatrix{
{\LRT(\la;R)} \ar[rr]^{-} \ar[ddd]_{\phib_R} \ar[dr]^{\iht}
& & {\LRT(\la^-;\Rhtb)} \ar[ddd]^{\phib_{\Rhtb}}  \\
 & {\LRT(\la;\Rht)} \ar[ur]^{-} \ar[d]^{\phib_{\Rht}} & \\
 & {\RC(\la^t;{\Rht}^t)} \ar[dr]^{\db} & \\
{\RC(\la^t;R^t)} \ar[rr]_{\del} \ar[ur]^{\jht}
& & {\RC(\la^{-t};\Rhtb^t).}
}
\end{equation*}
This diagram may be viewed as a prism whose top and bottom are triangles
and whose front is the diagram \eqref{del def} which must be proven.
The back left and back right faces commute by the definition of
$\phib$.  The top triangle obviously commutes.  It suffices
to show that the bottom triangle commutes.  This is done by computing
$\db \circ \jht$ explicitly.
Let $(\nu,J)\in\RC(\la^t;R^t)$.
Then $\jht(\nu,J)=(\nuht,\Jht)$ is obtained from $(\nu,J)$
by adding a singular string of length $1$ to each of the first
$\eta_L-1$ rigged partitions.  Let $\lb^{(k)}$ be the
lengths of the singular strings chosen by $\db$ acting on
$(\nuht,\Jht)$.  Since $\lb^{(0)}=1$ by definition,
the singular strings that were added by $\jht$, are selected by
$\db$, so $\lb^{(k)}=1$ for $1\le k \le \eta_L-1$.
Then it is obvious that for $k\ge \eta_L$ that $\db$
acting on $(\nuht,\Jht)$, selects the same strings that
$\del$ does, acting on $(\nu,J)$.  The result is now clear.
\end{proof}

\subsection{The new recurrence for $\phib_R$}
\label{new rec}

If the last rectangle of $R$ is a single row, then one may
use the commutative diagram \eqref{del def}
to express $\phib_R$ inductively in terms of $\phib_{\Rhtb}$.
If the last rectangle of $R$ has more than one row, then one
may apply the commutative diagram \eqref{great bij}
to express $\phib_R$ in terms of $\phib_{\Rg}$.
It will be shown in Section~\ref{sec trans} that this recurrence, which 
also defines $\phib_R$, is in a sense transpose to the usual definition 
of the bijection $\phib_{R^t}$.

\section{Transpose maps}
\label{sec trans}

Recall the LR-transpose bijection 
$\LRtr:\LRT(\la;R)\rightarrow \LRT(\la^t;R^t)$
of Definition \ref{def LRtr}. An analogous RC-transpose bijection exists 
for the set of rigged configurations denoted by 
$\RCtr:\RC(\la^t;R^t)\rightarrow\RC(\la;R)$, which was described
in \cite[Section 9]{KS}. In the following we recall its definition
and prove the Transpose Theorem~\ref{transpose thm}.

Let $(\nu,J)\in\RC(\la^t;R^t)$ and let $\nu$ have
the associated matrix $m$ with entries $m_{ij}$ as in \cite[(9.2)]{KS}
\begin{equation} \label{m def}
  m_{ij} = \alpha^{(i-1)}_j - \alpha^{(i)}_j
\end{equation}
for $i,j\ge 1$, where $\alpha^{(i)}_j$ is the size of the
$j$-th column of the partition $\nu^{(i)}$, recalling that
$\nu^{(0)}$ is defined to be the empty partition.
The configuration $\nu^t$ in $(\nu^t,J^t)=\RCtr(\nu,J)$ is defined
by its associated matrix $\mt$ given by
\begin{equation}\label{m transpose}
  \mt_{ij} = - m_{ji} +\chi((i,j)\in \la) -
	\sum_{a=1}^L \chi((i,j)\in R_a)
\end{equation}
for all $i,j \ge 1$. Here $(i,j) \in \la$ means that
the cell $(i,j)$ is in the Ferrers diagram of the partition $\la$
with $i$ specifying the row and $j$ the column.

For all $k,n\ge 1$, a rigging $J$ of $\nu$ 
determines a partition $J_n^{(k)}$ inside the rectangle of height
$m_n(\nu^{(k)})$ and width $P_n^{(k)}(\nu)$ given by the labels of the parts
of $\nu^{(k)}$ of size $n$. The partition $J_k^{t(n)}$ 
corresponding to $(\nu^t,J^t)=\RCtr(\nu,J)$ is defined as the transpose of 
the complementary partition to $J_n^{(k)}$ in the rectangle of 
height $m_n(\nu^{(k)})$ and width $P_n^{(k)}(\nu)$.

By \cite[(9.7)]{KS} and \cite[Lemma 10]{KS}, it follows that
\begin{equation} \label{Pm}
\begin{split}
  &\text{if $m_n(\nu^{(k)}) > 0$ then $P^{(k)}_n(\nu)=m_k(\nu^{t(n)})$} \\
  &\text{if $m_k(\nu^{t(n)}) > 0$ then $P^{(n)}_k(\nu^t)=m_n(\nu^{(k)})$.}
\end{split}
\end{equation}
This is weaker than the assertion \cite[(9.10)]{KS},
which does not seem to follow from only the assumption $m_n(\nu^{(k)}) > 0$
or $m_k(\nu^{t(n)})>0$.  Fortunately 
equation \eqref{Pm} still suffices to well-define the map $\RCtr$.
To see this, observe that it still follows from \cite[Section 9]{KS}
that $\RCtr$ is an involution on the level of configurations.
For $\RCtr$ to be well-defined for the riggings,
it is enough to show that the rectangle of height
$m_n(\nu^{(k)})$ and width $P^{(k)}_n(\nu)$
and the rectangle of height $m_k(\nu^{t(n)})$ and width
$P^{(n)}_k(\nu^t)$, are either transposes of each other or are
both empty.  But this follows from \eqref{Pm}.
Hence $\RCtr$ is an involution for rigged configurations.

We shall prove \cite[Conjecture 16]{KS}.
\begin{thm}[Transpose Theorem] \label{transpose thm} 
The following diagram commutes:
\begin{equation} \label{transpose bij}
\begin{CD}
  \LRT(\la;R) @>{\LRtr}>> \LRT(\la^t;R^t) \\
  @V{\phib_R}VV @VV{\phib_{R^t}}V \\
  \RC(\la^t;R^t) @>>{\RCtr}> \RC(\la;R).
\end{CD}
\end{equation}
\end{thm}
Theorem \ref{transpose thm} follows again from Lemma~\ref{nonsense}
by the usual arguments and requires the following two results.
Note that $R^{t>t}=\Rht$.
\begin{lem} \label{jht and jg}
The following diagram commutes:
\begin{equation} \label{RC tr hat}
\begin{CD}
  \RC(\la^t;R^t) @>{\jht}>> \RC(\la^t;{\Rht}^t) \\
  @V{\RCtr}VV @VV{\RCtr}V \\
  \RC(\la;R) @>>{\jg}> \RC(\la;\Rht).
\end{CD}
\end{equation}
\end{lem}
\begin{proof} 
Let $(\nu,J)\in\RC(\la^t;R^t)$ and set $(\nuht,\Jht)=\jht(\nu,J)$.
Recall that $(\nuht,\Jht)$ is obtained from $(\nu,J)$ by adding a 
singular string of length $\mu_L$ to each of the first $\eta_L-1$ 
rigged partitions in $(\nu,J)$.
Let $\mht$ be the matrix associated to $\nuht$. Then
\begin{equation*}
  \mht_{ij}=m_{ij}+\chi(1\le j\le \mu_L)
    \bigl\{\chi(i=\eta_L)-\chi(i=1)\bigr\}.
\end{equation*}
Using furthermore that
\begin{equation*}
\begin{split}
  \sum_{a=1}^{L+1}\chi((i,j)\in\Rht_a)&=
  \sum_{a=1}^L\chi((i,j)\in R_a)\\
   &+\chi(1\le i\le \mu_L)\bigl\{\chi(j=1)-\chi(j=\eta_L)\bigr\}
\end{split}
\end{equation*}
one finds by \eqref{m transpose} that ${\mht_{ij}}^t=m_{ij}^t$. Since
$\jg$ is an inclusion it follows that the configurations of
$\RCtr\circ\jht(\nu,J)$ and $\jg\circ\RCtr(\nu,J)$ coincide.
Since $\jht$ preserves the vacancy numbers by Lemma \ref{j hat}
also the riggings coincide.
\end{proof}

\begin{lem} \label{lem transpose}
Suppose the last rectangle of $R$ consists of a single column.
Then the following diagram commutes:
\begin{equation} \label{RC tr del}
\begin{CD}
  \RC(\la^t;R^t) @>{\jb}>> \RC(\la^{-t};\Rb^t) \\
  @V{\RCtr}VV @VV{\RCtr}V \\
  \RC(\la;R) @>>{\del}> \RC(\la^-;\Rb).
\end{CD}
\end{equation}
\end{lem}
The proof of this Lemma is given in Appendix \ref{App B}.

\begin{proof}[Proof of Theorem \ref{transpose thm}]
If $R$ is empty the result holds trivially.
If $R_L$ has more than one column consider the diagram:
\begin{equation*}
\xymatrix{
{\LRT(\la;R)} \ar[rrr]^{\LRtr} \ar[ddd]_{\phib_R} \ar[dr]^{\iht}
& & & {\LRT(\la^t;R^t)} \ar[ddd]^{\phib_{R^t}} \ar[dl]_{\ig} \\
 & {\LRT(\la;\Rht)} \ar[r]^{\LRtr} \ar[d]_{\phib_{\Rht}} &
	 {\LRT(\la^t;{\Rht}^t)} \ar[d]^{\phib_{{\Rht}^t}} & \\
 & {\RC(\la^t;{\Rht}^t)} \ar[r]_{\RCtr} & {\RC(\la;\Rht)} & \\
{\RC(\la^t;R^t)} \ar[ur]_{\jht} \ar[rrr]_{\RCtr} & & & 
{\RC(\la;R)} \ar[ul]^{\jg}
}
\end{equation*}
Since $\LRtr$ is an involution the top face commutes by \eqref{ig and iht}.
The bottom face commutes by Lemma \ref{jht and jg} and the back face by 
induction. The left face is the commutative diagram~\eqref{bij hat}
and the right face that of Lemma \ref{ig and jg}.
Since $\jg$ is injective the front face commutes by Lemma \ref{nonsense}.

Suppose the last rectangle is a single column. Consider the diagram:
\begin{equation*}
\xymatrix{
{\LRT(\la;R)} \ar[rrr]^{\LRtr} \ar[ddd]_{\phib_R} \ar[dr]^{\ib}
& & & {\LRT(\la^t;R^t)} \ar[ddd]^{\phib_{R^t}} \ar[dl]_{\ib} \\
 & {\LRT(\la^\ib;\Rb)} \ar[r]^{\LRtr} \ar[d]_{\phib_{\Rb}} &
	 {\LRT(\la^{\ib t};\Rb^t)} \ar[d]^{\phib_{\Rb^t}} & \\
 & {\RC(\la^{\ib t};\Rb^t)} \ar[r]_{\RCtr} & {\RC(\la^\ib;\Rb)} & \\
{\RC(\la^t;R^t)} \ar[ur]_{\jb} \ar[rrr]_{\RCtr} & & & 
{\RC(\la;R)} \ar[ul]^{\del}
}
\end{equation*}
$\del$ is injective, the back face commutes by induction,
the commutativity of the top is trivial and that of the bottom face
follows from Lemma \ref{lem transpose}. The left face commutes 
by~\eqref{bij bar} and the right face
by \eqref{del def} with $\la$ and $R$ replaced by their transposes.
So by Lemma \ref{nonsense} the front face commutes.
\end{proof}

\section{Embeddings}
\label{sec embed}

In  \cite[Section 6.1]{SW} and \cite[Section 2.3]{S2} embeddings
were given between sets of LR tableaux of the form $\CLRT(\la;R)$.
We restate these embeddings in terms of the tableaux $\LRT(\la;R)$
and show in Theorem~\ref{embed bij} that they are induced by
\textit{inclusions} of the corresponding sets of rigged configurations
under the coquantum version $\phit$ of the map $\phib$ 
(see~\eqref{phi coquantum}), thereby proving \cite[Conjecture 18]{KS}.

\subsection{Embedding definitions}
Define the partial order $\la \dom \mu$ on partitions
by $|\la|=|\mu|$ and
$\la_1+\dots+\la_i \ge \mu_1 + \dots + \mu_i$ for all $i$.
Let $R$ and $R'$ be two sequences of rectangles.
Recall that $\xi^{(k)}(R)$ is the partition whose parts
are the heights of the rectangles in $R$ of width $k$.
Say that $R \dom R'$ if $\xi^{(k)}(R) \dom \xi^{(k)}(R')$
for all $k\ge 1$.  Clearly $R\dom R'$ and $R'\dom R$ if and
only if $R'$ is a reordering of $R$.
Thus the relation $R \dom R'$ is a pseudo-order.
It is generated by the following two relations (see \cite[Section 2.3]{S2})
\begin{enumerate}
\item[(E1)] $R \dom \Rp$ where $R_i=\Rp_i$ for $i>2$, $R_1=(c^a)$,
$R_2=(c^b)$, $\Rp_1=(c^{a-1})$, $\Rp_2=(c^{b+1})$ for $a-1\ge b+1$ 
and $c$ a positive integer.
\item[(E2)] $R \dom s_p R$ where $s_p R$ denotes the sequence obtained 
from $R$ by exchanging the rectangles $R_p$ and $R_{p+1}$.  
\end{enumerate}

For the relation (E1) define the embedding
$\ip$ by the commutation of the diagram
\begin{equation} \label{i plus def}
\begin{CD}
  \LRT(\la;R) @>{\ip}>> \LRT(\la;\Rp) \\
  @V{\LRtr}VV @VV{\LRtr}V \\
  \LRT(\la^t;R^t) @>>{\mathrm{inclusion}}> \LRT(\la^t;\Rpt).
\end{CD}
\end{equation}
For the relation (E2) we have the following result.

\begin{defprop} For $1\le p \le L-1$ there are unique
bijections $\sigma_p:\LRT(\la;R)\rightarrow\LRT(\la;s_p R)$
satisfying the following properties:
\begin{enumerate}
\item If $p < L-1$ then $\sigma_p$ commutes with restriction to
the initial interval $B-B_L$ where $B$ and $B_L$ are as in the
definition of $\LRT(\la;R)$.
\item If $p=L-1$ then the following diagram commutes:
\begin{equation*}
\begin{CD}
  \LRT(\la;R) @>{\sigma_p}>> \LRT(\la;s_p R) \\
  @V{\ev}VV @VV{\ev}V \\
  \LRT(\la;R^\ev) @>>{\sigma_1}> \LRT(\la;s_1 (R^\ev))
\end{CD}
\end{equation*}
\end{enumerate}
\end{defprop}
\begin{proof} One may reduce to the case $p=L-1$ using 1,
then to the case $p=1$ using 2, and then to the case
$p=1$ and $L=2$ using 1 again.  In this case the
sets of tableaux are all empty or all singletons
by \cite[Prop. 33]{S} and the result holds trivially.
\end{proof}

\begin{rem} \label{autos} In the case that $R_j$ is a single row
of length $\eta_j$ for all $j$, by Example \ref{std ex} there
is a bijection of $\LRT(\la;R)$ with the set of column-strict tableaux
of shape $\la$ and content $\eta$.  The action of the
bijections $\sigma_p$ on the column-strict tableaux are the
automorphisms of conjugation defined by Lascoux and
Sch\"utzenberger \cite{LS}.
\end{rem}

Let $R \dom R'$.  Then $R'$ may be obtained from $R$
by a sequence of transformations of the form (E1) and (E2);
fix such a sequence.  Define the
embedding $i^{R'}_R:\LRT(\la;R)\to\LRT(\la;R')$
as the corresponding composition of embeddings of the form
$\ip$ and $\sigma_p$.  By \cite[Theorem 4]{S2} it follows that
the embedding $i^{R'}_R$ is independent of the sequence of transformations
(E1) and (E2) leading from $R$ to $R'$.  This fact also
follows immediately from Theorem \ref{embed bij} below.

For rigged configurations, it follows immediately
from the definitions that if $R\dom R'$ then there
is an inclusion
$\RC(\la^t;R^t) \subseteq \RC(\la^t;(R')^t)$
which shall be denoted by $j^{R'}_R$.

\begin{thm}[Embedding Theorem] \label{embed bij} 
Let $R\dom R'$. Then the diagram commutes:
\begin{equation} \label{embed}
\begin{CD}
  \LRT(\la;R) @>{i^{R'}_R}>> \LRT(\la;R') \\
  @V{\phit_R}VV @VV{\phit_{R'}}V \\
  \RC(\la^t;R^t) @>{j^{R'}_R}>> \RC(\la^t;(R')^t).
\end{CD}
\end{equation}
\end{thm}

Clearly it suffices to prove Theorem \ref{embed bij} in the
cases (E1) and (E2).

\subsection{The case (E1)}
Suppose $R \dom \Rp$ as in (E1).
Define the embedding $\jp:\RC(\la^t;R^t)\to\RC(\la^t;\Rpt)$
by the commutativity of the diagram
\begin{equation} \label{j plus def}
\begin{CD}
  \RC(\la^t;R^t) @>{\mathrm{inclusion}}>> \RC(\la^t;\Rpt) \\
  @V{\com_R}VV @VV{\com_{\Rp}}V \\
  \RC(\la^t;R^t) @>>{\jp}> \RC(\la^t;\Rpt).
\end{CD}
\end{equation}
Note that $\jp$ preserves colabels.
The following result immediately proves
Theorem \ref{embed bij} in the case (E1).

\begin{lem} \label{plus bijection} The diagram commutes:
\begin{equation*}
\begin{CD}
  \LRT(\la;R) @>{\ip}>> \LRT(\la;\Rp) \\
  @V{\phib_R}VV @VV{\phib_{\Rp}}V \\
  \RC(\la^t;R^t) @>>{\jp}> \RC(\la^t;\Rpt).
\end{CD}
\end{equation*}
\end{lem}
\begin{proof}
For $L\ge 3$ the proof of this lemma is the same as the proof of
Lemma~\ref{lem less bij} with $\Rl$ (resp. $\il,\jl$) replaced by
$\Rp$ (resp. $\ip,\jp$). Note that in the remaining case $L\le 2$
both $R$ and $\Rp$ have at most two rectangles so that $\LRT(\la;R)$
and $\LRT(\la;\Rp)$ are singletons, and the proof follows again by 
arguments similar to those of the proof of Lemma~\ref{lem less bij}.
\end{proof}

\subsection{The case (E2)}
The case (E2) of Theorem \ref{embed bij} is an immediate
consequence of the following result.

\begin{lem} \label{switch bij} The diagram commutes:
\begin{equation} \label{switch}
\begin{CD}
\LRT(\la;R) @>{\sigma_p}>> \LRT(\la;s_p R) \\
@V{\phib_R}VV @VV{\phib_{s_p R}}V \\
\RC(\la^t;R^t) @= \RC(\la^t;(s_p R)^t). \\
\end{CD}
\end{equation}
\end{lem}
\begin{proof} Recall the intervals of integers $B_j$ ($1\le j \le L$)
in the definition of the set $\LRT(\la;R)$.  It follows from the
definition of $\sigma_p$ that for all $S \in \LRT(\la;R)$,
$(\sigma_p S)|_{B_j} = S|_{B_j}$ for $j > p+1$.
Using this fact and the definition of $\phib_R$,
one may reduce to the case $L = p+1$.
So it may be assumed that $L=p+1$.
Obviously $(s_p R)^\ev = s_1(R^\ev)$.  Consider the diagram
\begin{equation*} 
\xymatrix{
{\LRT(\la;R)} \ar[rrr]^{\sigma_p} \ar[ddd]_{\phib_R} \ar[dr]^{\ev}
& & & {\LRT(\la;s_p R)} \ar[ddd]^{\phib_{s_p R}} \ar[dl]_{\ev} \\
 & {\LRT(\la;R^\ev)} \ar[r]^{\sigma_1} \ar[d]_{\phib_{R^\ev}} &
	 {\LRT(\la;s_1(R^\ev))} \ar[d]^{\phib_{s_1 (R^\ev)}} & \\
 & {\RC(\la^t;R^{\ev t})} \ar@{=}[r] & {\RC(\la^t;(s_1 R^{\ev t}))} & \\
{\RC(\la^t;R^t)} \ar[ur]_{\comev_R} \ar@{=}[rrr] & & & {\RC(\la^t;s_p R^t)}
	\ar[ul]^{\comev_{s_p R}}
}
\end{equation*}
The commutation of the bottom face is obvious.
The top face commutes by definition and the
left and right faces commute by Theorem \ref{ev}.
The desired commutation of the front face
is reduced to the commutation of the back face.
Replacing $R$ by $R^\ev$ it may now be assumed that $p=1$.  
Applying a previous reduction it may also be assumed that $L=2$.
But in this case the sets in \eqref{switch} are all singletons
or all empty, so the diagram \eqref{switch} commutes.
\end{proof}

\subsection{Connection with $\CLRT(\la;R)$}

In \cite[Section 10]{KS}, \cite[Section 6.1]{SW} and \cite[Section 2.3]{S2}
embeddings $\theta^{R'}_R:\CLRT(\la;R)\to\CLRT(\la;R')$ were
defined when $R\dom R'$. They are related to the embeddings 
$i^{R'}_R:\LRT(\la;R)\rightarrow \LRT(\la;R')$ by
the bijection $\beta_R:\CLRT(\la;R)\to\LRT(\la;R)$ of
Remark~\ref{LR defs}.

\begin{prop} \label{relabel embed} Let $R\dom R'$.  The diagram commutes:
\begin{equation} \label{beta theta}
\begin{CD}
  \CLRT(\la;R) @>{\theta^{R'}_R}>> \CLRT(\la;R') \\
  @V{\beta_R}VV @VV{\beta_{R'}}V \\
  \LRT(\la;R) @>>{i^{R'}_R}> \LRT(\la;R').
\end{CD}
\end{equation}
\end{prop}

For the proof recall the
evacuation map $\Ev:\CLRT(\la;R)\to\CLRT(\la;R^\ev)$.
Let $n=\sum_{j=1}^L \mu_j$ be the sum of heights
of rectangles in $R$.  There is a unique
involution $T\mapsto T^\Ev$ on column-strict tableaux of shape
$\la$ in the alphabet $[1,n]$ defined by
\begin{equation*}
  \shape((T^\Ev)|_{[1,i]}) = \shape(P(T|_{[n-i+1,n]}))
\end{equation*}
for all $1\le i\le n$.  The bijection
$\Ev$ restricts to a bijection
$\CLRT(\la;R)\rightarrow \CLRT(\la;R^\ev)$ \cite{KS}.

\begin{lem} \label{ev beta} The diagram commutes:
\begin{equation} \label{ev LRT}
\begin{CD}
  \CLRT(\la;R) @>{\Ev}>> \CLRT(\la;R^\ev) \\
  @V{\beta_R}VV @VV{\beta_{R^\ev}}V \\
  \LRT(\la;R) @>>{\ev}> \LRT(\la;R^\ev).
\end{CD}
\end{equation}
\end{lem}
\begin{proof}[Sketch of Proof] 
$\beta_R = \gamma_R^{-1} \circ \std$ by Remark \ref{LR defs}.
Using the well-known fact
$\std \circ \Ev = \ev \circ \std$
it is enough to show that the following diagram commutes:
\begin{equation*}
\begin{CD}
  \RLRT(\la;R) @>{\ev}>> \RLRT(\la;R^\ev) \\
  @V{\gamma_R}VV @VV{\gamma_{R^\ev}}V \\
  \LRT(\la;R) @>>{\ev}> \LRT(\la;R^\ev).
\end{CD}
\end{equation*}
By \cite{LS} $S^\ev = P(\# \word(S))$ for $S\in \ST(\la)$
where $\# w$ is the reverse of the complement of the word $w$
in the alphabet $[1,|\la|]$.
Let $\WR(R)$ (resp. $\WC(R)$) be the set of standard words $w$ such that
the standard tableau $P(w)$ is in $\RLRT(\la;R)$ (resp. $\LRT(\la;R)$).
The relabeling bijection $\gamma_R:\RLRT(\la;R)\to\LRT(\la;R)$ extends
to a map $\WR(R)\to\WC(R)$ in an obvious way.  This relabeling
map satisfies $P(\gamma_R(w))=\gamma_R(P(w))$ \cite[Lemma 36]{S2}.
It is not hard to see that $\# \circ \gamma_R = \gamma_{R^\ev} \circ \#$
as maps $\WR(R)\rightarrow\WC(R^\ev)$.  For all $S\in\RLRT(\la;R)$ we have
\begin{equation*}
\begin{split}
  (\gamma_R(S))^\ev &= P(\# \gamma_R(S))=P(\# \gamma_R(\word(S)))\\
  &= P(\gamma_{R^\ev}(\# \word(S))) = \gamma_{R^\ev}(P(\# \word(S))) \\
  &= \gamma_{R^\ev}(S^\ev).
\end{split}
\end{equation*}
\end{proof}

\begin{proof}[Proof of Proposition \ref{relabel embed}]
Again the result follows from the special cases
(E1) and (E2).  Consider the case (E1).
By the definition of the relabeling map $\beta_R$ and
the fact that in this case both $\theta^{R'}_R$ and
$i^{R'}_R$ only change the subtableaux corresponding to the
first two rectangles, one may reduce to the case that $L=2$.
By \cite[Prop. 33]{S} all sets are either singletons
or empty.  The vertical maps are bijections and the
horizontal maps are embeddings, so the diagram must
commute.

Denote by $s_p$ the rectangle switching bijection on
$\CLRT(\la;R)$ as defined in \cite[Section 7]{KS} (see also 
\cite[Section 2.5]{S}).
In the case (E2), neither $s_p$ nor $\sigma_p$ disturb
the subtableaux corresponding to the rectangles
$R_j$ for $j>p+1$.  Therefore it may be assumed that
$L=p+1$.  Now the evacuation map is applied to the
diagram \eqref{beta theta}:

\begin{equation*}
\xymatrix{
{\CLRT(\la;R)} \ar[rrr]^{s_p} \ar[ddd]_{\beta_R} \ar[dr]^{\Ev} & & &
{\CLRT(\la;s_p R)} \ar[ddd]^{\beta_{s_p R}} \ar[dl]_{\Ev} \\
& {\CLRT(\la;R^\ev)} \ar[r]^{s_1} \ar[d]_{\beta_{R^\ev}}
& {\CLRT(\la;s_1(R^\ev))} \ar[d]^{\beta_{s_1 (R^\ev)}} & \\
& {\LRT(\la;R^\ev)} \ar[r]_{\sigma_1} & {\LRT(\la;s_1(R^\ev))} & \\
{\LRT(\la;R)} \ar[rrr]_{\sigma_p} \ar[ur]_{\ev} & & &
{\LRT(\la;s_p R)} \ar[ul]_{\ev} 
}
\end{equation*}
The left and right faces commute by Lemma \ref{ev beta}.
The top face commutes by the proof of \cite[Lemma 5]{KS},
and the bottom commutes by the definition of $\sigma_p$.
To obtain the desired commutation of the front face,
it is enough to show that the back face commutes.
But the back face is the special case $p=1$
with $R$ replaced by $R^\ev$.
This reduces to the case $p=1$.
As before one may reduce to the case that $L=p+1=2$.
But in this case by \cite[Prop. 33]{S} the
sets are all singletons or all empty, so the diagram commutes.
\end{proof}

\subsection{Single rows}
\label{single row sec}

Let $\rows(R)$ be obtained by slicing all of the rectangles
of $R$ into single rows.  Clearly $R \dom \rows(R)$
and $\rows(R)$ is minimal with respect to the pseudo-order
$\dom$.  Define
$i_R := i^{\rows(R)}_R:\LRT(\la;R)\to\LRT(\la;\rows(R))$.
Recall that Schensted's standardization map
gives a bijection from column-strict tableaux of shape
$\la$ and content $(\eta_1^{\mu_1},\eta_2^{\mu_2},\dots,\eta_L^{\mu_L})$
onto $\LRT(\la;\rows(R))$, so that $i_R$ essentially embeds
$\LRT(\la;R)$ into column-strict tableaux.

\section{Statistics for LR tableaux and rigged configurations}
\label{sec statistics}

Recall that both the set of LR tableaux and the set of rigged 
configurations are endowed with statistics, given by the
charge $\charge_R(T)$ for $T\in\LRT(\la;R)$ and $\cc(\nu,J)$
for $(\nu,J)\in\RC(\la^t;R^t)$. The statistic $\cc(\nu,J)$ was 
given in~\eqref{RC charge} and an explicit expression for the 
generalized charge $\charge_R(T)$ can be found in~\cite[(2.1)]{S}.

The objective of this section is to prove that $\phit$ preserves
the statistics, which settles \cite[Conjecture 9]{KS}.

\begin{thm}\label{thm charge}
Let $T\in\LRT(\la;R)$. Then $\charge_R(T)=\cc(\phit_R(T))$.
\end{thm}
\begin{proof} Recall the embedding $i_R:\LRT(\la;R)\to\LRT(\la;\rows(R))$
of Section \ref{single row sec}.  By
\cite[Prop. 7]{S2}, for all $S\in \LRT(\la;R)$ we have
$\charge_{\rows(R)}(i_R(S))=\charge_R(S)$.
On the other hand by Theorem \ref{embed bij},
$\phit_{\rows(R)}(i_R(S))=\phit_R(S)$.
By replacing $R$ by $\rows(R)$ it may be assumed that
$R$ consists of single rows.

The transpose maps are applied to reduce to the case of single columns.
Recall the definition of the number $n(R)$ \cite[(2.7)]{KS}.
For all $S\in\LRT(\la;R)$ we have
$\charge_{R^t}(\LRtr(S)) = n(R) - \charge_R(S)$ 
by \cite[Prop. 24 and Theorem 26]{S2} and \cite[Lemma 6.5]{SW}.  
For all $(\nu,J)\in\RC(\la^t;R^t)$ we have
$\cc(\RCtr(\nu,J))=n(R)-\cc(\nu,J)$ by \cite[(9.12)]{KS}.
Moreover
\begin{equation*}
\begin{split}
  \RCtr \circ \phit_R &= \RCtr \circ \com_R \circ \phib_R =
	\com_{R^t} \circ \RCtr \circ \phib_R \\
  &= \com_{R^t} \circ \phib_{R^t} \circ \LRtr = \phit_{R^t} \circ \LRtr
\end{split}
\end{equation*}
by the definitions of $\phit_R$, $\RCtr$, $\com_R$, and
Theorem \ref{transpose thm}.  From these facts, by replacing
$R$ with $R^t$ it may be assumed that
$R$ consists of single columns.

Applying the first reduction again, it may be assumed that
$R$ consists of single cells. Observe that in this case
$\LRT(\la;R)=\ST(\la)$. It is possible to verify
$\charge_R(S)=\cc(\phit_R(S))$ for $S\in\ST(\la)$ directly.
It is known that for all $S\in\ST(\la)$,
$\charge_R(S)=\charge(S)$ where $\charge(S)$ is the ordinary charge
of a standard tableau $S$ \cite[Section III.6]{Mac}.
Let $\asc(S)$ denote the number of ascents in $S$, that
is, the number of indices $i$ such that $i+1$ is in a later
column in $S$ than $i$.  It follows immediately from the
definition of $\charge$ that $\charge(S^-) = \charge(S) - \asc(S)$.
By induction it suffices to show that if $\phit_R(S)=(\nu,J)$ then
$\asc(S)$ is equal to $\alpha_1^{(1)}$ and that 
$\cc(\nu,J)-\cc(\dt(\nu,J))=\alpha_1^{(1)}$.

It is first shown that $\asc(S)=\alpha^{(1)}_1$.
Denote $(\nut,\Jt)=\dt(\nu,J)$, let $\at_n^{(k)}$ be the size
of the $n$-th column in $\nut^{(k)}$ and
let $\lt^{(k)}$ be the length of the selected string of label zero
in $(\nu,J)^{(k)}$ in passing to $(\nut,\Jt)^{(k)}$.
By induction $\asc(S^-)=\at_1^{(1)}$ where
$\phit_{\Rt}(S^-)=\dt(\nu,J)=(\nut,\Jt)$.
(At this point the astute reader may be concerned that the relation
between $\phit_R$ and $\phib_R$ involves $\com_R$ whereas the relation
between $\jb$ and $\dt$ involves $\comev_R$ (compare \eqref{jt def});
however, since $R$ consists only of single boxes $\com_R=\comev_R$).
Now
\begin{equation*}
  \alpha_1^{(1)}- \at_1^{(1)} = \chi(\lt^{(1)}=1).
\end{equation*}
Thus it must be shown that $\lt^{(1)}=1$ if and only if
$N-1$ is an ascent of $S$ where $N=|\la|$. 
Let $c$ and $c^-$ be the column indices of the letters $N$ and $N-1$ in $S$
and $\ltt^{(k)}$ the length of the string selected by $\dt$ in
$(\nut,\Jt)^{(k)}$.
First suppose that $\lt^{(1)}=1$. Then $\ltt^{(k)}\ge \lt^{(k)}$ for all
$0\le k<c$. Hence $c>c^-$ by \eqref{dt vacancy} with $k=c-1$.
Conversely, if $\lt^{(1)}>1$ then $\ltt^{(k)}<\lt^{(k)}$ for all
$0\le k<c$ so that $c\le c^-$.
This proves $\asc(S)=\alpha_1^{(1)}$.

Now it is shown that $\cc(\nu,J)-\cc(\dt(\nu,J))=\alpha^{(1)}_1$
when $R$ consists of single cells. To this end we first compute 
\begin{equation} \label{config diff}
\begin{split}
 \cc(\nu)-\cc(\nut)
  &= \sum_{k,n\ge 1} \alpha^{(k)}_n(\alpha^{(k)}_n-\alpha^{(k+1)}_n) \\
  &- \sum_{k,n\ge 1} (\alpha^{(k)}_n-\delta_{n,\lt^{(k)}})
	(\alpha^{(k)}_n-\alpha^{(k+1)}_n-\delta_{n,\lt^{(k)}}
        +\delta_{n,\lt^{(k+1)}}) \\
  &= \sum_{k\ge 1} (2 \alpha^{(k)}_{\lt^{(k)}} - \alpha^{(k)}_{\lt^{(k+1)}}-
     \alpha^{(k+1)}_{\lt^{(k)}}-\chi(\lt^{(k)}<\lt^{(k+1)})).
\end{split}
\end{equation}
Next let us determine the difference of the sums of quantum numbers.
Recall that $\dt$ does not change the colabels of unselected strings,
and that the vacancy numbers change by \eqref{dt vacancy}.
The selected strings have label zero before and after being shortened.
Then
\begin{equation} \label{coq diff}
\begin{split}
  &\sum_{k,n\ge 1} |J_n^{(k)}| - \sum_{k,n\ge 1} |\Jt_n^{(k)}| \\
  &= \sum_{k,n\ge 1} (m_n(\nu^{(k)})-\delta_{n,\lt^{(k)}})
     (\chi(\lt^{(k-1)}\le n< \lt^{(k)}) - \chi(\lt^{(k)}\le n< \lt^{(k+1)})) \\
  &= \sum_{k\ge 1} (\alpha^{(k)}_{\lt^{(k-1)}}-\alpha^{(k)}_{\lt^{(k)}})
    -\sum_{k\ge 1} (\alpha^{(k)}_{\lt^{(k)}}-\alpha^{(k)}_{\lt^{(k+1)}})
    +\sum_{k\ge 1} \chi(\lt^{(k)}<\lt^{(k+1)}).
\end{split}
\end{equation}
By \eqref{config diff} and \eqref{coq diff} we have
$\cc(\nu,J)-\cc(\nut,\Jt)=\alpha^{(1)}_{\lt^{(0)}}=\alpha^{(1)}_1$.
\end{proof}

The cocharge of an LR tableau $T\in\LRT(\la;R)$ is defined as
$\co_R(T)=\|R\|-\charge_R(T)$ where $\|R\|=\sum_{i<j}|R_i\cap R_j|$.
Recall that for $\la$ and $\mu$ partitions $\CST(\la;\mu)$
denotes the set of column-strict tableaux of shape $\la$ and content
$\mu$. By Example~\ref{std ex} the set $\CST(\la;\mu)$ equals 
$\CLRT(\la;R)$ where $R=((\mu_1),\ldots,(\mu_L))$, and identifying 
the two sets the cocharge of a column-strict tableau $T$ is given 
by $\co(T)=\co_R(T)$.

Recall the map $\beta_R:\CLRT(\la;R)\to\LRT(\la;R)$ of Remark~\ref{LR defs}.
Then, identifying $\CST(\la;\mu)$ and $\CLRT(\la;R)$, the coquantum
version of the original bijection~\cite{KR} between column-strict tableaux
and rigged configurations is given by $\phit_{R^t}\circ\LRtr\circ\beta_R$.

\begin{cor}
Let $T$ be a column-strict tableau of shape $\la$ and partition content 
$\mu=(\mu_1,\dots,\mu_L)$ and $1\le r\le L$ an index such that
$\mu_r>\mu_{r+1}$ where $\mu_{L+1}=0$.  Let $T'$ be obtained from $T$ by
applying the automorphism of conjugation (see Remark \ref{autos})
$\sigma=\sigma_{L-1} \sigma_{L-2} \dots \sigma_r$ that changes
the content to $(\mu_1,\dots,\mu_{r-1},\mu_{r+1},\dots,\mu_L,\mu_r)$,
then removing the rightmost copy of the letter $L$, then
applying the automorphism of conjugation $\sigma^{-1}$ to
change the content to
$(\mu_1,\dots,\mu_{r-1},\mu_r - 1,\mu_{r+1},\dots,\mu_L)$.  Let
$R=((\mu_1),\ldots,(\mu_L))$ and 
$(\nu,J)=\phit_{R^t}\circ\LRtr\circ\beta_R(T)$. Then
\begin{equation*}
   \co(T) - \co(T') = \alpha^{(1)}_{\mu_r}
\end{equation*}
where $\alpha_{\mu_r}^{(1)}$ is the size of the $\mu_r$-th column of 
$\nu^{(1)}$.
\end{cor}

\begin{proof}
By Lemma \ref{switch bij} it suffices to prove the result
in the case that $\mu=(\mu_1,\dots,\mu_L)$ is any sequence of
positive integers, $r=L$, and $T'$ is obtained from $T$ by removing
the rightmost copy of the letter $L$.
Set $S=\LRtr\circ\beta_R(T)$, $S'=\LRtr\circ\beta_{\overline{R^t}^t}(T')$,
$(\nu,J)=\phit_{R^t}(S)$ and $(\nu',J')=\phit_{\overline{R^t}}(S')$. 
By \cite[Lemma 6.5]{SW} and \cite[Theorem 26]{S2} 
$\co_R(\LRtr(S))=\charge_{R^t}(S)$. Hence
\begin{align*}
\co(T)-\co(T')&=\charge_{R^t}(S)-\charge_{\overline{R^t}}(S')&\\
&=\cc(\phit_{R^t}(S))-\cc(\phit_{\overline{R^t}}(S')&
\text{by Theorem~\ref{thm charge}}\\
&=\cc(\nu,J)-\cc(\nu',J').&
\end{align*}
Notice that $(\nu,J)\in\RC(\la;R)$ where $R^t$ is a sequence of single 
columns. Similar calculations to~\eqref{config diff} and~\eqref{coq diff}
yield $\cc(\nu,J)-\cc(\nu',J')=\alpha^{(1)}_{\mu_L}$.
\end{proof}

\appendix
\section{Proof of Lemma \ref{bar tilde}}\label{App A}

The proof of Lemma \ref{bar tilde} is given here
by Lemmas \ref{db and dt} and \ref{label commute}.
After the Lemmas it is explained
why the statements regarding $\rkb$ and $\rkt$ in
Lemma \ref{bar tilde} follow from Lemma \ref{db and dt}.

In this section the following notation is used.
Let $R=(R_1,\ldots,R_L)$ be a sequence of rectangles with
$R_i=(\eta_i^{\mu_i})$, such that $R_1$ and $R_L$ are single columns
and $|R|\ge 2$, let $\la$ be a partition and $(\nu,J)\in\RC(\la^t;R^t)$.
Write
\begin{equation*}
\begin{split}
  \db(\nu,J) &= (\nub,\Jb) \\
  \dt(\nu,J) &= (\nut,\Jt) \\
  \dt\circ\db(\nu,J) &= (\nubt,\Jbt) \\
  \db\circ\dt(\nu,J) &= (\nutb,\Jtb).
\end{split}
\end{equation*}
Furthermore, let $\lb^{(k)}$,  $\lt^{(k)}$,
$\lbt^{(k)}$ and $\ltb^{(k)}$ denote the lengths of the strings
that are shortened in the transformations
$(\nu,J)\mapsto (\nub,\Jb)$,
$(\nu,J)\mapsto (\nut,\Jt)$,
$(\nub,\Jb)\mapsto (\nubt,\Jbt)$ and
$(\nut,\Jt)\mapsto (\nutb,\Jtb)$, respectively.

\begin{rem} \label{undisturbed}
Note that, except for the strings that change length in the transformations,
a string has label zero before applying $\db$ if and only it does
afterwards, and a string is singular before applying $\dt$ if and only
if it is afterwards.
\end{rem}

If $\db$ and $\dt$ select the same string in $\nu^{(k)}$
(that is, both select a string with the same length and same
label), say that the doubly singular case holds for $\nu^{(k)}$.
Otherwise say that the generic case holds for $\nu^{(k)}$.

The next Lemma shows that $\nutb=\nubt$.

\begin{lem} \label{db and dt}
If the generic case holds for $\nu^{(k)}$ then
$\lb^{(k)}=\ltb^{(k)}$ and $\lt^{(k)}=\lbt^{(k)}$.
Otherwise suppose the doubly singular case holds for
$\nu^{(k)}$.  Let $\ell := \lb^{(k)}=\lt^{(k)}$ be the common
string length.  Then
\begin{enumerate}
\item If $\lbt^{(k)} < \ell$ (or $\ltb^{(k)} < \ell$) then
$\lbt^{(k)} = \ltb^{(k)} = \ell - 1$, and
$m_{\ell-1}(\nu^{(k+1)})=0$.
\item If $\lbt^{(k)} = \ell$ (or $\ltb^{(k)} = \ell$) then
$\ltb^{(k)} = \lbt^{(k)}= \ell$. 
\item If $\lbt^{(k)} > \ell$ (or $\ltb^{(k)} > \ell$) then
$\lbt^{(k)} = \ltb^{(k)}$, $\lt^{(k+1)}\ge \lbt^{(k)}$ and
$\lb^{(k+1)} \ge \lbt^{(k)}$.  Moreover, if
$\lbt^{(k)} = \infty$ then $\la^t_k = \la^t_{k+1}$.
\end{enumerate}
\end{lem}

\begin{lem} \label{label commute} $\Jbt=\Jtb$.
\end{lem}

Together Lemmas \ref{db and dt} and \ref{label commute} 
establish the desired result that $(\nutb,\Jtb)=(\nubt,\Jbt)$.

Given Lemma \ref{db and dt}, the assertions of
Lemma \ref{bar tilde} on the relationships between
$\rb=\rkb(\nu,J)$, $\rt=\rkt(\nu,J)$,
$\rbt=\rkt(\db(\nu,J))$, and $\rtb=\rkb(\dt(\nu,J))$ are now established.

Suppose $\rb<\rt$.  Then $\lb^{(\rb)}=\infty>\lt^{(\rb)}$ and the
generic case holds for $(\nu,J)^{(\rb)}$.  It follows that
for all $k\ge \rb$, $\lbt^{(k)}=\lt^{(k)}$ and $\ltb^{(k)}=\lb^{(k)}$.
Therefore $\rbt = \rt$ and $\rtb \le \rb$.  Suppose $\rtb<\rb$.
Then $\ltb^{(\rtb)}=\infty > \lb^{(\rtb)}$.  It follows that
Case 3 holds for $(\nu,J)^{(\rtb)}$.  But then
$\lbt^{(\rtb)}=\ltb^{(\rtb)}=\infty$ so that
$\rbt \le \rtb < \rb < \rt$, which is a contradiction.
Therefore $\rtb=\rb$ and $\rbt=\rt$ as desired.

If $\rb>\rt$ a similar proof shows $\rtb=\rb$ and $\rbt=\rt$.
So it may be assumed that $r=\rt=\rb$.
Notice that $\lt^{(r)}=\lb^{(r)}=\infty$ implies that
$\lbt^{(r)}=\ltb^{(r)}=\infty$ so that $\rbt \le r$ and
$\rtb \le r$.
Suppose Case 3 does not hold for $(\nu,J)^{(r-1)}$.  Then
$\lbt^{(r-1)} \le \lt^{(r-1)}<\infty$ and
$\ltb^{(r-1)} \le \lb^{(r-1)}<\infty$ so that $\rbt \ge r$ and
$\rtb\ge r$. Thus $\rtb=\rbt=r$ as desired.
Otherwise suppose Case 3 holds for $(\nu,J)^{(r-1)}$.
Recall that in Case 3 one has $\lbt^{(r-1)}=\ltb^{(r-1)}$.
If $\lbt^{(r-1)} < \infty$ then $\rbt=\rtb=r$ as desired.
Otherwise $\lbt^{(r-1)}=\ltb^{(r-1)}=\infty$ so that
$\rbt\le r-1$, $\rtb\le r-1$ and furthermore
$\la^t_{r-1}=\la^t_r$. Suppose that $\lbt^{(r-2)}=\infty$. Then it 
is clear from Lemma \ref{db and dt} that either $\lt^{(r-2)}=\infty$
or $\lt^{(r-1)}=\infty$ contradicting the definition of $r$.
Similarly $\ltb^{(r-2)}=\infty$ leads to a contradiction. Hence
$\rbt=\rtb=r-1$ as desired.

The rest of this section is devoted to the proofs of 
Lemmas \ref{db and dt} and \ref{label commute}.

\begin{proof}[Proof of Lemma \ref{db and dt}]
The proof proceeds by induction on $k$.
There is nothing to prove unless at least
one of $\lb^{(k)}$ and $\lt^{(k)}$ is finite.
If one is finite and the other infinite then obviously $\db$ and $\dt$
choose different strings and $\lbt^{(k)}=\lt^{(k)}$ and $\ltb^{(k)}=\lb^{(k)}$.
So it is assumed that both $\lb^{(k)}$ and $\lt^{(k)}$ are finite.

For the base case $k=0$, observe that the doubly singular
case holds precisely when $R=((1^{\mu_1}))$. Then $\RC(\la^t;R^t)$
is empty unless $\la=(1^{\mu_1})$, and in that case consists
of the empty rigged configuration $(\emptyset,\emptyset)$.  Then
$\lb^{(0)}=\lt^{(0)} = \mu_1$ and $\lbt^{(0)}=\ltb^{(0)}=\mu_1 - 1$.
But $\rkb(\emptyset,\emptyset)=\rkt(\emptyset,\emptyset)=1$
so this case is never used in the inductive step.
In the generic case $L \ge 2$, $\lbt^{(0)}=\lt^{(0)}=\mu_1$ and
$\ltb^{(0)}=\lb^{(0)}=\mu_L$.

Now assume $k \ge 1$. Note that
\begin{equation}\label{ineq for lb}
\begin{aligned}
\lb^{(k)}& \ge \ltb^{(k-1)}\\
\lt^{(k)}& \ge \lbt^{(k-1)}.
\end{aligned}
\end{equation}
If the generic case, Case 1 or Case 2 occurs at $k-1$, this follows 
immediately from $\lb^{(k)}\ge \lb^{(k-1)}\ge \ltb^{(k-1)}$ and
$\lt^{(k)}\ge \lt^{(k-1)}\ge \lbt^{(k-1)}$. For Case 3 at $k-1$
\eqref{ineq for lb} holds by induction hypothesis.

\subsection*{Generic case}

Observe that $\lt^{(k)} = \lbt^{(k)}$ is obtained from
$\lb^{(k)} = \ltb^{(k)}$ under conjugation by the involution $\comev_R$,
so we shall only prove the latter.
By Remark \ref{undisturbed},
the singular string in $\nu^{(k)}$ of length $\lb^{(k)}$
remains singular in passing to $\nut^{(k)}$.  Since
$\ltb^{(k-1)} \le \lb^{(k)}$ by \eqref{ineq for lb}, it follows that
$\ltb^{(k)} \le \lb^{(k)}$.

If $\ltb^{(k)}=\lb^{(k)}$ we are done. By induction hypothesis
$\ltb^{(k)}\ge \ltb^{(k-1)}\ge \lb^{(k-1)}-1$. Let us first assume
that $\lb^{(k-1)}\le \ltb^{(k)} < \lb^{(k)}$.
By Remark \ref{undisturbed} this is only possible if
the string selected by $\db$ acting on $\nut^{(k)}$
is the string shortened by $\dt$ acting on $\nu^{(k)}$.
This string in $\nut^{(k)}$ has length $\lt^{(k)}-1$ and label $0$.
We show that this cannot occur.  For this it is enough to show that
\begin{equation}\label{P generic}
P_{\lt^{(k)}-1}^{(k)}(\nut)>0 \quad \text{if} \quad
 \lb^{(k-1)}<\lt^{(k)}\leq \lb^{(k)} \quad \text{and} \quad 
\ltb^{(k-1)}<\lt^{(k)}.
\end{equation}
Otherwise we may assume that $\lb^{(k-1)}-1=\ltb^{(k)}<\lb^{(k)}$. This
means that Case 3 occurs at $k-1$ so that $m_{\lb^{(k-1)}-1}(\nu^{(k)})=0$
and $\lt^{(k-1)}=\lb^{(k-1)}$. Hence $\ltb^{(k)}=\lb^{(k-1)}-1$ can only
occur if $\lt^{(k)}=\lt^{(k-1)}=\lb^{(k-1)}$. To prove that this cannot
happen it suffices to show that
\begin{equation}\label{P Case 3}
P_{\lt^{(k)}-1}^{(k)}(\nut)>0 \quad \text{if $m_{\lt^{(k-1)}-1}(\nu^{(k)})=0$
and $\lb^{(k-1)}=\lt^{(k-1)}=\lt^{(k)}\le \lb^{(k)}$}.
\end{equation}

By \eqref{dt vacancy} with $n=\lt^{(k)}-1$,
\begin{equation} \label{v1}
   P^{(k)}_{\lt^{(k)}-1}(\nut) = P^{(k)}_{\lt^{(k)}-1}(\nu) -
	\chi(\lt^{(k-1)} \le \lt^{(k)} - 1).
\end{equation}
We proceed by considering various cases that exhaust all possibilities, and
show that \eqref{P generic} and \eqref{P Case 3} both hold in each case.

First assume that $m_{\lt^{(k)}-1}(\nu^{(k)})=0$ and
$\lt^{(k-1)}=\lt^{(k)}$.  
Since $\lb^{(k-1)} \le \lt^{(k)} \le \lb^{(k)}$ by assumption,
$P^{(k)}_{\lt^{(k)}}(\nu) \ge 1$, for otherwise
$\lb^{(k)}=\lt^{(k)}$ and $\db$ and $\dt$ would select the same string
in $\nu^{(k)}$.  By Lemma \ref{lower bound} with $n=\lt^{(k)} - 1$
it follows that $P^{(k)}_{\lt^{(k)}-1}(\nu)\ge 1$, which by \eqref{v1}
implies $P^{(k)}_{\lt^{(k)}-1}(\nut)\ge 1$.
This proves in particular \eqref{P Case 3}.

For all remaining cases \eqref{P Case 3} holds vacuously.
Hence it may be assumed that the hypotheses of \eqref{P generic} hold.

Assume that $m_{\lt^{(k)}-1}(\nu^{(k)})>0$ and
$\lt^{(k-1)}<\lt^{(k)}$.
Since $\lt^{(k-1)} \le \lt^{(k)} -1 < \lt^{(k)}$
and $\lb^{(k-1)} \le \lt^{(k)} - 1 < \lb^{(k)}$, there cannot be strings
in $\nu^{(k)}$ of length $\lt^{(k)}-1$ that have label zero or are singular.
Since there are strings in $\nu^{(k)}$ of length $\lt^{(k)}-1$,
there must be an available label that is neither zero nor maximum.
Thus $P_{\lt^{(k)}-1}^{(k)}(\nu)\geq 2$, and
by \eqref{v1}, $P_{\lt^{(k)}-1}^{(k)}(\nut)\ge 1$.

Now assume that $m_{\lt^{(k)}-1}(\nu^{(k)})>0$ and $\lt^{(k-1)}=\lt^{(k)}$.
Since $\lb^{(k-1)}<\lt^{(k)} \le \lb^{(k)}$ there cannot be a
singular string of length $\lt^{(k)}-1$ in $\nu^{(k)}$.
Hence $P^{(k)}_{\lt^{(k)}-1}(\nu)\ge 1$ and
by \eqref{v1} $P^{(k)}_{\lt^{(k)}-1}(\nut)\ge 1$.

Finally consider the case $m_{\lt^{(k)}-1}(\nu^{(k)})=0$ and
$\lt^{(k-1)}<\lt^{(k)}$.  By \eqref{v1} if
$P^{(k)}_{\lt^{(k)}-1}(\nu) \ge 2$ then we are done, so
it may be assumed that $P^{(k)}_{\lt^{(k)}-1}(\nu) = 1$.
Let $\ell$ be maximal such that
$\ell<\lt^{(k)}$ and $m_\ell(\nu^{(k)}) \ge 1$. If no such $\ell$ exists
set $\ell=0$.

Suppose that $P^{(k)}_{\lt^{(k)}}(\nu) = 0$.
By the definition of $\lt^{(k)}$, $m_{\lt^{(k)}}(\nu^{(k)}) \ge 1$.
Hence there is a string of length $\lt^{(k)}$ in $\nu^{(k)}$,
which is singular since its vacancy number is zero.
Due to the assumption that $\lb^{(k-1)} < \lt^{(k)} \le \lb^{(k)}$,
the definition of $\lb^{(k)}$ yields $\lb^{(k)}=\lt^{(k)}$.
This means that the doubly singular case holds for $\nu^{(k)}$,
which is a contradiction.

Suppose $P^{(k)}_{\lt^{(k)}}(\nu) = 1$.  By Lemma \ref{lower bound}
we have $P^{(k)}_n(\nu) = 1$ for $\ell < n \le \lt^{(k)}$ and
$0\le P_{\ell}^{(k)}(\nu)\le 1$.
By \eqref{vacancy ddiff}, $m_n(\nu^{(k-1)})=0$ for $\ell+1 < n < \lt^{(k)}$
and $0\le m_{\ell+1}(\nu^{(k-1)})\le 1$.

First consider $m_{\ell+1}(\nu^{(k-1)})=0$. If $\ell=0$ there is no string
of length smaller than $\lt^{(k)}$ in $\nu^{(k-1)}$ so that
$\lt^{(k-1)}\ge \lt^{(k)}$. This contradicts the assumptions. So assume that
$\ell>0$. Since $\lt^{(k-1)}\le \ell$ and $\lb^{(k-1)}\le \ell$
this requires that the strings of length $\ell$ in $\nu^{(k)}$ are neither
singular nor have label zero so that $P_{\ell}^{(k)}(\nu)\ge 2$. This is
a contradiction.

Hence assume that $m_{\ell+1}(\nu^{(k-1)})=1$. By \eqref{vacancy ddiff}
at $n=\ell+1$ this implies that $P_{\ell}^{(k)}(\nu)=0$. Since 
$\lt^{(k)}>\ell$ and $\lb^{(k)}>\ell$ this requires that 
$\lt^{(k-1)}=\lb^{(k-1)}=\ell+1$ and hence $P_{\ell+1}^{(k-1)}(\nu)=0$.
Since $m_{\ell+1}(\nu^{(k-1)})=1$, Case 1 or Case 3 occurs at $k-1$.
If Case 3 occurs, $\lbt^{(k-1)}=\ltb^{(k-1)}\ge \lt^{(k)}$ which
contradicts the assumptions. Therefore $\nu^{(k-1)}$ must be in Case 1
and $\ell>0$. By induction $m_{\ell}(\nu^{(k)})=0$ which contradicts the
definition of $\ell$.

Now suppose $P^{(k)}_{\lt^{(k)}}(\nu) = 2$.  Again~\eqref{vacancy ddiff}
and Lemma~\ref{lower bound} fail unless $\lt^{(k)}=\ell+2$ and 
$P^{(k)}_\ell(\nu)=0$.
By \eqref{vacancy ddiff} with $n=\ell+1=\lt^{(k)}-1$,
$m_{\ell+1}(\nu^{(k-1)}) = 0$.  If $\ell=0$ then
$\lt^{(k-1)} < \lt^{(k)} = 2$ which forces $\lt^{(k-1)}=1$,
but there is no string of length $1$ in $\nu^{(k-1)}$, which
is a contradiction.  So suppose $\ell>0$.  Since
$\lb^{(k-1)}<\lt^{(k)}$, $\lt^{(k-1)}<\lt^{(k)}$
and $m_{\lt^{(k)}-1}(\nu^{(k-1)}) = 0$, one finds
$\lb^{(k-1)} \le \ell$ and $\lt^{(k-1)} \le \ell$.
Also there is a string of length $\ell$ in $\nu^{(k)}$,
which is both singular and has label zero since
$P^{(k)}_\ell(\nu) = 0$.  But $\ell < \lt^{(k)} \le \lb^{(k)}$,
which contradicts the definition of $\lb^{(k)}$ and $\lt^{(k)}$.

If $P^{(k)}_{\lt^{(k)}}(\nu) > 2$ then there is an immediate
contradiction since \eqref{convex} fails for $n=\lt^{(k)}-1$.

This completes the proof of \eqref{P generic} and \eqref{P Case 3}
and hence the proof of the generic case for $\nu^{(k)}$.

\subsection*{Doubly singular case}
Since there is a string of length $\ell$ in $\nu^{(k)}$ that is both
singular and has label zero, it must have vacancy number zero, that is,
$P_{\ell}^{(k)}(\nu)=0$.

\subsection*{Case 1: $\lbt^{(k)} < \ell$}

By induction we have $\lbt^{(k)}\ge \lbt^{(k-1)}\ge \lt^{(k-1)}-1$.
First assume that $\lt^{(k-1)} \le \lbt^{(k)} < \ell$. In the light
of Remark \ref{undisturbed}, $\dt$ must select the string in $\nub^{(k)}$
that was shortened by $\db$ in the transformation $(\nu,J)\mapsto (\nub,\Jb)$,
so that $\lbt^{(k)} = \ell - 1$.  This string in $\nub^{(k)}$ is singular
since it was shortened by $\db$ and has label zero since it is selected by
$\dt$, so $P^{(k)}_{\ell-1}(\nub)=0$.
The case $\lt^{(k-1)}-1=\lbt^{(k)}<\ell$ can only occur for Case 1 at
$k-1$. By induction this implies that $m_{\lt^{(k-1)}-1}(\nu^{(k)})=0$.
For $\lt^{(k-1)}-1=\lbt^{(k)}$ one needs $m_{\lt^{(k-1)}-1}(\nub^{(k)})>0$
so that $\ell=\lt^{(k-1)}$. Hence $\lbt^{(k)}=\ell-1$ and
$P^{(k)}_{\ell-1}(\nub)=0$ as before.

The goal is to show that $\ltb^{(k)} = \ell - 1$.
Since $\lt^{(k)} = \ell$, it follows that $m_{\ell-1}(\nut^{(k)}) \ge 1$.
It suffices to show that $\ltb^{(k-1)} \le \ell - 1$
and $P^{(k)}_{\ell-1}(\nut) = 0$.  For then
$\ltb^{(k)} < \ell$, and by the same arguments as before
it follows that $\ltb^{(k)} = \ell - 1$.

By \eqref{db vacancy}, \eqref{dt vacancy}, and $P^{(k)}_{\ell-1}(\nub)=0$,
\begin{equation} \label{theeq}
\begin{split}
  P^{(k)}_{\ell-1}(\nu) &= P^{(k)}_{\ell-1}(\nut) +
  	\chi(\lt^{(k-1)} \le \ell - 1) \\
  &= \chi(\lb^{(k-1)} \le \ell - 1).
\end{split}
\end{equation}

Suppose that $\ltb^{(k-1)} \ge \ell$.
Now $\lbt^{(k-1)} \le \lbt^{(k)} = \ell - 1$
so $\lbt^{(k-1)} \not= \ltb^{(k-1)}$.
By induction, the generic case holds for $\nu^{(k-1)}$, and
$\lbt^{(k-1)} = \lt^{(k-1)}$ and
$\ltb^{(k-1)} = \lb^{(k-1)}$.  So
$\lb^{(k-1)} = \ltb^{(k-1)} \ge \ell$
and $\lt^{(k-1)} = \lbt^{(k-1)} \le \ell -1$.
This leads to a contradiction in evaluating
\eqref{theeq}, so $\ltb^{(k-1)} \le \ell - 1$.

Suppose $P^{(k)}_{\ell-1}(\nut) \ge 1$.  Then
by \eqref{theeq}, $\lt^{(k-1)} \ge \ell$ and
$\lb^{(k-1)} \le \ell - 1$.  
Since $\lt^{(k-1)}\not=\lb^{(k-1)}$, by induction
the generic case holds for $\nu^{(k-1)}$
and $\lbt^{(k-1)} = \lt^{(k-1)} \ge \ell$,
which contradicts $\ell > \lbt^{(k)} \ge \lbt^{(k-1)}$.

Therefore $P_{\ell-1}^{(k)}(\nut)=0$ and $\ltb^{(k)}=\lbt^{(k)}=\ell-1$.

To finish Case 1 it suffices to show that
$m_{\ell-1}(\nu^{(k+1)}) = 0$.

Since it has been shown that $P^{(k)}_{\ell-1}(\nut)=0$,
\eqref{theeq} becomes
\begin{equation} \label{neweq}
  P^{(k)}_{\ell-1}(\nu) = \chi(\lt^{(k-1)} \le \ell-1)
  =\chi(\lb^{(k-1)} \le \ell-1).
\end{equation}

Suppose that $P^{(k)}_{\ell-1}(\nu) = 0$.
By \eqref{neweq}, $\lt^{(k-1)}=\ell$ and $\lb^{(k-1)}=\ell$.
Now $\lbt^{(k-1)} \le \lbt^{(k)} = \ell - 1 < \lt^{(k-1)}$,
so Case 1 holds for $\nu^{(k-1)}$.  By induction
$\lbt^{(k-1)} = \lt^{(k-1)} - 1$ and
$m_{\lt^{(k-1)}-1}(\nu^{(k)})=0$, that is,
$m_{\ell-1}(\nu^{(k)}) = 0$.
By \eqref{vacancy ddiff} for $n=\ell-1$, it follows that
$m_{\ell-1}(\nu^{(k+1)}) = 0$, as desired.

Otherwise $P^{(k)}_{\ell-1}(\nu) = 1$.  Here $\ell\ge 2$.
By \eqref{neweq}, $\lt^{(k-1)} \le \ell-1$ and $\lb^{(k-1)} \le \ell-1$.
By the minimality of $\lt^{(k)}$ and $\lb^{(k)}$,
there cannot be strings in $\nu^{(k)}$ of length $\ell-1$ that
are singular or have label zero, so $m_{\ell-1}(\nu^{(k)}) = 0$.
Applying \eqref{vacancy ddiff} at $n=\ell-1$ and using the fact that
$P^{(k)}_\ell(\nu) = 0$ (since the doubly singular case holds for
$\nu^{(k)}$) one obtains
\begin{equation} \label{NPinq}
  P^{(k)}_{\ell-2}(\nu) + m_{\ell-1}(\nu^{(k-1)}) +
  m_{\ell-1}(\nu^{(k+1)}) \le 2.
\end{equation}

If $m_{\ell-1}(\nu^{(k-1)}) = 2$ then by
\eqref{NPinq} $m_{\ell-1}(\nu^{(k+1)})=0$ as desired.

Suppose $m_{\ell-1}(\nu^{(k-1)})=1$.  
If $P^{(k)}_{\ell-2}(\nu) = 1$ then again we conclude
$m_{\ell-1}(\nu^{(k+1)})=0$ by \eqref{NPinq}.
So assume that $P^{(k)}_{\ell-2}(\nu) = 0$.
If $m_{\ell-2}(\nu^{(k)})=0$ then by
\eqref{vacancy ddiff} with $n=\ell-2$,
$P^{(k)}_{\ell-1}(\nu)=0$ which is a contradiction.
Otherwise $m_{\ell-2}(\nu^{(k)})\ge 1$.  Then there
is a singular string of length $\ell-2$ in $\nu^{(k)}$.
By the definition of $\lb^{(k)}$ we have
$\lb^{(k-1)} > \ell-2 $, that is, $\lb^{(k-1)} = \ell-1$.
Similarly $\lt^{(k-1)}=\ell-1$.
Since $m_{\ell-1}(\nu^{(k-1)})=1$, by induction
the doubly singular case holds for $\nu^{(k-1)}$.
Now $\lbt^{(k-1)} \le \lbt^{(k)}=\ell-1=\lt^{(k-1)}$,
so Case 3 is impossible.  Since $m_{\ell-1}(\nu^{(k-1)})=1$,
Case 2 is also impossible.  So Case 1 holds for $\nu^{(k-1)}$.
It follows that $m_{\lt^{(k-1)}-1}(\nu^{(k)})=0$,
that is, $m_{\ell-2}(\nu^{(k)})=0$.  But this is a contradiction.

Suppose that $m_{\ell-1}(\nu^{(k-1)}) = 0$.
Then $\lb^{(k-1)} \le \ell - 2$ and $\lt^{(k-1)} \le \ell - 2$.
This yields a contradiction unless $\ell > 2$.
By the minimality of $\lb^{(k)}$ and $\lt^{(k)}$, there
cannot be strings in $\nu^{(k)}$ of length $\ell-2$
that are either singular or have label zero, so it follows that
either $m_{\ell-2}(\nu^{(k)}) = 0$ or $P^{(k)}_{\ell-2}(\nu) \ge 2$.
Using \eqref{NPinq} the latter immediately yields the desired result
$m_{\ell-1}(\nu^{(k+1)})=0$,
so assume that $m_{\ell-2}(\nu^{(k)}) = 0$ and
$P^{(k)}_{\ell-2}(\nu) \le 1$.

If $P^{(k)}_{\ell-2}(\nu)=0$,
by \eqref{vacancy ddiff} with $n=\ell-2$ it follows that
$P^{(k)}_{\ell-1}(\nu)=0$, which is a contradiction.

So $P^{(k)}_{\ell-2}(\nu) = 1$.

Let $p<\ell-1$ be maximal such that $m_p(\nu^{(k)})\ge 1$;
if no such $p$ exists set $p=0$.  Then by Lemma \ref{lower bound},
$P^{(k)}_n(\nu)=1$ for $p < n < \ell$ and $P^{(k)}_p(\nu) \le 1$.
By \eqref{vacancy ddiff} it follows that
$m_n(\nu^{(k-1)})=0$ for $p+1 < n < \ell-1$
and $m_{\ell-1}(\nu^{(k-1)})+m_{\ell-1}(\nu^{(k+1)})\le 1$. 
If $m_{\ell-1}(\nu^{(k-1)})=1$ then $m_{\ell-1}(\nu^{(k+1)})=0$ as desired.
Hence assume that $m_{\ell-1}(\nu^{(k-1)})=0$.

Suppose $P^{(k)}_p(\nu)=1$ which implies $p>0$. Then by
\eqref{vacancy ddiff} it follows that $m_{p+1}(\nu^{(k-1)})=0$.
Since $\lt^{(k-1)} < \ell-1$ and $\lb^{(k-1)} < \ell-1$, it
follows that $\lt^{(k-1)} \le p$ and $\lb^{(k-1)} \le p$.
Also $m_p(\nu^{(k)}) \ge 1$ since $p>0$.
Then there is a singular string of length $p$ in $\nu^{(k)}$
or one of label zero, contradicting the definition of $\lb^{(k)}$
or $\lt^{(k)}$.

Therefore $P^{(k)}_p(\nu)=0$.
By \eqref{vacancy ddiff} for $n=p+1$,
we have $m_{p+1}(\nu^{(k-1)})+m_{p+1}(\nu^{(k+1)})\le 1$.

Suppose that $m_{p+1}(\nu^{(k-1)}) = 0$.  If $p>0$ then the
proof proceeds as before.  If $p=0$ then there are no strings
in $\nu^{(k-1)}$ of length less than $\ell$.
This contradicts $\lt^{(k-1)} \le \ell - 1$.

So assume that $m_{p+1}(\nu^{(k-1)})=1$.

Certainly $\lt^{(k-1)} \le p+1$ and $\lb^{(k-1)} \le p+1$.
If either $\lt^{(k-1)} \le p$ or $\lb^{(k-1)}\le p$ then
there is a string in $\nu^{(k)}$
of length $p$ that (due to
$P^{(k)}_p(\nu)=0$) is either singular or has label zero.
But $p < \ell$, contradicting the definition of $\lb^{(k)}$ or
$\lt^{(k)}$.  

So $\lt^{(k-1)} = p+1$ and $\lb^{(k-1)} = p+1$.
Since $m_{p+1}(\nu^{(k-1)})=1$,
$m_n(\nut^{(k-1)})=0$ for $p < n < \ell$.

Now $\ltb^{(k-1)} \le \ltb^{(k)} = \ell-1$,
so $\ltb^{(k-1)} \le p < p+1 = \lb^{(k-1)}$.
So Case 1 holds for $\nu^{(k-1)}$.
By induction $\ltb^{(k-1)}=\lb^{(k-1)}-1=p$
and $m_p(\nu^{(k)})=0$ which contradicts the construction of $p$.

This concludes the proof of Case 1.

Using the involution $\comev_R$, the above argument also shows
that if $\ltb^{(k)} < \ell$ then $\ltb^{(k)}=\lbt^{(k)}=\ell-1$.

\subsection*{Case 2: $\lbt^{(k)} = \ell$}
It will be shown that $\ltb^{(k)} = \ell$.
By Case 1, $\ltb^{(k)} \ge \ell$.

The assumption $\lbt^{(k)}=\ell$ implies that
$m_\ell(\nub^{(k)}) \ge 1$.  Since $\lb^{(k)}=\ell$
one part of size $\ell$ is shortened in passing from
$\nu^{(k)}$ to $\nub^{(k)}$, so that $m_\ell(\nu^{(k)}) \ge 2$.
Now $P^{(k)}_\ell(\nu) = 0$, so there is at least one singular
string in $\nu^{(k)}$ that is not selected by $\dt$ acting on
$(\nu,J)$.  By Remark \ref{undisturbed} this string remains singular
of length $\ell$ in passing to $\nut^{(k)}$.  This shows that
there is a singular string of length $\ell$ in $\nut^{(k)}$.
Thus to prove $\ltb^{(k)} = \ell$ it suffices to show
that $\ltb^{(k-1)} \le \ell$.  

If $\ltb^{(k-1)} \le \lb^{(k-1)}$, then
$\ltb^{(k-1)} \le \lb^{(k-1)} \le \lb^{(k)} = \ell$.
Otherwise $\ltb^{(k-1)} > \lb^{(k-1)}$, so Case 3
holds for $\nu^{(k-1)}$, and by induction
$\ltb^{(k-1)}=\lbt^{(k-1)}\le\lbt^{(k)}=\ell$.

By applying $\comev_R$ this also shows that
if $\ltb^{(k)} = \ell$ then $\lbt^{(k)} = \ell$.

\subsection*{Case 3: $\lbt^{(k)} > \ell$}

By Cases 1 and 2 it follows that $\ltb^{(k)} > \ell$.
It will be shown that $\ltb^{(k)} = \lbt^{(k)}$ and that
\begin{align}
\label{cond1}
P_{\ell+1}^{(k)}(\nu)&=0\\
\label{cond2}
P_{\ell-1}^{(k)}(\nu)&=2-m_\ell(\nu^{(k-1)}) \\
\label{cond3}
m_\ell(\nu^{(k+1)}) &= 0.
\end{align}

Suppose $m_\ell(\nu^{(k)}) \ge 2$.  Since
$P^{(k)}_\ell(\nu) = 0$, there is a string in
$\nu^{(k)}$ of length $\ell$ and label zero that
is not selected by $\db$ in passing to $\nub^{(k)}$.
By Remark \ref{undisturbed} it follows that
there is a string in $\nub^{(k)}$ of length $\ell$
with label zero.  Since $\lbt^{(k)} > \ell$, this
string cannot be selected, that is,
$\lbt^{(k-1)} > \ell$.  Now
$\lt^{(k-1)} \le \lt^{(k)} = \ell < \lbt^{(k-1)}$.
By induction Case 3 holds at $k-1$. This implies in particular
$\lt^{(k)} \ge \lbt^{(k-1)} > \ell$, which is a
contradiction.

Therefore $m_\ell(\nu^{(k)}) = 1$.
By \eqref{vacancy ddiff} with $n=\ell$,
$P^{(k)}_\ell(\nu)=0$, and $m_\ell(\nu^{(k)})=1$, we have
\begin{equation} \label{Pinq}
  P^{(k)}_{\ell-1}(\nu) + P^{(k)}_{\ell+1}(\nu) +
  m_\ell(\nu^{(k-1)}) + m_\ell(\nu^{(k+1)}) \le 2.
\end{equation}

We distinguish the three cases $m_\ell(\nu^{(k-1)})\in\{0,1,2\}$.

We start with $m_\ell(\nu^{(k-1)})=0$.  Recall that
$\lb^{(k-1)} \le \lb^{(k)} = \ell$ and
$\lt^{(k-1)} \le \lt^{(k)} = \ell$.  However the inequalities
must be strict since there are no strings of length $\ell$ in
$\nu^{(k-1)}$.  So
\begin{equation} \label{small length}
\begin{split}
  \lb^{(k-1)} &\le \ell - 1 \\
  \lt^{(k-1)} &\le \ell - 1.
\end{split}
\end{equation}

If $m_{\ell-1}(\nu^{(k)})>0$ then necessarily 
$P_{\ell-1}^{(k)}(\nu)\ge 2$ since otherwise a string in
$\nu^{(k)}$ of length $\ell-1$ is selected by $\db$ or $\dt$,
contradicting $\lb^{(k)}=\lt^{(k)}=\ell$.
By \eqref{Pinq} we conclude that
$P_{\ell-1}^{(k)}(\nu)=2$, $P^{(k)}_{\ell+1}(\nu)=0$, and
$m_\ell(\nu^{(k+1)})=0$.

Now assume $m_{\ell-1}(\nu^{(k)})=0$.  By \eqref{db vacancy} with
$n=\ell-1$, $\lb^{(k-1)}\le \ell-1$, and $\lb^{(k)}=\ell$, we have
\begin{equation} \label{dv}
  P^{(k)}_{\ell-1}(\nu) = P^{(k)}_{\ell-1}(\nub) + 1.
\end{equation}

Suppose that $P^{(k)}_{\ell-1}(\nub)=0$.  Since $m_{\ell-1}(\nub^{(k)}) = 1$
there is a string of length $\ell-1$ and label zero in $\nub^{(k)}$.
But $\lbt^{(k)} > \ell$, so the only way that this string is not selected
is if $\lbt^{(k-1)} > \ell - 1$.  But by \eqref{small length}
$\ell - 1 \ge \lt^{(k-1)}$, so $\lbt^{(k-1)} > \lt^{(k-1)}$.  This is Case 3
for $\nu^{(k-1)}$.  By induction
$\lbt^{(k-1)} = \ltb^{(k-1)} > \lt^{(k-1)}=\lb^{(k-1)}$
and $\lt^{(k)} \ge \lbt^{(k-1)}$.  But $\lt^{(k)}=\ell$
and $\lbt^{(k-1)} \ge \ell$ so $\lbt^{(k-1)} = \ell$.
Then $m_\ell(\nub^{(k-1)}) \ge 1$.  However $\lb^{(k-1)} \le \ell - 1$
by \eqref{small length} so $m_\ell(\nu^{(k-1)}) \ge 1$,
which is a contradiction.

Therefore $P^{(k)}_{\ell-1}(\nub) \ge 1$.
By \eqref{dv}, $P^{(k)}_{\ell-1}(\nu) \ge 2$, so by
\eqref{Pinq}, $P^{(k)}_{\ell-1}(\nu)=2$, $P_{\ell+1}^{(k)}(\nu)=0$ and
$m_{\ell}(\nu^{(k+1)})=0$ as above.

Next consider $m_\ell(\nu^{(k-1)})=1$.
Suppose that $P^{(k)}_{\ell-1}(\nu)=0$.
By \eqref{db vacancy}
\begin{equation} \label{v3}
P^{(k)}_{\ell-1}(\nu)=P^{(k)}_{\ell-1}(\nub)+\chi(\lb^{(k-1)}\le \ell-1).
\end{equation}
Therefore $\lb^{(k-1)} \ge \ell = \lb^{(k)} \ge \lb^{(k-1)}$,
so that $\lb^{(k-1)}=\ell$.  Similarly $\lt^{(k-1)}=\ell$.
Then in $\nu^{(k-1)}$ there is only one string of length $\ell$
and both $\db$ and $\dt$ select it.
So for $\nu^{(k-1)}$ we are in Case 1 or Case 3.
Suppose it is Case 3.  Then
$\lbt^{(k-1)} = \ltb^{(k-1)} > \lt^{(k-1)} = \lb^{(k-1)} = \ell$
and $\lt^{(k)} \ge \lbt^{(k-1)} > \ell$ which is a
contradiction.  Suppose it is Case 1.
Then $\lbt^{(k-1)} = \lt^{(k-1)} - 1 = \ell - 1$.
By \eqref{v3}, $P^{(k)}_{\ell-1}(\nub) = 0$.
Since $\lb^{(k)}=\ell$, $m_{\ell-1}(\nub^{(k)}) \ge 1$.
So there is a singular string of length $\ell-1$ in $\nub^{(k)}$
and $\lbt^{(k-1)}=\ell-1$.  By definition,
$\lbt^{(k)} = \ell-1$, which contradicts the assumption $\lbt^{(k)}>\ell$.

Therefore $P^{(k)}_{\ell-1}(\nu) \ge 1$.  By \eqref{Pinq} it follows
that $P^{(k)}_{\ell-1}(\nu) = 1$, $P^{(k)}_{\ell+1}(\nu)=0$ and
$m_\ell(\nu^{(k+1)})=0$.

Finally consider $m_\ell(\nu^{(k-1)})=2$.  By \eqref{Pinq},
$P_{\ell-1}^{(k)}(\nu)=0$, $P_{\ell+1}^{(k)}(\nu)=0$,
and $m_\ell(\nu^{(k+1)})=0$.

So \eqref{cond1}, \eqref{cond2}, and \eqref{cond3} are proven.

Since $P^{(k)}_\ell(\nu)=P^{(k)}_{\ell+1}(\nu)=0$ (see \eqref{cond1})
it follows from Lemma \ref{lower bound}, that if
$p>\ell$ and $m_n(\nu^{(k)})=0$ for all $\ell < n < p$,
then $P^{(k)}_n(\nu)=0$ for $\ell \le n \le p$. Equation \eqref{vacancy ddiff}
furthermore implies that $m_n(\nu^{(k+1)}) = 0$ for $\ell < n < p$.

Suppose $\nu^{(k)}$ has a string longer than $\ell$.
Let $p$ be minimal such that $p>\ell$ and $m_p(\nu^{(k)}) \ge 1$.
Since the string in $\nu^{(k)}$ of length $p$
is selected by neither $\db$ nor $\dt$ but has vacancy number zero,
its length remains $p$, its label stays zero in $\nub^{(k)}$,
and it remains singular in $\nut^{(k)}$.
Neither $\nub^{(k)}$ nor $\nut^{(k)}$ have strings of length $n$ for
$\ell \le n < p$. Since $\lbt^{(k-1)}\le \ell$ and $\ltb^{(k-1)}\le\ell$
by~\eqref{ineq for lb} and since by assumption $\lbt^{(k)} > \ell$ and
$\ltb^{(k)} > \ell$, it follows that $\lbt^{(k)} = p = \ltb^{(k)}$.
Moreover, by the previous paragraph and \eqref{cond3},
$m_n(\nu^{(k+1)})=0$ for $\lt^{(k)}=\lb^{(k)}=\ell \le n < p$, so that
$\lt^{(k+1)} \ge p$ and $\lb^{(k+1)} \ge p$.

Otherwise there is no string in $\nu^{(k)}$ longer than $\ell$.
Then $\lbt^{(k)}=\ltb^{(k)} = \infty$ and $\rkb(\nut,\Jt)=
\rkt(\nub,\Jb)=k$.
Moreover, the above result holds for all $p>\ell$, so that
$m_n(\nu^{(k+1)})=0$ for $n>\ell$.  But by \eqref{cond3}
$m_\ell(\nu^{(k+1)})=0$ so that $\lt^{(k+1)}=\lb^{(k+1)}=\infty$
and $\rkb(\nu,J)=\rkt(\nu,J)=k+1$.

By the appendix in \cite{KS} it follows that
\begin{equation*}
  \la^t_k-\la^t_{k+1}=\lim_{n\rightarrow\infty} P^{(k)}_n(\nu).
\end{equation*}
But the right-hand side is zero so that $\la^t_k=\la^t_{k+1}$.
\end{proof}

\begin{proof}[Proof of Lemma \ref{label commute}]
$\nubt=\nutb$ by Lemma \ref{db and dt},
whose entire proof will be used repeatedly without additional mention.

\subsection*{Selected strings}
Consider a string in $(\nu,J)^{(k)}$
that is either selected by $\db$ or $\dt$,
or is such that its image under $\db$  (resp. $\dt$) is selected
by $\dt$ (resp. $\db$).  It is shown that
the image of any such string under both $\dt\circ\db$ and
$\db\circ\dt$, has the same label.

\subsection*{Selected strings, generic case}
In the generic case for $\nu^{(k)}$,
Remark \ref{undisturbed} shows that the string
$(\lt^{(k)},0)\in (\nu,J)^{(k)}$ is sent to the string
$(\lt^{(k)}-1,0)$ under either $\db \circ \dt$ or $\dt\circ \db$,
and that a singular string of length $\lb^{(k)}$ in
$(\nu,J)^{(k)}$ is sent to a singular string of length
$\lb^{(k)}-1$ under either $\db\circ \dt$ or $\dt\circ \db$.

\subsection*{Selected strings, doubly singular case}
Write $\ell=\lt^{(k)}=\lb^{(k)}$.
The string $(\ell,0)\in(\nu,J)^{(k)}$ is also singular.

\subsection*{Selected strings, Case 1}
Here only the string $(\ell,0)\in (\nu,J)^{(k)}$
and its images under $\db$ and $\dt$ are selected.
Moreover $P^{(k)}_{\ell-1}(\nub)=P^{(k)}_{\ell-1}(\nut)=0$
and $\lbt^{(k)}=\ltb^{(k)}=\ell-1$.
The string $(\ell,0)$ is sent to a string of length $\ell-2$
and singular label under $\db \circ \dt$ and
zero label under $\dt\circ \db$.
Hence it must be shown that $P^{(k)}_{\ell-2}(\nutb)=0$.
Applying \eqref{db vacancy} and \eqref{dt vacancy},
\begin{equation} \label{ds1v}
  P^{(k)}_{\ell-2}(\nutb)=P^{(k)}_{\ell-2}(\nu)
-\chi(\lt^{(k-1)} \le \ell-2)-\chi(\ltb^{(k-1)}\le \ell-2).
\end{equation}
We divide into cases as in the proof of Case 1 in Lemma \ref{db and dt}.
Suppose first that $P^{(k)}_{\ell-1}(\nu)=0$.  Then
$m_{\ell-1}(\nu^{(k)})=0$, and applying \eqref{vacancy ddiff}
with $n=\ell-1$ one obtains $P^{(k)}_{\ell-2}(\nu) = 0$.
The admissibility of $\nutb$ and
\eqref{ds1v} imply $P^{(k)}_{\ell-2}(\nutb)=0$.

By \eqref{neweq} the only alternative is
$P^{(k)}_{\ell-1}(\nu)=1$, which implies that $\lt^{(k-1)} \le \ell-1$
and $\lb^{(k-1)} \le \ell-1$.
By \eqref{NPinq} $m_{\ell-1}(\nu^{(k-1)}) \le 2$.
As in the proof of Case 1 in Lemma \ref{db and dt},
we distinguish the three cases given by
$m_{\ell-1}(\nu^{(k-1)}) \in \{0,1,2\}$.

Suppose $m_{\ell-1}(\nu^{(k-1)})=2$.  By \eqref{NPinq}
and \eqref{ds1v} $P^{(k)}_{\ell-2}(\nu)=0$ and $P^{(k)}_{\ell-2}(\nutb)=0$.

Suppose $m_{\ell-1}(\nu^{(k-1)})=1$.
By \eqref{NPinq} $P^{(k)}_{\ell-2}(\nu) \le 1$.
By \eqref{ds1v} it suffices to show that either
$\lt^{(k-1)} \le \ell-2$ or $\ltb^{(k-1)} \le \ell-2$.
Suppose neither holds.  Then
$\lt^{(k-1)} > \ell-2$ and $\ltb^{(k-1)} > \ell-2$.
Thus $\lt^{(k-1)}=\ell-1$.
Since $m_{\ell-1}(\nu^{(k-1)})=1$ and
$\lt^{(k-1)}=\ell-1$, it follows that
$m_{\ell-1}(\nut^{(k-1)})=0$.  But then
$\ell-1=\ltb^{(k)}\ge \ltb^{(k-1)}> \ell-2$, so
$\ltb^{(k-1)}=\ell-1$.  However there are no
strings of length $\ell-1$ in $\nut^{(k-1)}$, which
is a contradiction.

Suppose $m_{\ell-1}(\nu^{(k-1)})=0$.
By \eqref{NPinq} $P^{(k)}_{\ell-2}(\nu) \le 2$, so by \eqref{ds1v} 
it is enough to show that
$\lt^{(k-1)} \le \ell-2$ and $\ltb^{(k-1)}\le \ell-2$.
In this subcase $\lt^{(k-1)} \le \ell-2$ and $\lb^{(k-1)}\le \ell-2$.
Suppose $\ltb^{(k-1)} > \ell-2$.
Now $\ell-1=\ltb^{(k)}\ge\ltb^{(k-1)} > \ell-2$ so
$\ltb^{(k-1)}=\ltb^{(k)}=\ell-1$.
But $m_{\ell-1}(\nu^{(k-1)})=0$ and $\lt^{(k-1)} \le \ell-2$,
so $m_{\ell-1}(\nut^{(k-1)})=0$, contradicting $\ltb^{(k-1)} = \ell-1$.

\subsection*{Selected strings, Case 2}
Here there are two copies of the singular string
$(\ell,0)\in (\nu,J)^{(k)}$.  If we think of $\db$ 
as selecting one of them and $\dt$ the other,
then the proof is the same as in the generic case.

\subsection*{Selected strings, Case 3}
Let $p=\lbt^{(k)}=\ltb^{(k)}$.  It satisfies
$p > \ell$ and $P^{(k)}_p(\nu) = 0$.
Moreover $\lt^{(k+1)} \ge p$ and $\lb^{(k+1)} \ge p$.
The strings that will be selected are singular strings
$(\ell,0)$ and $(p,0)$ in $(\nu,J)^{(k)}$.

The string $(\ell,0)$ maps to a string of length $\ell-1$,
with label zero under $\db\circ \dt$ and singular label
under $\dt\circ \db$.
The string $(p,0)$ is sent to a string of length $p-1$,
with singular label under $\db\circ \dt$ and 
to zero label under $\dt\circ \db$.

It must be shown that
\begin{align}
\label{ds3v1}
  P^{(k)}_{\ell-1}(\nubt)&=0 \\
\label{ds3v2}
  P^{(k)}_{p-1}(\nutb)&=0.
\end{align}

First \eqref{ds3v1} is established.
By \eqref{db vacancy} and \eqref{dt vacancy}
\begin{equation} \label{ds3v}
  P^{(k)}_{\ell-1}(\nubt) =
  P^{(k)}_{\ell-1}(\nu) - \chi(\lbt^{(k-1)} \le \ell-1)
  -\chi(\lb^{(k-1)} \le \ell-1).
\end{equation}
By \eqref{cond2} $m_\ell(\nu^{(k-1)}) \le 2$.
We divide into cases for the choices of
$m_\ell(\nu^{(k-1)}) \in \{0,1,2\}$.

Suppose $m_\ell(\nu^{(k-1)})=2$.  By \eqref{cond2}
$P^{(k)}_{\ell-1}(\nu) = 0$.  It follows immediately
from the admissibility of $\nubt$ and \eqref{ds3v} that
$P^{(k)}_{\ell-1}(\nubt) = 0$.

Suppose $m_\ell(\nu^{(k-1)})=1$.
By \eqref{cond2} $P^{(k)}_{\ell-1}(\nu) = 1$,
so by \eqref{ds3v} it is enough to show that either
$\lb^{(k-1)} \le \ell-1$ or $\lbt^{(k-1)} \le \ell-1$.
Suppose neither holds.  Then
$\ell=\lb^{(k)} \ge \lb^{(k-1)} > \ell-1$ so $\lb^{(k-1)}=\ell$.
Also by \eqref{ineq for lb}
$\ell=\lt^{(k)} \ge \lbt^{(k-1)} > \ell-1$ so that
$\lbt^{(k-1)} = \ell$.  Now $m_\ell(\nu^{(k-1)})=1$
and $\lb^{(k-1)} = \ell$ so $m_\ell(\nub^{(k-1)})=0$,
contradicting $\lbt^{(k-1)}=\ell$.

Suppose $m_\ell(\nu^{(k-1)}) = 0$.
By \eqref{cond2} $P^{(k)}_{\ell-1}(\nu) = 2$,
so by \eqref{ds3v}, to prove \eqref{ds3v1} it is enough to show that
$\lb^{(k-1)} \le \ell-1$ and $\lbt^{(k-1)} \le \ell-1$.
By \eqref{small length}
we have $\lb^{(k-1)} \le \ell-1$ and $\lt^{(k-1)} \le \ell-1$.
Suppose $\lbt^{(k-1)} > \ell-1$.  
Since $\lbt^{(k-1)} \le \lt^{(k)}=\ell$
by \eqref{ineq for lb}, it follows that
$\lbt^{(k-1)} = \ell$.
Now $m_\ell(\nu^{(k-1)})=0$ and $\lb^{(k-1)} \le \ell-1$
so $m_\ell(\nub^{(k-1)})=0$.
This contradicts $\lbt^{(k-1)}=\ell$.

Now let us prove \eqref{ds3v2}.  Using
$p = \lbt^{(k)} > \lt^{(k)}$
and $\lt^{(k+1)} \ge \lbt^{(k)} = p$, by
\eqref{db vacancy} and \eqref{dt vacancy} we have
\begin{equation*}
  P^{(k)}_{p-1}(\nutb) 
  = P^{(k)}_{p-1}(\nu) - \chi(\ltb^{(k-1)} \le p-1) + 1.
\end{equation*}
Since $P^{(k)}_{p-1}(\nu) = 0$, it must be shown that
$\ltb^{(k-1)} \le p-1$.  This is certainly the case,
for otherwise by \eqref{ineq for lb},
$\ell=\lb^{(k)} \ge \ltb^{(k-1)} \ge p$, which is a
contradiction.

\subsection*{Unselected strings}
For the rest of the proof, assume that
$(n,x)$ is a string in $(\nu,J)^{(k)}$ that
is not selected by $\dt$ or $\db$, and is such that
its image under $\db$ (resp. $\dt$) is not selected
by $\dt$ (resp. $\db$).

Using the fact that $\db$ preserves labels and
$\dt$ preserves colabels, it is enough to show that
$P^{(k)}_n(\nu) - P^{(k)}_n(\nut) =
P^{(k)}_n(\nub) - P^{(k)}_n(\nubt)$,
which by \eqref{dt vacancy} is equivalent to
\begin{equation} \label{qd1}
\begin{split}
 &\chi(\lt^{(k-1)}\le n <\lt^{(k)}) -
 \chi(\lt^{(k)}\le n <\lt^{(k+1)}) \\
 =\qquad &
 \chi(\lbt^{(k-1)}\le n <\lbt^{(k)}) -
 \chi(\lbt^{(k)}\le n <\lbt^{(k+1)}).
\end{split}
\end{equation}
Observe that for $a\le b$,
$\chi(a \le n < b) = \chi(a \le n) - \chi(b \le n)$.
Consider the functions
\begin{equation*}
\begin{split}
  \Delta^{(k)}_n &= \chi(\lt^{(k)} \le n) - \chi(\lbt^{(k)} \le n) \\
  \bm^{(k)}_n &= \chi(m_n(\nu^{(k+1)}) \ge 1) \Delta^{(k)}_n \\
  \beq^{(k)}_n &= \chi(m_n(\nu^{(k)}) \ge 1) \Delta^{(k)}_n \\
  \bp^{(k)}_n &= \chi(m_n(\nu^{(k-1)}) \ge 1) \Delta^{(k)}_n.
\end{split}
\end{equation*}
For parts $n$ that occur in $\nu^{(k)}$, \eqref{qd1} is equivalent to
\begin{equation}\label{qd2}
  \bm^{(k-1)}_n - \beq^{(k)}_n = \beq^{(k)}_n - \bp^{(k+1)}_n.
\end{equation}
It will be shown that 
$\bm^{(k)}_n=\beq^{(k)}_n=\bp^{(k)}_n=0$ if $n$ is an unselected string
in $\nu^{(k+1)}$, $\nu^{(k)}$, $\nu^{(k-1)}$, respectively, which implies
\eqref{qd2}.

\subsection*{Unselected strings, generic or Case 2}
In this case $\lbt^{(k)}=\lt^{(k)}$, so that
$\Delta^{(k)}_n=0$, and $\bm^{(k)}_n = \beq^{(k)}_n=\bp^{(k)}_n=0$.

For the rest of the proof write $\ell=\lt^{(k)}=\lb^{(k)}$.

\subsection*{Unselected strings, Case 1}
Here $\lbt^{(k)}=\ell-1$, so that
$\Delta^{(k)}_n = -\chi(n=\ell-1)$.  Moreover
$m_{\ell-1}(\nu^{(k)})=0$ and $m_{\ell-1}(\nu^{(k+1)})=0$.  We have
\begin{equation*}
\begin{split}
  \beq^{(k)}_n&= -\chi(n=\ell-1)\chi(m_{\ell-1}(\nu^{(k)})\ge 1) = 0 \\
  \bm^{(k)}_n &= -\chi(n=\ell-1)\chi(m_{\ell-1}(\nu^{(k+1)})\ge1) = 0 \\
  \bp^{(k)}_n &= -\chi(n=\ell-1)\chi(m_{\ell-1}(\nu^{(k-1)}) \ge 1).
\end{split}
\end{equation*}
It must be shown that $\bp^{(k)}_n = 0$ if $n$ is an unselected string
in $\nu^{(k-1)}$ or equivalently that $m_{\ell-1}(\nu^{(k-1)})=0$
if $\ell-1$ is an unselected string in $\nu^{(k-1)}$.
Again we follow the division into cases of the proof of
Lemma \ref{db and dt}.

If $P^{(k)}_{\ell-1}(\nu)=0$, it also follows by
\eqref{vacancy ddiff} with $n=\ell-1$ that $m_{\ell-1}(\nu^{(k-1)})=0$.

Otherwise $P^{(k)}_{\ell-1}(\nu) = 1$.  We have
$m_{\ell-1}(\nu^{(k)}) = 0$, $\lt^{(k-1)} \le \ell-1$
and $\lb^{(k-1)} \le \ell-1$.  Also by \eqref{NPinq}
\begin{equation} \label{NNPinq}
  P^{(k)}_{\ell-2}(\nu) \le 2 - m_{\ell-1}(\nu^{(k-1)}).
\end{equation}
Consider the three cases $m_{\ell-1}(\nu^{(k-1)})\in \{0,1,2\}$.

If $m_{\ell-1}(\nu^{(k-1)})=0$ there is nothing to show.

Suppose $m_{\ell-1}(\nu^{(k-1)}) = 1$.
It must be shown that there are no unselected strings
of length $\ell-1$ in $\nu^{(k-1)}$.  Thus it suffices to show
that either $\lt^{(k-1)}=\ell-1$ or $\lb^{(k-1)}=\ell-1$, for then
the unique string in $\nu^{(k-1)}$ of length $\ell-1$ is selected.
So assume that $\lt^{(k-1)}\le \ell-2$ and $\lb^{(k-1)} \le \ell-2$.
By \eqref{NNPinq} $P^{(k)}_{\ell-2}(\nu) \le 1$.
Let $p$ be maximal such that $p<\ell-1$ and $m_p(\nu^{(k)}) \ge 1$;
if it does not exist set $p=0$.
If $p=\ell-2$ it follows from $P_{\ell-2}^{(k)}(\nu)\le 1$, 
$\lt^{(k-1)}\le \ell-2$ and $\lb^{(k-1)}\le \ell-2$ that either
$\dt$ or $\db$ selects a string of length $\ell-2$ in $\nu^{(k)}$
which contradicts the minimality of $\lt^{(k)}$ and $\lb^{(k)}$.
Hence assume that $p<\ell-2$.
By Lemma \ref{lower bound}, $P^{(k)}_n(\nu) = 1$ for
$p+1 \le n \le \ell-1$ and $P^{(k)}_p(\nu) \le 1$.
Moreover by \eqref{vacancy ddiff} $m_n(\nu^{(k-1)})=m_n(\nu^{(k+1)})=0$
for $p+2 \le n \le \ell-2$ and
\begin{equation} \label{Peq}
P^{(k)}_p(\nu) + m_{p+1}(\nu^{(k-1)}) + m_{p+1}(\nu^{(k+1)}) \le 1.
\end{equation}

Suppose $P^{(k)}_p(\nu)=1$.  This case can only occur if $p \ge 1$.
By \eqref{Peq} $m_{p+1}(\nu^{(k-1)})=0$,
so that $\lt^{(k-1)} \le p$ and $\lb^{(k-1)} \le p$.  But
$m_p(\nu^{(k)}) \ge 1$, so there is either a string in $\nu^{(k)}$
of length $p<\ell$ that is singular or of label $0$, contradicting
the minimality of $\lt^{(k)}$ and $\lb^{(k)}$.

Therefore $P^{(k)}_p(\nu) = 0$.
If $\lb^{(k-1)} \le p$ (which can only happen if $p\ge 1$),
since $m_p(\nu^{(k)}) \ge 1$ there is a singular string $(p,0)$
in $(\nu,J)^{(k)}$, contradicting the definition of $\lb^{(k)}=\ell$.
So $\lb^{(k-1)} > p$ and similarly $\lt^{(k-1)} > p$.
Since it is assumed that $\lb^{(k-1)}\le \ell-2$ and
$\lt^{(k-1)} \le \ell-2$, and there are no strings in
$(\nu,J)^{(k-1)}$ of lengths $n$ such that $p+2 \le n \le \ell-2$,
it follows that 
$\lb^{(k-1)}=\lt^{(k-1)}=p+1$ and $m_{p+1}(\nu^{(k-1)})=1$.
	
This rules out the generic case and Case 2 for $\nu^{(k-1)}$.
Case 1 for $\nu^{(k-1)}$ is also ruled out, for otherwise
if $p\ge 1$ then by induction $m_p(\nu^{(k)})=0$ which is a contradiction,
and if $p=0$ then $\lbt^{(k-1)} = \lb^{(k-1)} - 1 = p = 0$
which is also a contradiction.

Suppose Case 3 holds for $\nu^{(k-1)}$, that is,
$\ltb^{(k-1)}=\lbt^{(k-1)} > \lt^{(k-1)}=\lb^{(k-1)} = p+1$.
But $\lb^{(k-1)}=p+1$ and $m_n(\nu^{(k-1)})=0$ for $p+2 \le n \le \ell-2$,
so $m_n(\nub^{(k-1)})=0$ for $p+2\le n\le\ell-2$.
This implies $\lbt^{(k-1)} \ge \ell-1$.
Also $\ell-1=\lbt^{(k)}\ge \lbt^{(k-1)}$ so
$\lbt^{(k-1)}=\ell-1$.   Thus the unique string in
$(\nu,J)^{(k-1)}$ of length $\ell-1$ is selected while passing
from $\db(\nu,J)$ to $\dt(\db(\nu,J))$.

Suppose $m_{\ell-1}(\nu^{(k-1)})=2$.  By \eqref{NNPinq},
$P^{(k)}_{\ell-2}(\nu) = 0$.  If $m_{\ell-2}(\nu^{(k)}) = 0$ then
by \eqref{vacancy ddiff} it follows that
$P^{(k)}_{\ell-1}(\nu)=0$, which is a contradiction.
Thus $m_{\ell-2}(\nu^{(k)})\ge 1$.
By the minimality of $\lb^{(k)}=\lt^{(k)}=\ell$,
it follows that $\lt^{(k-1)} \ge \ell-1$
and $\lb^{(k-1)} \ge \ell-1$.  Since it is assumed that
$\lt^{(k-1)} \le \ell-1$ and $\lb^{(k-1)} \le \ell-1$,
we have $\lb^{(k-1)}=\lt^{(k-1)}=\ell-1$.
It must be shown that both of the strings
of length $\ell-1$ in $(\nu,J)^{(k-1)}$ get selected.
Since $m_{\ell-1}(\nu^{(k-1)})=2$,
Case 3 does not hold for $\nu^{(k-1)}$.
Case 1 does not hold either, for otherwise
$m_{\ell-2}(\nu^{(k)}) = 0$, which is a contradiction.
Therefore either the generic case or Case 2 holds
for $\nu^{(k-1)}$ and either way,
both strings in $(\nu,J)^{(k-1)}$ of length $\ell-1$ are selected.

\subsection*{Unselected strings, Case 3}
Here $\Delta^{(k)}_n = \chi(\ell \le n < \lbt^{(k)})$.
It follows from the proof of Lemma~\ref{db and dt} that
$m_n(\nu^{(j)}) = 0$ for $k-1 \le j \le k+1$ and
$\ell < n < \lbt^{(k)}$.  Hence
$\beq^{(k)}_n = \chi(n=\ell)\chi(m_\ell(\nu^{(k)})\ge 1)$,
$\bm^{(k)}_n = \chi(n=\ell)\chi(m_\ell(\nu^{(k+1)}) \ge 1)$, and
$\bp^{(k)}_n = \chi(n=\ell)\chi(m_\ell(\nu^{(k-1)}) \ge 1)$.
However in Case 3, $m_{\ell}(\nu^{(k)})=1$, so that
there are no other strings in $(\nu,J)^{(k)}$ of length $\ell$
other than the one selected by $\db$ and $\dt$.
So $\beq^{(k)}_n = 0$ for all applicable strings $(n,J)$.
By \eqref{cond3} also $\bm_n^{(k)}=0$.

It is enough to show that $\bp^{(k)}_n=0$, that is,
$m_\ell(\nu^{(k-1)})=0$ if $\ell$ is an unselected string.
By \eqref{cond2} we distinguish the cases 
$m_{\ell}(\nu^{(k-1)})\in\{0,1,2\}$.
If $m_\ell(\nu^{(k-1)})=0$ we are done.

Suppose that $m_\ell(\nu^{(k-1)})=1$.
By definition $\lb^{(k-1)}\le \ell$ and $\lt^{(k-1)}\le \ell$.
If $\lt^{(k-1)}=\ell$ or $\lb^{(k-1)}=\ell$ then the
only string of length $\ell$ in $(\nu,J)^{(k-1)}$ is selected.
So suppose $\lb^{(k-1)} \le \ell - 1$.  Then by~\eqref{cond2} 
and~\eqref{v3} $P^{(k)}_{\ell-1}(\nub)=0$.
Since $m_{\ell-1}(\nub^{(k)})\ge 1$ and $\lbt^{(k)}>\ell$,
it follows that $\lbt^{(k-1)} \ge \ell$.
On the other hand $\lbt^{(k-1)} \le \lt^{(k)}=\ell$ by
\eqref{ineq for lb},  so $\lbt^{(k-1)}=\ell$.
Thus the lone string in $\nu^{(k-1)}$ of length $\ell$
is selected in passing from $\nub^{(k-1)}$ to
$\nubt^{(k-1)}$.

Finally suppose $m_\ell(\nu^{(k-1)})=2$.  By \eqref{cond2}
$P^{(k)}_{\ell-1}(\nu)=0$.
It follows from \eqref{v3} that $P^{(k)}_{\ell-1}(\nub)=0$
and $\lb^{(k-1)} \ge \ell$.  Since $\lb^{(k-1)}\le\lb^{(k)}=\ell$,
we have $\lb^{(k-1)}=\ell$.
Similarly $\lt^{(k-1)} = \ell$.  It is enough to show that
either the generic case or Case 2 holds for $\nu^{(k-1)}$
for then both strings of length $\ell$ in $\nu^{(k-1)}$ are selected.
By \eqref{ineq for lb} $\lt^{(k-1)}=\ell=\lt^{(k)} \ge \lbt^{(k-1)}$
so Case 3 for $\nu^{(k-1)}$ is ruled out.
Suppose Case 1 holds for $\nu^{(k-1)}$.  Then
$\lb^{(k-1)}=\lt^{(k-1)}=\ell$ and
$\lbt^{(k-1)} = \lt^{(k-1)} - 1 = \ell - 1$.
By \eqref{db vacancy} with $n=\ell-1$,
$P^{(k)}_{\ell-1}(\nu) = P^{(k)}_{\ell-1}(\nub) = 0$.
Now $m_{\ell-1}(\nub^{(k)}) \ge 1$, so there is a string
$(\ell-1,0)\in (\nub,\Jb)^{(k)}$.  But $\ell-1 = \lbt^{(k-1)}$
so this is a contradiction to the definition of $\lbt^{(k)} > \ell$.
\end{proof}

\section{Proof of Lemma \ref{lem transpose}} \label{App B}

Lemma \ref{lem transpose} is established here by Corollary
\ref{trans config} and Lemma \ref{trans rigging}.

Let $R=(R_1,\dots,R_L)$ such that $R_j=(\eta_j^{\mu_j})$ with $\eta_L=1$, 
$\la$ a partition and $(\nu,J)\in\RC(\la^t;R^t)$.
Denote $(\nub,\Jb)=\jb(\nu,J)\in\RC(\rho^t;\Rb^t)$ and 
$(\nu^t,J^t)=\RCtr(\nu,J)\in\RC(\la;R)$.
Let $\lb^{(k)}$ be the lengths of the strings selected by $\db$
and $r =\rkb(\nu,J)$, the minimum index such that $\lb^{(r)}=\infty$.
By Proposition~\ref{j bar} the partition $\rho$ is obtained by removing 
the corner cell $(\la^t_r,r)$ of $\la$.

Recall that to every configuration $\nu$ there is an associated 
matrix \eqref{m def}.
Let $\db m$ be the matrix of $\nub$.  Then for all $i,j \ge 1$,
\begin{equation*}
  m_{ij} - (\db m)_{ij} =
  \chi(i-1 \ge 1) \chi(j = \lb^{(i-1)}) - \chi(j=\lb^{(i)}).
\end{equation*}

The matrix of the RC-transpose of $\db m$ can be calculated as follows.
\begin{equation*}
\begin{split}
((\db m)^t)_{ij} &= - (\db m)_{ji} + \chi((i,j)\in \rho)
	-\sum_{a=1}^L \chi((i,j)\in\Rb_a) \\
&= - m_{ji} + \chi(j-1 \ge 1) \chi(i=\lb^{(j-1)})-\chi(i=\lb^{(j)})
+ \chi((i,j)\in\la) \\
&\,\,-\chi((i,j)=(\la^t_r,r))-\sum_{a=1}^L
\chi((i,j) \in R_a)+\chi((i,j)=(\mu_L,1)) \\
&= \mt_{ij} + \chi(i=\lb^{(j-1)})-\chi(i=\lb^{(j)})-
	\chi((i,j)=(\la^t_r,r)).
\end{split}
\end{equation*}
Let $\del \mt$ be the matrix of the configuration of
$\del(\nu^t,J^t)$ and $\ell^{(k)}$ the
lengths of the strings selected by $\del$.  Then
\begin{equation*}
  (\del \mt)_{ij} = \mt_{ij} - \chi(i-1 \ge \mu_L)\chi(j=\ell^{(i-1)})+
  	\chi(i\ge \mu_L) \chi(j=\ell^{(i)}).
\end{equation*}
Thus the configurations of $\del \circ \RCtr(\nu,J)$ and
$\RCtr \circ \db(\nu,J)$ coincide, if and only if
\begin{equation} \label{tr eq}
\begin{split}
  &\,\,- \chi(i-1 \ge \mu_L)\chi(j=\ell^{(i-1)})+
  	\chi(i\ge \mu_L) \chi(j=\ell^{(i)}) \\
  &= \chi(i=\lb^{(j-1)})-\chi(i=\lb^{(j)})-\chi((i,j)=(\la^t_r,r)).
\end{split}
\end{equation}

\begin{rem} \label{exclusion}
By the definition of the RC-transpose map, there cannot be
a singular string of length $n$ in $\nu^{(k)}$
and a singular string of length $k$ in $\nu^{t(n)}$.
\end{rem}

\begin{rem} \label{P equals m} Let $i,j\ge 1$.
Suppose $(i,j+1)\in \la$, that is, $\la_{j+1}^t \ge i$.
Then by \cite[(9.7)]{KS} for the admissible $(\la^t;R^t)$-configuration
$\nu$, we have $m_j(\nu^{t(i)}) = P^{(j)}_i(\nu)$.

Similarly if $(i+1,j)\in\la$ then $m_i(\nu^{(j)}) = P^{(i)}_j(\nu^t)$.
\end{rem}

\begin{lem} \label{precise limit} Suppose that $\nu$ is an admissible
configuration of type $(\la^t;R^t)$.  Then
for all $n\ge \la^t_k$, $P^{(k)}_n(\nu) = \la^t_k - \la^t_{k+1}$.
\end{lem}
\begin{proof} By~\cite[(11.2)]{KS}, we have
\begin{equation} \label{big vac}
  P^{(k)}_n(\nu) = \la^t_k - \la^t_{k+1} +
  	\sum_{a=1}^L \min(0,n-\mu_a) \delta_{\eta_a,k}
\end{equation}
whenever $n \ge \nu^{(k-1)}_1$, $n\ge\nu^{(k)}_1$, and
$n\ge\nu^{(k+1)}_1$.  But in the proof of \cite[Lemma 10]{KS}
it is shown that $\la^t_{p+1} \ge \nu^{(p)}_1$
for all $p$, and in particular for $k-1 \le p \le k+1$.
Thus for $n\ge\la^t_k$, \eqref{big vac} holds.
It suffices to show that the sum over $a$ is zero.
Suppose not.  Then there is a rectangle $R_a$ in $R$ whose width
is $\eta_a=k$ and its height $\mu_a$ satisfies $\mu_a > n \ge \la^t_k$.
This means that $R_a$ is not contained in $\la$, which contradicts
the assumption that $\RC(\la^t;R^t)$ is nonempty.
\end{proof}

Define the sets of pairs of positive integers
\begin{equation*}
\begin{split}
H_\delta &= \{(i,j) \mid i = \lb^{(j)} \} \\
H_\del &= \{(i,j) \mid j = \ell^{(i)} \}.
\end{split}
\end{equation*}

Let $H = H_\delta \cup H_\del$ and $H^+ = H \cup \{(\la^t_r,r)\}$.

\begin{rem} \label{path subset}
$H\subset \la$.  To see this, by \cite[Lemma 10]{KS},
$m_n(\nu^{(k)}) > 0$ implies that $n \le \la^t_{k+1}$,
that is, $(n,k+1)\in \la$ and hence $(n,k)\in \la$.
Likewise, $m_k(\nu^{t(n)}) > 0$ implies $k \le \la_{n+1}$,
that is, $(n+1,k)\in \la$ and hence $(n,k)\in \la$.
This also shows that every cell of $H_\delta$ is not
at the end of its row in $\la$ and every cell of $H_\del$
is not at the bottom of its column in $\la$.
\end{rem}

\begin{lem} \label{path}
$(\la^t_r,r)\not\in H$ and
$H^+$ is a rookwise connected subset of the Ferrers diagram of $\la$ 
that can be viewed as leading from the cell $(1,1)$ to the
cell $(\la^t_r,r)$ where cells of $H_\delta$ are viewed as
commands to proceed to the east by one column and
cells of $H_\del$ are commands to proceed south by a row.
\end{lem}
\begin{proof} Define the cell $s_i=(i,1)$ for
$1\le i \le \mu_L-1$; $s_i \in H_\del$ by definition.
Also $s_i \in \la$ since $\RC(\la^t;R^t)$ is assumed to be nonempty
which by the Littlewood-Richardson rule implies that $\la$ contains
$R_L=(1^{\mu_L})$.  By induction suppose that the cell
$s_{m-1}\in H \cap \la$ has been defined for $m\ge \mu_L$.
Let $s_m$ be the cell just east of $s_{m-1}$ if $s_{m-1} \in H_\delta$
and the cell just south of $s_{m-1}$ if $s_{m-1}\in H_\del$.
This is well-defined since $H_\delta \cap H_\del=\emptyset$ by 
Remark \ref{exclusion}.
Moreover $s_m\in \la$ by Remark \ref{path subset}.
If $\mu_L>1$ the algorithm defines $s_{\mu_L}=(\mu_L,1)$.
If $\mu_L=1$ we define $s_1=(1,1)$.

We shall show that one of three possibilities must occur,
which are mutually disjoint by
Remarks \ref{exclusion} and \ref{path subset}.
\begin{enumerate}
\item $s_m$ is a corner cell of $\la$.
\item $s_m \in H_\delta$.
\item $s_m \in H_\del$.
\end{enumerate}
Write $s_m=(i,j)$.

Suppose that there is a singular string in $(\nu,J)^{(j)}$
of length $i$.  Then it will be shown that
$i=\lb^{(j)}$ and $(i,j)\in H_\delta$.
By construction either $s_{m-1}=(i,j-1)$ or
$s_{m-1}=(i-1,j)$.  Suppose that $s_{m-1}=(i,j-1)$, that is,
$i=\lb^{(j-1)}$.
Then by definition $i=\lb^{(j)}$ and $(i,j)\in H_\delta$.
Otherwise $s_{m-1}=(i-1,j)$.  By definition
$(\lb^{(j-1)},j-1)\in H_\delta$ and by the construction of $s_1$ through
$s_m$ we have $(i',j)\in H_\del$ for $\lb^{(j-1)}\le i'\le i-1$,
for otherwise the path would have proceeded into column $j+1$.
In particular for such $i'$ there are no singular strings of length
$i'$ in $(\nu,J)^{(j)}$ by Remark \ref{exclusion}, so that
by definition $i=\lb^{(j)}$.

Similarly if there is a singular string in $(\nu^t,J^t)^{(i)}$
of length $j$ then $j=\ell^{(i)}$ and $(i,j)\in H_\del$.

Otherwise it may be assumed that there is not a singular string
in $(\nu,J)^{(j)}$ of length $i$ and there is not a singular string in
$(\nu^t,J^t)^{(i)}$ of length $j$.

Under these assumptions, the following situations lead to
contradictions:
\begin{align}
\label{imposs1} m_i(\nu^{(j)})&=P^{(j)}_i(\nu)=0 \\
\label{imposs2} m_j(\nu^{t(i)})&=P^{(i)}_j(\nu^t)=0.
\end{align}
For suppose \eqref{imposs1} holds.
In the base case $m=\mu_L$, we have $(i,j)=(\mu_L,1)$.
Recall that $R_L=(1^{\mu_L})$.  This implies that
$m_{\mu_L}(\xi^{(1)}(R))>0$.  In particular
it follows immediately from \eqref{vacancy ddiff}
with $k=j=1$ and $n=i=\mu_L$ that \eqref{imposs1} is impossible.

Otherwise assume that we are in the
inductive case $m>\mu_L$.
Applying \eqref{vacancy ddiff} with $n=i$ and $k=j$, it follows that
$m_i(\nu^{(j-1)}) = 0$ and $P^{(j)}_{i-1}(\nu) = 0$.
The former immediately implies $(i,j-1)\not\in H$.
By construction, since $(i,j-1)\not\in H$,
it must be the case that $s_{m-1}=(i-1,j)\in H_\del$.
Thus there is a singular string of length $j$ in
$(\nu^t,J^t)^{(i-1)}$, so that $m_j(\nu^{t(i-1)})>0$.
By \eqref{Pm}, $m_{i-1}(\nu^{(j)})=P^{(i-1)}_j(\nu^t)$.
Suppose this common value is nonzero.
By \eqref{Pm}, $P^{(j)}_{i-1}(\nu) = m_j(\nu^{t(i-1)}) > 0$,
which is a contradiction.  Therefore $m_{i-1}(\nu^{(j)})=0$.
But also $P^{(j)}_{i-1}(\nu)=0$, so that
\eqref{imposs1} holds for $(i-1,j)$.  But this is a contradiction
by induction since $(i-1,j)$ is earlier in the path than $(i,j)$.

Similarly \eqref{imposs2} leads to a contradiction.

Now the proof divides into four cases depending on whether the
cells $(i+1,j)$ and $(i,j+1)$ are in $\la$ or not.

Suppose first that both $(i,j+1)\in\la$ and
$(i+1,j)\in\la$.  By Remark \ref{P equals m}
we have $a:=m_i(\nu^{(j)})=P^{(i)}_j(\nu^t)$
and $b:=m_j(\nu^{t(i)})=P^{(j)}_i(\nu)$.
If $a>0$ and $b>0$ then by the definition of the RC-transpose
map on riggings, there is either a singular string of length
$i$ in $(\nu,J)^{(j)}$ or a singular string of length $j$
in $(\nu^t,J^t)^{(i)}$, contradicting a previous assumption.
If $a>0$ and $b=0$ then there is a string of length
$i$ in $(\nu,J)^{(j)}$ with vacancy number zero, so it is
singular, which is a contradiction.
A similar contradiction is reached if $a=0$ and $b>0$.
If $a=b=0$ then one obtains the impossibility \eqref{imposs1}.

Next suppose $(i,j+1)\in\la$ and $(i+1,j)\not\in\la$ so that
$\la^t_j=i=\la^t_{j+1}$.
By Lemma \ref{precise limit} with $k=j$ and $n=i=\la^t_j$ we have
$P^{(j)}_i(\nu) = 0$.  If $m_i(\nu^{(j)})>0$ then
there is a singular string of length $i$ in $(\nu,J)^{(j)}$
which is a contradiction, and if $m_i(\nu^{(j)})=0$ one obtains
the impossibility \eqref{imposs1}.

In a similar manner one rules out the case
$(i,j+1)\not\in\la$ and $(i+1,j)\in\la$.

The last remaining case is $(i,j+1)\not\in\la$ and $(i+1,j)\not\in\la$.
But $(i,j)\in\la$ so $(i,j)$ is a corner cell of $\la$.
This finishes the proof of the trichotomy.

Now it is enough to show that if $s_m=(i,j)$ is a corner cell of
$\la$ then $j=r=\rkb(\nu,J)$.  There cannot be a cell of $H_\delta$ in
the $j$-th column or the path would have crossed into the $(j+1)$-st
column.  Thus $j \ge r$.
If $j>1$ then by construction there is a cell of
$H_\delta$ in the $(j-1)$-st column, so 
$j\le r$ and therefore $j=r$.
If $j=1$ then $H_\delta$ is empty and $\lb^{(1)}=\infty$,
and again $j=r$ since both equal $1$.
\end{proof}

\begin{cor} \label{trans config}
The configurations of $\RCtr\circ\db(\nu,J)$ and
$\del\circ\RCtr(\nu,J)$ coincide.
\end{cor}
\begin{proof} Follows directly from \eqref{tr eq}
and Lemma \ref{path}.
\end{proof}

\begin{lem} \label{trans rigging}
The riggings of $\RCtr\circ\db(\nu,J)$ and
$\del\circ\RCtr(\nu,J)$ coincide.
\end{lem}
\begin{proof}
Let $(\nub,\Jb)=\jb(\nu,J)$, $(\nu^t,J^t)=\RCtr(\nu,J)$,
$(\nub^t,\Jb^t)=\RCtr\circ \jb(\nu,J)$ and
$(\nutrb,\Jtrb)=\del\circ \RCtr(\nu,J)$. 
Denote the corresponding partitions by $\Jb_n^{(k)}$, $J^{t(k)}_n$,
$\Jb^{t(k)}_n$ and $\Jtrb_n^{(k)}$, respectively.
The aim is to show that $\Jb_k^{t(n)}=\Jtrb_k^{(n)}$ for all $n,k\ge 1$.
Since the transpose of a rigging depends on the vacancy number and
the vacancy number changes under $\jb$ according to \eqref{db vacancy} 
we need to distinguish several ranges for $n$.

First suppose that $\lb^{(k-1)}\le n<\lb^{(k)}-1$. 
In this case the riggings are not changed by $\jb$
so that $J_n^{(k)}=\Jb_n^{(k)}$ and furthermore the largest part of 
$J_n^{(k)}$ is smaller than $P_n^{(k)}(\nu)$ since by the definition
of $\lb^{(k)}$ there are no singular strings in the range of $n$.
The rectangle corresponding to $\Jb_n^{(k)}$ has height
$m_n(\nub^{(k)})=m_n(\nu^{(k)})$ and width 
$P_n^{(k)}(\nub)=P_n^{(k)}(\nu)-1$ by \eqref{db vacancy}.
Hence, compared to $J_k^{t(n)}$, $\Jb_k^{t(n)}$ has one part
of length $m_n(\nu^{(k)})$ less. 
It follows from Lemma \ref{path} that $\ell^{(n)}=k$. 
Since furthermore $\mu_L\le \lb^{(k-1)}\le n$, this implies that 
$\del$ removes one part of length $P_k^{(n)}(\nu^t)$ from $J_k^{t(n)}$.
It also ensures that $m_k(\nu^{t(n)})>0$.
By \eqref{Pm}, $P_k^{(n)}(\nu^t)=m_n(\nu^{(k)})$ which proves
$\Jb_k^{t(n)}=\Jtrb_k^{(n)}$.

Next suppose that $\lb^{(k)}<n<\lb^{(k+1)}$. 
Again by the choice of $n$ we have $J_n^{(k)}=\Jb_n^{(k)}$,
and the rectangle containing $\Jb_n^{(k)}$ has height $m_n(\nub^{(k)})=
m_n(\nu^{(k)})$ and width $P_n^{(k)}(\nub)=P_n^{(k)}(\nu)+1$ 
by \eqref{db vacancy}. Hence, compared to $J_k^{t(n)}$,
$\Jb_k^{t(n)}$ has an extra part of size $m_n(\nu^{(k)})$.
By Lemma \ref{path}, $k=\ell^{(n)}-1$ so that $\del$ adds an extra part 
of size $P_k^{(n)}(\nutrb)$ to $J_k^{t(n)}$ and $m_k(\nutrb^{(n)})>0$. 
By Corollary \ref{trans config} this also implies that $m_k(\nub^{t(n)})>0$.
Hence one finds by \eqref{Pm} that
$P_k^{(n)}(\nutrb)=P_k^{(n)}(\nub^t)=m_n(\nub^{(k)})=m_n(\nu^{(k)})$.
This shows that $\Jb_k^{t(n)}=\Jtrb_k^{(n)}$.

Now assume that $n=\lb^{(k)}$ and $\lb^{(k)}<\lb^{(k+1)}$.
Since $n=\lb^{(k)}$ the largest part of
$J_n^{(k)}$ equals $P_n^{(k)}(\nu)$ which is removed under $\jb$.
The rectangle corresponding to $\Jb_n^{(k)}$ has height $m_n(\nub^{(k)})
=m_n(\nu^{(k)})-1$ and width $P_n^{(k)}(\nub)=P_n^{(k)}(\nu)+1$.
Hence, compared to $J_k^{t(n)}$, $\Jb_k^{t(n)}$ has an extra part of
size $m_n(\nu^{(k)})-1$. From Lemma \ref{path} we find that $\ell^{(n)}=k+1$
so that $\del$ adds a part of
size $P_k^{(n)}(\nutrb)$ to $J_k^{t(n)}$.  Since $\del$ added a part,
$m_k(\nub^{t(n)})=m_k(\nutrb^{(n)})>0$ which by
\eqref{Pm} implies that $P_k^{(n)}(\nutrb)=m_n(\nu^{(k)})-1$.  Hence
$P_k^{(n)}(\nutrb)=P_k^{(n)}(\nub^t)=m_n(\nub^{(k)})=m_n(\nu^{(k)})-1$
by \eqref{Pm}.

Next assume that $n=\lb^{(k)}=\lb^{(k+1)}$.
Again $\Jb_n^{(k)}$ is obtained from $J_n^{(k)}$ by removing
a part of size $P_n^{(k)}(\nu)$, and this time the rectangle
corresponding to $\Jb_n^{(k)}$ has height $m_n(\nub^{(k)})=m_n(\nu^{(k)})-1$
and width $P_n^{(k)}(\nub)=P_n^{(k)}(\nu)$. Therefore
$\Jb_k^{t(n)}=J_k^{t(n)}$.
By Lemma \ref{path} we have $\ell^{(n-1)}\le k<\ell^{(n)}-1$ so that
also $\Jtrb_k^{(n)}=J_k^{t(n)}$ which proves the assertion.

Consider the case $n=\lb^{(k)}-1$ and $\lb^{(k-1)}<\lb^{(k)}$.
By the assumptions on $n$ there is no singular string of length $n$ in 
$\nu^{(k)}$ so that the largest part of $J_n^{(k)}$ is smaller than
$P_n^{(k)}(\nu)$. Under $\jb$ a part of size $P_n^{(k)}(\nub)
=P_n^{(k)}(\nu)-1$ gets added to $J_n^{(k)}$, and the rectangle
corresponding to $\Jb_n^{(k)}$ has height $m_n(\nub^{(k)})=m_n(\nu^{(k)})+1$
and width $P_n^{(k)}(\nu)-1$. Hence, compared to $J_k^{t(n)}$, a part
of size $m_n(\nu^{(k)})$ is missing in $\Jb_k^{t(n)}$. 
By Lemma \ref{path} $\ell^{(n)}=k$.  Hence
$\del$ removes a part of size $P_k^{(n)}(\nu^t)$ from
$J_k^{t(n)}$.  But $P_k^{(n)}(\nu^t)=m_n(\nu^{(k)})$ by \eqref{Pm}
because $m_k(\nu^{t(n)})>0$ since $k=\ell^{(n)}$.

Now assume that $n+1=\lb^{(k)}=\lb^{(k-1)}$.
The partition $\Jb_n^{(k)}$ is obtained from $J_n^{(k)}$ by adding
a part of size $P_n^{(k)}(\nub)=P_n^{(k)}(\nu)$ and its rectangle has
height $m_n(\nub^{(k)})=m_n(\nu^{(k)})+1$ and width
$P_n^{(k)}(\nu)$. Hence $\Jb_k^{t(n)}=J_k^{t(n)}$.
By Lemma \ref{path}, $\ell^{(n)}<k<\ell^{(n+1)}$ if $k>1$ which implies
that $\del$ preserves the rigging so that $\Jtrb_k^{(n)}=J_k^{t(n)}$ 
proving the assertion. If $k=1$ then $n=\lb^{(0)}-1=\mu_L-1$. In this case 
$\del$ also does not change the rigging so that $\Jtrb_k^{(n)}=J_k^{t(n)}$.

Finally assume that none of the above cases holds for $n$ so that
$m_n(\nub^{(k)})=m_n(\nu^{(k)})$, $P_n^{(k)}(\nub)=P_n^{(k)}(\nu)$ and
$\Jb_n^{(k)}=J_n^{(k)}$ which implies that $\Jb_k^{t(n)}=J_k^{t(n)}$.
If $n\ge \lb^{(k+1)}$ one finds by Lemma \ref{path} that $k\le \ell^{(n)}-2$.
If $n<\lb^{(k-1)}$ for $k>1$, Lemma \ref{path} implies that $k>\ell^{(n)}$.
The remaining case is $n<\lb^{(0)}=\mu_L$. In all three cases
$\del$ does not alter $J_k^{t(n)}$ so that indeed $\Jtrb_k^{(n)}=\Jb_k^{t(n)}$.
\end{proof}

\end{document}